\documentclass[11pt]{article}

\usepackage{latexsym,amsfonts,amsmath,amssymb,graphicx,hyperref,verbatim,enumerate}

\newcommand{\SA}{\mbox{\bf S}}
\newcommand{\rA}{\mbox{\bf r}}
\newcommand{\RA}{\mbox{\bf R}}
\newcommand{\HA}{\mbox{\bf H}}

\newcommand{\St}{\mathbb{S}}
\newcommand{\EE}{\mbox{\bf E}\,}
\newcommand{\PP}{\mathbb{P}}

\newcommand{\R}{\mathbb{R}}
\newcommand{\C}{\mathbb{C}}
\newcommand{\Q}{\mathbb{Q}}
\newcommand{\HH}{\mathbb{H}}
\newcommand{\N}{\mathbb{N}}
\newcommand{\D}{\mathbb{D}}
\newcommand{\TT}{\mathbb{T}}
\newcommand{\Z}{\mathbb{Z}}
\newcommand{\A}{\mathbb{A}}

\newcommand{\pa}{\partial}

\newcommand{\F}{{\cal F}}

\newcommand{\no}{\noindent}

\newcommand{\BGE}{\begin{equation}}
\newcommand{\BGEN}{\begin{equation*}}
\newcommand{\EDE}{\end{equation}}
\newcommand{\EDEN}{\end{equation*}}
\def\conf{\stackrel{\rm Conf}{\twoheadrightarrow}}

\def\eps{\varepsilon}
\def\til{\widetilde}
\def\ha{\widehat}
\def\sem{\setminus}

\def\lin{\overline}

\DeclareMathOperator{\ccap}{cap} 
\DeclareMathOperator{\sign}{sign} 
\DeclareMathOperator{\dist}{dist} 
 \DeclareMathOperator{\id}{id}
\DeclareMathOperator{\Imm}{Im } \DeclareMathOperator{\Ree}{Re }

\DeclareMathOperator{\modd}{mod} \DeclareMathOperator{\JP}{JP}
\DeclareMathOperator{\PV}{P.V.} 
\DeclareMathOperator{\mA}{m} \DeclareMathOperator{\Res}{Res}

\def\h0{{\bf h}}

\newtheorem{Lemma}{Lemma}[section]
\newtheorem{Theorem}{Theorem}[section]
\newtheorem{Definition}{Definition}[section]
\newtheorem{Corollary}{Corollary}[section]
\newtheorem{Proposition}{Proposition}[section]

\numberwithin{equation}{section}

\begin{document}
\title{\bf Reversibility of Whole-Plane SLE}
\date{\today}
\author{Dapeng Zhan\footnote{Supported by NSF grant DMS-0963733} \\ Michigan State University}
\maketitle

\begin{abstract}
The main result of this paper is that, for $\kappa\in(0,4]$, whole-plane SLE$_\kappa$ satisfies reversibility, which means that the time-reversal of a whole-plane SLE$_\kappa$ trace is still a whole-plane SLE$_\kappa$ trace. In addition, we find that the time-reversal of a radial SLE$_\kappa$ trace for $\kappa\in(0,4]$ is a disc SLE$_\kappa$ trace with a marked boundary point. The main tool used in this paper is a stochastic coupling technique, which is used to couple two whole-plane SLE$_\kappa$ traces so that they overlap.  Another tool used is the Feynman-Kac formula, which is used to solve a PDE. The solution of this PDE is then used to construct the above coupling.
\end{abstract}


\section{Introduction}
The stochastic Loewner evolution (SLE) introduced by Oded Schramm (\cite{S-SLE}) describes some random fractal curves in plane
domains that satisfy conformal invariance and Domain Markov Property. These two properties make SLEs the most suitable candidates
for the scaling limits of many two-dimensional lattice models at criticality. These models are proved or conjectured to converge to SLE
with different parameters (e.g., \cite{LSW-8/3}\cite{LSW-2}\cite{SS-4}\cite{SS}\cite{SS-6}\cite{SS-3}).
For basics of SLE, the  reader may refer to  \cite{LawSLE} and \cite{RS-basic}.

There are several different versions of SLEs, among which  chordal SLE and radial SLE are the most well-known.
A chordal or radial SLE trace is a random fractal curve that grows in a simply connected plane domain from a boundary point.
A chordal SLE trace ends at another boundary point, and a radial SLE trace ends an interior point.
Their behaviors both depend on a positive parameter $\kappa$. When $\kappa\in(0,4]$, both traces are simple curves, and all points on the trace other than the initial and final points lie inside the domain. When $\kappa>4$, the traces have self-intersections.

A stochastic coupling technique was introduced in \cite{reversibility} to prove that, for $\kappa\in(0,4]$,
chordal SLE$_\kappa$ satisfies reversibility, which means that if $\beta$ is a chordal SLE$_\kappa$ trace in a domain $D$ from
$a$ to $b$, then after a time-change, the time-reversal of $\beta$ becomes a chordal SLE$_\kappa$ trace in $D$ from $b$ to $a$.
The  technique was later used (\cite{duality}\cite{duality2}) to prove Duplantier's duality conjecture, which says that,
for $\kappa>4$, the boundary of the hull generated by a chordal SLE$_\kappa$ trace looks locally like an SLE$_{16/\kappa}$ trace.
The technique was also used to prove that the radial or chordal SLE$_2$ can be obtained by erasing loops on a planar Brownian motion (\cite{LEBM}), and the chordal SLE$(\kappa,\rho)$ introduced in \cite{LSW-8/3} also satisfies reversibility for $\kappa\in(0,4]$ and $\rho\ge \kappa/2-2$ (\cite{kappa-rho}).

Since the initial point and final point of a radial SLE are topologically different, the time-reversal of a radial SLE  trace can not be a radial SLE trace. However, we may consider whole-plane SLE instead, which describes a random fractal curve in the Riemann sphere $\ha\C=\C\cup\{\infty\}$ that grows from one interior point to another interior point. Whole-plane SLE is related to radial SLE as follows: conditioned on the initial part of a whole-plane SLE$_\kappa$ trace, the rest part of such trace has the distribution of a radial SLE$_\kappa$ trace that grows in the complementary domain of the initial part of this trace. The main result of this paper is the following theorem.

\begin{Theorem} Whole-plane SLE$_\kappa$ satisfies reversibility for $\kappa\in(0,4]$. \label{Main-Thm}
\end{Theorem}

The theorem in the case $\kappa=2$  has been proved  in \cite{int-LERW}. The proof used the reversibility of loop-erased random walk (LERW for short, see \cite{Law}) and the convergence of LERW to whole-plane SLE$_2$. In this paper we will obtain a slightly more general result: the reversibility of whole-plane SLE$(\kappa,s)$ process, which is defined by adding a constant drift  to the driving function for the whole-plane SLE$_\kappa$ process. This is the statement of Theorem \ref{Main-Thm-skew}.

To get some idea of the proof, let's first review the proof of the reversibility of chordal SLE$_\kappa$ in \cite{reversibility}. We constructed a pair of chordal SLE$_\kappa$ traces $\gamma_1$ and $\gamma_2$ in a simply connected domain $D$, where $\gamma_1$ grows from a boundary point $a_1$ to another boundary point $a_2$,  $\gamma_2$ grows from $a_2$ to $a_1$, and these two traces commute in the following sense. Fix $j\ne k\in\{1,2\}$, if $T_k$ is a stopping time for $\gamma_k$, then conditioned on $\gamma_k(t)$, $t\le T_k$, the part of $\gamma_j$ before hitting $\gamma_k(t)((0,T_k])$ has the distribution of a chordal SLE$_\kappa$ trace that grows from $a_j$ to $\gamma_k(T_k)$ in $D_k(T_k)$, which is a component of $D\sem \gamma_k(t)((0,T_k])$. In the case $\kappa\le 4$, a.s.\ $\gamma_j$ hits $\gamma_k(t)((0,T_k])$ exactly at $\gamma_k(T_k)$, so $\gamma_j$ visits $\gamma_k(T_k)$ before any $\gamma_k(t)$, $t<T_k$. Since this holds for any stopping time $T_k$ for $\gamma_k$, the two traces a.s.\ overlap, which implies the reversibility.

To prove the reversibility of whole-plane SLE$_\kappa$, we want to construct two whole-plane SLE$_\kappa$ traces in $D=\ha\C$, one is $\gamma_1$ from $a_1$ to $a_2$, the other is $\gamma_2$ from $a_2$ to $a_1$, so that $\gamma_1$ and $\gamma_2$ commute. Here we can not expect that they commute in exactly the same sense as in the above paragraph. Note that conditioned on $\gamma_k(t)$, $t\le T_k$, the part of $\gamma_j$ before hitting $\gamma_k(t)$, $t\le T_k$, can not have the distribution of a whole-plane SLE$_\kappa$ trace in $D_k(T_k)$ from $a_1$ to $\gamma_k(T_k)$ because now the complementary domain $D_k(T_k)$ is topologically different from $\ha\C$, while whole-plane SLEs are only defined in $\ha\C$. Since the conditional curve grows from an interior point to a boundary point, it is neither a radial SLE trace nor a chordal SLE trace.

Thus, we need to define SLE traces in simply connected domains that grow from an interior point to a boundary point.
We use the idea of defining whole-plane SLE using radial SLE. The situation here is a little different: after a positive initial part, the rest part of the curve grows in a doubly connected domain. Another difference is that there is a marked point on the boundary of the initial domain.
In this paper, we use the annulus Loewner equation introduced in \cite{Zhan} together with an annulus drift function $\Lambda=\Lambda(t,x)$ to define the so-called annulus SLE$(\kappa,\Lambda)$ process in a doubly connected domain $D$, which starts from a point $a\in\pa D$, and whose growth is affected by a marked point $b\in\pa D$. In the case when $a$ and $b$ lie on different boundary components, by shrinking the boundary component containing $a$ to a singlet, we get the so called disc SLE$(\kappa,\Lambda)$, which describes a random curve that grows in a simply connected domain and starts from an interior point.

We find that if $\Lambda=\kappa\frac{\Gamma'}{\Gamma}$, where $\Gamma$ is a positive differentiable function defined on $(0,\infty)\times\R$ that solves a linear PDE and satisfies some periodic condition (see (\ref{dot-Gamma=}) and (\ref{Gamma-period-s})), then using the coupling technique we could construct a coupling of two whole-plane SLE$_\kappa$ traces: $\gamma_1$ and $\gamma_2$, which commute in the sense that, conditioned on one curve up to a finite stopping time $T$, the other curve is a disc SLE$(\kappa,\Lambda)$ trace in the remaining domain, and its marked point is the tip point of the first curve at $T$. 

The main new idea in the current paper is an application of a Feynman-Kac representation, which is used to get a formal solution of the PDE for $\Gamma$ in the case $\kappa\in(0,4]$. 
Using Fubini's Theorem, It\^o's formula, and some estimations, we prove that the formal solution $\Gamma_\kappa$ is smooth and solves the PDE. We then find that $\Lambda_\kappa:=\kappa\frac{\Gamma_\kappa'}{\Gamma_\kappa}$ satisfies the property that the marked point for an annulus or disc SLE$(\kappa,\Lambda_\kappa)$ process is a subsequential limit point of the trace. This property implies that, if two whole-plane SLE$_\kappa$ traces commute as in the previous paragraph, then they must overlap. So the main theorem is proved. Moreover, from the relation between whole-plane SLE$_\kappa$ and radial SLE$_\kappa$, we conclude that, for $\kappa\in(0,4]$, the time-reversal of a radial SLE$_\kappa$ trace is a disc SLE$(\kappa,\Lambda_\kappa)$ trace. 

The marked point and the initial point of an annulus SLE$(\kappa,\Lambda)$ process could also lie on the same boundary component. In this case, if $\Lambda=\kappa\frac{\Gamma'}{\Gamma}$, and $\Gamma$ satisfies a similar linear PDE (see (\ref{dot-Gamma=2})),
then for a doubly connected domain $D$ with two boundary points $a_1$ and $a_2$ on the same boundary component, we can construct a pair of annulus SLE$(\kappa,\Lambda)$ traces $\gamma_1$ and $\gamma_2$ in $D$, which commute with each other.
If an SLE process in a doubly connected domain is the scaling limit of some random path in a lattice, which satisfies reversibility at the discrete level, then such SLE should satisfy reversibility. We hope that the work in this paper will shed some light on the study of these processes.

The study on the commutation relations of SLE in doubly connected domains continues the work in \cite{Julien-Comm} by Dub\'edat, who used some tools from Lie Algebra to obtain commutation conditions of SLE in simply connected domains. 
The annulus SLE$(\kappa,\Lambda_{\kappa})$ process used to prove the reversibility of whole-plane SLE$_\kappa$ was later studied in \cite{ann-restriction}. When $\kappa=8/3$, the process satisfies the restriction property, which is similar to the restriction property for chordal SLE$_{8/3}$ (see \cite{LSW-8/3}). For $\kappa\in(0,4]\sem \{8/3\}$, it satisfies some ``weak'' restriction property. 

G.\ F.\ Lawler  (\cite{Law-restriction}) used a different method to define annulus SLE$_\kappa$ processes for $\kappa\in(0,4]$, which agree with our  annulus SLE$(\kappa,\Lambda_{\kappa})$ processes. His construction uses the Brownian loop measures. The ``strong'' ($\kappa=8/3$) and ``weak'' ($\kappa\ne 8/3$) restriction properties of Lawler's annulus SLE processes are immediate from the definition; and the reversibility of these processes follows from the chordal reversibility. However, the reversibility of whole-plane SLE is not proved in \cite{Law-restriction}. To get the whole-plane reversibility, some additional work is required based on Lawler's work. In this paper, the reversibility of annulus SLE$(\kappa,\Lambda_{\kappa})$ and the reversibility of whole-plane SLE$_\kappa$ are proved separately, and the coupling technique is applied in both proofs, which are similar though.

J.\ Miller and S.\ Sheffield recently proved the reversibility of whole-plane SLE (\cite{MS4}) for all $\kappa\in[0,8]$. Their proof uses the imaginary geometry of Gaussian free field developed in their earlier papers (c.f.\ \cite{MS1}).

This paper is organized as follows. In Section \ref{prelim}, we introduce some symbols and notations. In Section \ref{Loewner}, we review  several versions of Loewner equations. In Section \ref{section-reference}, we define annulus SLE$(\kappa,\Lambda)$ and disc SLE$(\kappa,\Lambda)$ processes, whose growth is affected by one marked boundary point. In Section \ref{coupling} we prove that when $\Gamma$ solves (\ref{dot-Gamma=}) or (\ref{dot-Gamma=2}), there is a commutation coupling of two annulus SLE$(\kappa,\Lambda)$ processes, where $\Lambda=\kappa\frac{\Gamma'}{\Gamma}$. In Section \ref{coupling-deg}, we  construct a coupling of two whole-plane SLE processes, which is similar to the coupling in the previous section. In Section \ref{Section-PDE}, we solve PDE (\ref{dot-Gamma=}) using a Feynman-Kac expression, and the solution is then used to prove the reversibility of whole-plane SLE$_\kappa$ process. In fact, we obtain a slightly more general result: the reversibility of skew whole-plane SLE$_\kappa$ processes for $\kappa\in(0,4]$.
In the last section, we find some solutions to the PDE for $\Gamma$ and $\Lambda$ when $\kappa\in\{0,2,3,4,16/3\}$, which can be expressed in terms of well-known special functions.

\vskip 3mm

\no {\bf Acknowledgements.} I would like to thank Alexander Volberg for his suggestions on transforming PDE (\ref{dot-Gamma-c=}) into (\ref{PDE-Psi*}), and thank Zhen-Qing Chen for his help on Lemma \ref{Lemma-inq}. I also thank Gregory Lawler for his valuable comments and suggestions.

\section{Preliminary} \label{prelim}
\subsection{Symbols}
Throughout this paper, we will use the following symbols. Let $\ha\C=\C\cup\{\infty\}$, $\D=\{z\in\C:|z|<1\}$, $\TT=\{z\in\C:|z|=1\}$, and $\HH=\{z\in\C:\Imm z>0\}$. For $p>0$, let $\A_p=\{z\in\C:1>|z|> e^{-p}\}$ and $\St_p=\{z\in\C:0< \Imm z<p\}$. For $p\in\R$, let $\TT_p=\{z\in\C:|z|=e^{-p}\}$ and $\R_p=\{z\in\C:\Imm z=p\}$. Then $\pa\D=\TT$, $\pa\HH=\R$, $\pa\A_p=\TT\cup\TT_p$, and $\pa \St_p=\R\cup\R_p$. Let $e^i$ denote the map $z\mapsto e^{iz}$. Then $e^i$ is a covering map from $\HH$ onto $\D$, and from $\St_p$ onto $\A_p$; and it maps $\R$ onto $\TT$ and maps $\R_p$ onto $\TT_p$. For a doubly connected domain $D$, we use $\modd(D)$ to denote its modulus. For example, $\modd(\A_p)=p$.

A conformal map in this paper is a univalent analytic function. A conjugate conformal map is defined to be the complex conjugate of a conformal map. Let $I_0(z)=1/\lin{z}$ be the reflection w.r.t.\ $\TT$. Then $I_0$ is a conjugate conformal map from $\ha\C $ onto itself, fixes $\TT$, and interchanges $0$ and $\infty$. Let $\til I_0(z)=\lin{z}$ be the reflection w.r.t.\ $\R$. Then $\til I_0$ is a conjugate conformal map from $\C$ onto itself and satisfies $e^i\circ \til I_0=I_0\circ e^i$. For $p>0$, let $I_p(z):=e^{-p}/\lin{z}$ and $\til I_p(z)=ip+\lin{z}$.
Then $I_p$ and $\til I_p$ are conjugate conformal automorphisms of $\A_p$ and $\St_p$, respectively.
Moreover, $I_p$ interchanges $\TT_p$ and $\TT$, $\til I_p$ interchanges $\R_p$ and $\R$, and $I_p\circ e^i
=e^i\circ \til I_p$.

We will frequently use functions $\cot(z/2)$, $\tan(z/2)$, $\coth(z/2)$, $\tanh(z/2)$, $\sin(z/2)$, $\cos(z/2)$, $\sinh(z/2)$, and $\cosh(z/2)$. For simplicity, we write $2$ as a subscript. For example, $\cot_2(z)$ means $\cot(z/2)$, and we have $\cot_2'(z)=-\frac 12\sin_2^{-2}(z)$.

An increasing function in this paper will always be strictly increasing.
For a real interval $J$, we use $C(J)$ to denote the space of real continuous functions on $J$.
The maximal solution to an ODE or SDE with initial value is the solution with the biggest definition domain.

Many functions in this paper depend on two variables. In some of these functions, the first variable represents time or modulus, and the second variable does not. In this case, we use $\pa_t$ and $\pa_t^{n}$ to denote the partial derivatives w.r.t.\ the first variable, and use $'$, $''$, and the superscripts $(h)$ to denote the partial derivatives w.r.t.\ the second variable. For these functions, we say that it has period $r$ (resp.\ is even or odd) if it has period $r$ (resp.\ is even or odd) in the second variable when the first variable is fixed. Some functions in Section \ref{coupling} and Section \ref{coupling-deg} depend on two variables: $t_1$ and $t_2$, which both represent time. In this case we use $\pa_j$ to denote the partial derivative w.r.t.\ the $j$-th variable, $j=1,2$.

In this paper, a domain is a connected open subset of $\ha\C$, and continuum is a connected compact subset of $\ha\C$ that contains more than one point. A continuum $K$ is called a hull in $\C$ if $K\subset\C$ and $\ha\C\sem K$ is connected. In this case, there is a unique conformal map $f_K$ from $\ha\C\sem\lin{\D}$ onto $\ha\C\sem K$ and satisfies $\lim_{z\to\infty} f_K(z)/z=a_K$ for some positive number $a_K$. Then $a_K$ is called the capacity of $K$, and is denoted by $\ccap(K)$.

A doubly connected domain in this paper is a domain whose complement is a disjoint union of two continuums. Let $D$ be a doubly connected domain. If $K$ is a relatively closed subset of $D$, has positive distance from one boundary component of $D$, and if $D\sem K$ is also doubly connected, then we call $K$ a hull in $D$, and call the number $\modd(D)-\modd(D\sem K)$ the capacity of $K$ in $D$, and let it be denoted by $\ccap_{D}(K)$.

\subsection{Brownian motions} \label{section-BM}
Throughout this paper, a Brownian motion means a standard one-dimensional Brownian motion, and $B(t)$, $0\le t<\infty$, will always be used to denote a Brownian motion. This means that $B(t)$ is continuous, $B(0)=0$, and $B(t)$ has independent increment with $B(t)-B(s)\sim N(0,t-s)$ for $t\ge s\ge 0$. For $\kappa\ge 0$, the rescaled Brownian motion $\sqrt\kappa B(t)$ will be used to define annulus SLE$_\kappa$. The symbols $B_*(t)$, $\ha B_*(t)$, or $\til B_*(t)$ will also be used to denote a Brownian motion, where the $*$ stands for subscript. Let $(\F_t)_{t\ge 0}$ be a filtration. By saying that $B(t)$ is an $(\F_t)$-Brownian motion, we mean that $(B(t))$ is $(\F_t)$-adapted, and for any fixed $t_0\ge 0$, $B(t_0+t)-B(t_0)$, $t\ge 0$, is a Brownian motion independent of $\F_{t_0}$.

\begin{Definition}
  Let $\kappa>0$ and $(\F_t)_{t\in\R}$ be a right-continuous filtration. A process $B^{(\kappa)}(t)$, $t\in\R$, is called a pre-$(\F_t)$-$(\TT;\kappa)$-Brownian motion if $(e^i(B^{(\kappa)}(t)))$ is $(\F_t)$-adapted, and for any $t_0\in\R$, \BGE B_{t_0}(t):=\frac{1}{\sqrt\kappa}\Big(B^{(\kappa)}(t_0+t)-B^{(\kappa)}(t_0)\Big),\quad 0\le t<\infty,\label{B-c}\EDE
  is an $(\F_{t_0+t})$-Brownian motion. If $(\F_t)$ is generated by $(e^i(B^{(\kappa)}(t)))$, then we simply call $(B^{(\kappa)}(t))$ a pre-$(\TT;\kappa)$-Brownian motion.
\end{Definition}

\no{\bf Remark.} The name of the  pre-$(\TT;\kappa)$-Brownian motion comes from the fact that $B_{\TT}(t):=e^i(B^{(\kappa)}(t))$, $t\in\R$, is a Brownian motion on $\TT$ with speed $\kappa$: for every $t_0\in\R$,  $B_{\TT}(t_0)$ is uniformly distributed on $\TT$; and $B_{\TT}(t_0+t)/B_{\TT}(t_0)$, $t\ge 0$, has the distribution of $e^i(\sqrt\kappa B(t))$, $t\ge 0$, and is independent of $B_{\TT} (t)$, $t\le t_0$. One may construct $B^{(\kappa)}(t)$ as follows. Let $B_+(t)$ and $B_-(t)$, $t\ge 0$, be two independent Brownian motions. Let ${\bf x}$ be a random variable uniformly distributed on $[0,2\pi)$, which is independent of $(B_\pm(t))$. Let $B^{(\kappa)}(t) ={\bf x}+\sqrt \kappa B_{\sign(t)}(|t|)$ for $t\in\R$. Then $B^{(\kappa)}(t)$, $t\in\R$, is a pre-$(\TT;\kappa)$-Brownian motion.

\begin{Definition}
Let $B^{(\kappa)}(t)$, $t\in\R$, be a pre-$(\F_t)$-$(\TT;\kappa)$-Brownian motion, where $(\F_t)$ is right-continuous, and every $\F_t$ contains all eligible events w.r.t.\ the process $(e^i(B^{(\kappa)}(t)))$. Suppose $T$ is an $(\F_t)$-stopping time, and $T>t_0$ for a deterministic number $t_0\in\R$. We say that $X(t)$ satisfies the $(\F_t)$-adapted SDE
$$ dX(t)=a(t) dB^{(\kappa)}(t)+b(t)dt, \quad -\infty<t<T,$$ 
if $e^i(X(t))$, $a(t)$, and $b(t)$ are continuous and $(\F_t)$-adapted, and if for any deterministic number $t_0$ with $t_0<T$, $X_{t_0}(t):=X(t_0+t)-X(t_0)$ satisfies the following $(\F_{t_0+t})_{t\ge 0}$-adapted SDE with the traditional meaning in \cite{RY}: $$ dX_{t_0}(t)=a_{t_0}(t)\sqrt\kappa dB_{t_0}(t)+b_{t_0}(t) dt,\quad 0\le t<T-{t_0},$$
where $B_{t_0}(t)$ is given by (\ref{B-c}), $a_{t_0}(t):=a({t_0}+t)$, and $b_{t_0}(t):=b({t_0}+t)$. Note that $B_{t_0}(t)$ is an $(\F_{{t_0}+t})_{t\ge 0}$-Brownian motion, $X_{t_0}(t)$, $a_{t_0}(t)$ and $b_{t_0}(t)$ are all $(\F_{{t_0}+t})_{t\ge 0}$-adapted. \label{SDE-T}
\end{Definition}

\subsection{Special functions} \label{special}
We now introduce some functions that will be used to define annulus Loewner equations. For $t>0$, define $$\SA(t,z)=\lim_{M\to\infty}\sum_{k=-M}^M \frac{e^{2kt}+z}{e^{2kt}-z}
=\PV\sum_{2\mid n} \frac{e^{nt}+z}{e^{nt}-z},$$ 
$$\HA(t,z)=-i\SA(t,e^i(z))=-i\PV\sum_{2\mid n} \frac{e^{nt}+e^{iz}}{e^{nt}-e^{iz}}=\PV\sum_{2\mid n}\cot_2(z-int).$$
Then $\HA(t,\cdot)$ is a meromorphic function in $\C$, whose poles are $\{2m\pi+i2kt:m,k\in\Z\}$, which are
all simple poles with residue $2$. Moreover, $\HA(t,\cdot)$ is an odd function and takes real values on $\R\sem\{\mbox{poles}\}$;
$\Imm \HA(t,\cdot)\equiv -1$ on $\R_t$;  $\HA(t,z+2\pi)=\HA(t,z)$ and $\HA(t,z+i2t)=\HA(t,z)-2i$
for any $z\in\C\sem\{\mbox{poles}\}$.

The power series expansion of $\HA(t,\cdot)$ near $0$ is
 \BGE \HA(t,z)=\frac 2z+ \rA(t)z+O(z^3),\label{Taylor}\EDE
where $ \rA(t)=\sum_{k=1}^\infty \sinh^{-2}(kt)-\frac 16$. 
As $t\to \infty$, $\SA(t,z)\to \frac{1+z}{1-z}$,  $\HA(t,z)\to \cot_2(z)$, and $\rA(t)\to -\frac16$.
So we define $\SA(\infty,z)=\frac{1+z}{1-z}$, $\HA(\infty,z)=\cot_2(z)$, and $\rA(\infty)=-\frac16$.  Then
$\rA$ is continuous on $(0,\infty]$, and (\ref{Taylor}) still holds when $t=\infty$. In fact, we have $\rA(t)-\rA(\infty)
=O(e^{-t})$ as $t\to \infty$, so we may define $\RA$ on $(0,\infty]$ by
$ \RA(t) =-\int_t^\infty (\rA(s)-\rA(\infty))ds$. 
Then $\RA$ is continuous on $(0,\infty]$, $\RA(t)=O(e^{-t})$ as $t\to\infty$, and for $0<t<\infty$, \BGE \RA'(t)=\rA(t)-\rA(\infty).\label{RA}\EDE

Let $\SA_I(t,z)=\SA(t,e^{-t}z)-1$ and $\HA_I(t,z)=-i\SA_I(t,e^{iz})=\HA(t,z+it)+i$.
It is easy to check:  \BGE\SA_I(t,z)=\PV\sum_{2\nmid n} \frac{e^{nt}+z}{e^{nt}-z}, \quad
\HA_I(t,z) =\PV\sum_{2\nmid n} \cot_2(z-int). \label{HA-I}\EDE
So $\HA_I(t,\cdot)$ is a meromorphic function in $\C$ with poles $\{2m\pi+i(2k+1)t:m,k\in\Z\}$, which are
all simple poles with residue $2$;  $\HA_I(t,\cdot)$ is an odd function and takes real values on $\R$; and
 $\HA_I(t,z+2\pi)=\HA_I(t,z)$, $\HA_I(t,z+i2t)=\HA_I(t,z)-2i$ for any $z\in\C\sem\{\mbox{poles}\}$.

It is possible to express $\HA$ and $\HA_I$ using classical functions. Let $\theta(\nu,\tau)$ and $\theta_k(\nu,\tau)$, $k=1,2,3$, be the Jacobi theta functions defined in \cite{elliptic}. Define $\Theta(t,z)=\theta(\frac{z}{2\pi},\frac{it}\pi)$ and $\Theta_I(t,z)=\theta_2(\frac{z}{2\pi},\frac{it}\pi)$. Then $\Theta_I$ has period $2\pi$, $\Theta$ has antiperiod $2\pi$, and
\BGE \HA=2\,\frac{\Theta'}{\Theta},\quad \HA_I=2\,\frac{\Theta_I'}{\Theta_I}.\label{HA-Theta}\EDE
These follow from the product representations of $\Theta$ and $\Theta_I$. For example,
\BGE \Theta_I(t,z)=\prod_{m=1}^\infty (1-e^{-2mt})(1-e^{-(2m-1)t}e^{iz})(1-e^{-(2m-1)t}e^{-iz}).\label{Theta4}\EDE
Both $\Theta$ and $\Theta_I$ solve the heat equation \BGE \pa_t\Theta=\Theta'',\quad  \pa_t\Theta_I=\Theta_I''.\label{heat}\EDE
So $\HA$ and $\HA_I$ solve the PDE: \BGE \pa_t\HA=\HA''+\HA'\HA,\quad \pa_t\HA_I=\HA_I''+\HA_I'\HA_I.\label{pat-HA}\EDE

We rescale the functions $\HA$ and $\HA_I$ as follows. For $t>0$ and $z\in\C$, let \BGE \ha \HA(t,z)=\frac \pi t\HA\Big(\frac{\pi^2}t,\frac\pi t z\Big)+\frac zt,\qquad \ha\HA_I(t,z)=\frac \pi t\HA_I\Big(\frac{\pi^2}t,\frac\pi t z\Big)+\frac zt.\label{ha-HA-0*}\EDE
Since $\ha\HA$ and $\ha\HA_I$ have period $2\pi$,
\BGE \ha\HA(t,z+2k t)=\ha\HA(t,z)+2k,\quad\ha\HA_I(t,z+2k t)=\ha\HA_I(t,z)+2k,\quad k\in\Z.\label{ha-HA-I-period}\EDE
From the identities for $\theta$ in \cite{elliptic} or formula (3) in \cite{ann-prop}, we see
$\HA(t,z)=i\frac\pi t\HA(\frac{\pi^2}t,i\frac\pi t z)-\frac zt$. So
\BGE \ha\HA(t,z)=-i\HA(t,-iz)=\PV\sum_{2\mid n}\coth_2(z-nt).\label{ha-HA*}\EDE
Since $\HA_I(t,z)=\HA(t,z+it)+i$,
\BGE \ha\HA_I(t,z)=\ha\HA(t,z+\pi i) =\PV\sum_{2\mid n}\tanh_2(z-nt).\label{ha-HA-I}\EDE
From (\ref{pat-HA}) and (\ref{ha-HA-0*}) we may check that
\BGE -\pa_t\ha\HA=\ha\HA''+\ha\HA'\ha\HA,\quad -\pa_t\ha\HA_I=\ha\HA_I''+\ha\HA_I'\HA_I.\label{pat-ha-HA}\EDE
From (\ref{ha-HA*}) and (\ref{ha-HA-I}) we see that $\ha\HA(t,\cdot)\to\coth_2$ and $\ha\HA_I(t,\cdot)\to\tanh_2$ as $t\to\infty$.

From (\ref{HA-I}) we see that  as $t\to \infty$, $\HA_I(t,z)\to 0$, so its derivatives about $z$ also tend to $0$. The following lemma gives some estimations of these limits.
\begin{Lemma} If  $|\Imm z|<t$, then
\BGE |\HA_I(t,z)|\le \frac{4e^{|\Imm z|-t}}{(1-e^{|\Imm z|-t})^2(1-e^{2(|\Imm z|-t)})}.\label{est-of-HA-I}\EDE
If $t\ge |\Imm z|+2$, then $|\HA_I(t,z)|< 5.5 e^{|\Imm z|-t}$. For any $h\in\N$, if $t\ge |\Imm z|+h+2$, then
$|\HA_I^{(h)}(t,z)|< 15\sqrt h e^{|\Imm z|-t}$.
\label{estimation}
\end{Lemma}
{\bf Proof.} From (\ref{HA-I}), if $|\Imm z|<t$, then
$$ |\HA_I(t,z)|=\left|\sum_{k=0}^{\infty}\Big( \frac{e^{(2k+1)t}+e^{iz}}{e^{(2k+1)t}-e^{iz}}
+\frac{e^{-(2k+1)t}+e^{iz}}{e^{-(2k+1)t}-e^{iz}}\Big)\right|$$
\BGE=\left|\sum_{k=0}^\infty\frac{2\sin(z)}{\cosh((2k+1)t)-\cos(z)}\right|\le \sum_{k=0}^\infty\frac{2e^{|\Imm z|}}{\cosh((2k+1)t)-\cosh(|\Imm z|)}.
\label{HA-I=}\EDE
Here we use the facts that $|\sin(z)|  \le e^{|\Imm z|}$ and $|\cos(z)|\le\cosh(|\Imm z|)<\cosh(t)$. Let $h_0=t-|\Imm z|>0$. Then for $k\ge 0$,
$$\cosh((2k+1)t)-\cosh(|\Imm z|)= 2\sinh_2((2k+1)t+|\Imm z|)\sinh_2((2k+1)t-|\Imm z|) $$
$$=\frac 12 e^{((2k+1)t+|\Imm z|)/2}(1-e^{-(2k+1)t-|\Imm z|})e^{((2k+1)t-|\Imm z|)/2}(1-e^{-(2k+1)t+|\Imm z|})$$
$$\ge \frac 12 e^{((2k+1)t+|\Imm z|)/2}e^{((2k+1)t-|\Imm z|)/2}(1-e^{-h_0})^2=\frac 12 e^{(2k+1)t}(1-e^{-h_0})^2.$$
So the RHS of (\ref{HA-I=}) is not bigger than
$$\sum_{k=0}^\infty\frac{4e^{|\Imm z|}e^{-(2k+1)t}}{(1-e^{-h_0})^2}=\frac{4e^{|\Imm z|-t}}{(1-e^{-h_0})^2(1-e^{-2t})}\le\frac{4e^{-h_0}}{(1-e^{-h_0})^2(1-e^{-2h_0})}.$$
So we proved (\ref{est-of-HA-I}).

If $t\ge |\Imm z|+2$, then $4/({(1-e^{|\Imm z|-t})^2(1-e^{2(|\Imm z|-t)})})\le 4/((1-e^{-2})^2(1-e^{-4}))<5.5$. From (\ref{est-of-HA-I}) we have $|\HA_I(t,z)|< 5.5 e^{|\Imm z|-t}$. Now we assume $h\in\N$ and $t\ge |\Imm z|+h+2$. Then for any $w\in\C$ with $|w-z|=h$, we have
$t\ge |\Imm w|+2$, so $|\HA_I(t,w)|< 5.5 e^{|\Imm w|-t}\le 5.5 e^he^{|\Imm z|-t}$. From Cauchy's integral formula
and Stirling's formula, we have
$$|\HA_I^{(h)}(t,z)|\le 5.5\frac{h!e^h}{h^h} e^{|\Imm z|-t}\le 5.5\sqrt{2\pi h}e^{1/(12h)} e^{|\Imm z|-t}<15\sqrt h e^{|\Imm z|-t}.
\quad \Box$$

\section{Loewner Equations} \label{Loewner}
\subsection{Whole-plane Loewner equation} \label{section-interior}
  Let $T\in(-\infty,\infty]$ and $\xi\in C((-\infty,T))$. Let $g_I(t,\cdot)$, $-\infty<t<T$, be the solution of
  \BGE\pa_t  g_I(t,z)= g_I(t,z)\,\frac{e^{i\xi(t)}+ g_I(t,z)}{e^{i\xi(t)}- g_I(t,z)},\quad -\infty<t<T,\label{radial-eqn}\EDE
$$ \lim_{t\to -\infty} e^t g_I(t,z)=z,\quad z\in\C\sem\{0\}.$$ 
For each $t\in(-\infty,T)$, let $K_I(t)$ be the set of $z\in\C$ at which $g_I(t,\cdot)$ is not defined. Then $K_I(t)$ and $g_I(t,\cdot)$, $-\infty<t<T$, are called the whole-plane Loewner hulls and maps driven by $\xi$.
\vskip 4mm
\no{\bf Remark.} From Proposition 4.21 in \cite{LawSLE} we know that $K_I(t)$ and $g_I(t,\cdot)$ exist and are determined by $e^{i\xi(s)}$, $-\infty<s\le t$. Moreover, each $g_I(t,\cdot)$ maps $\ha\C\sem K_I(t)$ conformally onto $\ha\C\sem \lin{\D}$ and fixes $\infty$, and $g_I(t,z)=e^{-t}z+O(1)$ near $\infty$. So each $K_I(t)$ is a hull in $\C$ with $\ccap(K_I(t))=e^t$.
\vskip 4mm
We say that $\xi$ generates a whole-plane Loewner trace $\beta_I$ if
$$\beta_I(t):=\lim_{|z|>1,z\to e^{i\xi(t)}} g_I(t,\cdot)^{-1}(z)$$
exists for $t\in(-\infty,T)$, and if $\beta_I(t)$, $-\infty\le t<T$, is a continuous curve in $\C$.
 Such a trace, if it exists, starts from $0$, i.e., $\beta_I(-\infty):=\lim_{t\to-\infty} \beta_I(t)=0$. The trace is called simple if $\beta_I(t)$, $-\infty\le t<T$, has no self-intersection.  If $\xi$ generates a whole-plane Loewner trace $\beta_I$, then for each $t$, $\C\sem K_I(t)$ is the unbounded component of $\C\sem \beta_I([-\infty,t])$. In particular, if $\beta_I$ is simple, then $K_I(t)=\beta_I([-\infty,t])$ for each $t$.
\vskip 4mm
Let $\kappa>0$. A pre-$(\TT;\kappa)$-Brownian motion a.s.\ generates a whole-plane Loewner trace, which is called a standard whole-plane SLE$_\kappa$ trace. The trace goes from $0$ to $\infty$, i.e., $\lim_{t\to\infty}\beta_I(t)=\infty$. If $\kappa\in(0,4]$, the trace is simple. If the driving function is the sum of a pre-$(\TT;\kappa)$-Brownian motion and $s_0t$ for some constant $s_0\in\R$, then we also get a whole-plane Loewner trace, which is called a standard whole-plane SLE$(\kappa,s_0)$ trace. The trace still goes from $0$ to $\infty$ as $t\to\infty$, and is simple when $\kappa\le 4$. For any $z_1\ne z_2\in\ha\C$, we may define whole-plane SLE$_\kappa$ and SLE$(\kappa,s_0)$ trace from $z_1$ to $z_2$ via M\"obius transform.
\vskip 4mm
\no{\bf Remark.}
Whole-plane SLE$_\kappa$ is related to radial SLE in the way that, if $T\in\R$ is fixed, then conditioned on $K_I(t)$, $-\infty<t\le T$, the curve $\beta_I(T+t)$, $t\ge 0$, is the radial SLE$_\kappa$ trace in $\ha\C\sem K_I(T)$ from $\beta_I(T)$ to $\infty$. Whole-plane SLE$(\kappa,s_0)$ is related to radial SLE$(\kappa,-s_0)$ (the radial Loewner process driven by $\sqrt\kappa B(t)-s_0t$) in a similar way. The additional negative sign is due to the convention of definitions.
\vskip 4mm

We will need the following inverted whole-plane Loewner process, which grows from $\infty$.
For $-\infty<t<T$, let $K(t)=I_0(K_I(t))$ and $g(t,\cdot)=I_0\circ g_{I}(t,\cdot)\circ I_0$.
Then for each $t$, $g(t,\cdot)$ maps $\ha\C\sem K(t)$ conformally onto $\D$ and fixes $0$. Moreover, $g(t,\cdot)$ satisfies
 (\ref{radial-eqn}) with some initial value at $-\infty$.
We call $K(t)$ and $g(t,\cdot)$  the inverted whole-plane Loewner hulls and maps driven by $\xi$. If $\xi$ generates a whole-plane Loewner trace $\beta_I$, then $\beta(t):=I_0\circ\beta_I(t)$ is a continuous curve in $\ha\C$ that satisfies $\beta(-\infty)=\infty$ and
$\beta(t)=\lim_{\D\ni z\to e^{i\xi(t)}} g(t,\cdot)^{-1}(z)$, $-\infty<t<T$. 
We call $\beta$ the  inverted whole-plane Loewner trace driven by $\xi$.

Let $K_I(t)$ and $g_I(t,\cdot)$, $-\infty<t<T$, be as before. Let $\til K_I(t)=(e^i)^{-1}(K_I(t))$, $-\infty<t<T$. It is easy to see that there exists a unique family $\til g_I(t,\cdot)$, $-\infty<t<T$,  such that, $\til g_I(t,\cdot)$
maps $\C\sem \til K_I(t)$ conformally onto $-\HH$, $e^i\circ\til g_I(t,\cdot)=g_I(t,\cdot)\circ e^i$, and $\til g_I$ satisfies
\BGE \pa_t{\til g_I}(t,z)=\cot_2(\til g_I(t,z)-\xi(t)),\quad -\infty<t<T, \label{radial-eqn-covering}\EDE
and the initial value at $-\infty$:
$$ \lim_{t\to -\infty} (\til g_I(t,z)-it)=z.$$ 
Then we call $\til K_I(t)$ and $\til g_I(t,\cdot)$ the covering whole-plane Loewner hulls and maps driven by $\xi$.

For $-\infty<t<T$, let $\til K(t)=\til I_0(\til K_I(t))$ and
$\til g(t,\cdot)=\til I_0\circ \til g_I(t,\cdot)\circ \til I_0$. 
Then $\til K(t)=(e^i)^{-1}(K(t))$ and $e^i\circ \til g(t,\cdot)=g(t,\cdot)\circ e^i$.
We call $\til K(t)$ and $\til g(t,\cdot)$ the inverted covering whole-plane Loewner hulls and maps driven by $\xi$. Then for each $t\in(-\infty,T)$, $\til g(t,\cdot)$ maps $\C\sem \til K(t)$ conformally onto $\HH$, and satisfies (\ref{radial-eqn-covering}) for $t\in(-\infty,T)$  and the initial value at $-\infty$:
\BGE \lim_{t\to -\infty} (\til g(t,z)+it)=z.\label{invert-covering-whole-limit}\EDE

\subsection{Annulus Loewner equation}
The annulus Loewner equation was introduced in \cite{Zhan}. Fix $p\in(0,\infty)$. Let $\xi\in C([0,T))$ where $0<T\le p$. The annulus Loewner equation of modulus $p$ driven by $\xi$ is
$$ \pa_t  g(t,z)= g(t,z) \SA(p-t, g(t,z)/e^{i\xi(t)}),\quad  g(0,z)=z.$$ 
For $0\le t<T$, let $K(t)$ denote the set of $z\in\A_p$ such that the
solution $ g(s,z)$ blows up before or at time $t$. We call $K(t)$ and $ g(t,\cdot)$, $0\le t<T$, the annulus
Loewner hulls and maps of modulus $p$ driven by $\xi$. For each $t$, $K(t)$ is a hull in $\A_p$ with $\ccap_{\A_p}(K(t))=t$, and $g(t,\cdot)$ maps $\A_p\sem K(t)$ conformally onto $\A_{p-t}$ and maps $\TT_p$ onto $\TT_{p-t}$. It is clear that, for any $n\in\Z$, $\xi$ and $\xi+2n\pi$ generate the same annulus Loewner hulls and maps.

We say that $\xi$ generates an annulus Loewner trace $\beta$ of modulus $p$ if \BGE \beta(t):=\lim_{\A_{p-t}\ni z\to e^{i\xi(t)}} g(t,\cdot)^{-1}(z)\label{trace-2}\EDE exists for all $t\in[0,T)$, and if $\beta(t)$, $0\le t<T$, is a continuous curve. The curve lies in $\A_p\cup\TT$ and starts from $e^{i\xi(0)}\in\TT$. The trace is called simple if $\beta$ has no self-intersection and stays inside $\A_p$ for $t>0$.

\vskip 4mm
\no{\bf Remarks.}
\begin{enumerate}
  \item If $\xi$ generates an annulus Loewner trace $\beta$, then for each $t$, $\A_p\sem K(t)$ is the component of $\A_p\sem \beta((0,t])$ whose boundary contains $\TT_p$. If the trace is simple, then $K(t)=\beta((0,t])$ for each $t$.
\item  Let $\beta(t)$, $0\le t<T$, be a simple curve with $\beta(0)\in\TT$ and $\beta(t)\in\A_p$ for $0<t<T$. If it is parameterized by capacity in $\A_p$, i.e., $\ccap_{\A_p}(\beta((0,t]))=t$ for each $t$, then it is an annulus Loewner trace of modulus $p$. In the general case, let $u(t)=\ccap_{\A_p}(\beta((0,t]))$. Then $\beta(u^{-1}(t))$ is an annulus Loewner trace of modulus $p$.
  \item If $\xi(t)=\sqrt\kappa B(t)$, $0\le t<p$, then a.s.\ $\xi$ generates an annulus Loewner trace. If $\kappa\in(0,4]$, the trace is simple. From Girsanov theorem, the above still hold if $\xi$ is a semimartingale, whose stochastic part is $\sqrt\kappa B(t)$, and whose drift part is a continuously differentiable function.
\end{enumerate}

The covering annulus Loewner equation of modulus $p$ driven by the above $\xi$ is
\BGE \pa_t{\til g}(t,z)=\HA(p-t,\til g(t,z)-\xi(t)),\quad \til g(0,z)=z.\label{annulus-eqn-covering}\EDE
For $0\le t<T$, let $\til K(t)$ denote the set of $z\in\St_p$ such that the
solution $\til g(s,z)$ blows up before or at time $t$.
Then for $0\le t<T$, $\til g(t,\cdot)$ maps $\St_p\sem \til K(t)$ conformally onto $\St_{p-t}$ and
maps $\R_p$ onto $\R_{p-t}$. We call $\til K(t)$ and $\til  g(t,\cdot)$, $0\le t<T$, the covering annulus
Loewner hulls and maps of modulus $p$ driven by $\xi$. Let $K(t)$ and $ g(t,\cdot)$ be the notations appeared above.
 Then we have $\til K(t)=(e^i)^{-1}(K(t))$ and $e^i\circ \til g(t,\cdot) = g(t,\cdot)\circ e^i$ for $0\le t<T$.

Since $\til g(t,\cdot)$ maps $\R_p$ onto $\R_{p-t}$ and $\HA_I(t,z)=\HA(t,z+it)+i$, we have
$$ \pa_t \Ree \til g(t,z)=\HA_I(p-t,\Ree\til g(t,z)-\xi(t)),\quad z\in\R_p.$$ 
Differentiating the above formula w.r.t.\ $z$, we obtain
\BGE \pa_t {\til g}'(t,z)=\til g'(t,z)\HA_I'(p-t,\Ree\til g(t,z)-\xi(t)),\quad z\in\R_p.\label{deriv2*}\EDE

If $\xi$ generates an annulus Loewner trace $\beta$ of modulus $p$, then a.s.\
$$ \til\beta(t):=\lim_{\St_{p-t}\ni z\to \xi(t)} \til g(t,\cdot)^{-1}(z)$$
exists for $0\le t<T$, and $\til\beta(t)$, $0\le t<T$, is a continuous curve in $\St_p\cup\R$ started from $\til\beta(0)=\xi(0)\in\R$. Such $\til\beta$ is called the covering annulus Loewner trace of modulus $p$ driven by $\xi$. And we have $\beta=e^i\circ \til \beta$. If $\beta$ is a simple trace, then $\til\beta$ has no self-intersection, stays inside $\St_p$ for $t>0$, and $\til K(t)=\til\beta((0,t])+2\pi\Z$ for each $t$.


\vskip 3mm

Let $K_I(t)=I_p(K(t))$, $g_I(t,\cdot)=I_{p-t}\circ g(t,\cdot)\circ I_p$, $\til K_I(t)=\til I_p(\til K(t))$, and $\til g_I(t,\cdot)=\til I_{p-t}\circ \til g(t,\cdot)\circ \til I_p$.Then  $K_I(t)$ is a hull in $\A_p$ with $\ccap_{\A_p}(K_I(t))=t$, and $ g_I(t,\cdot)$ maps $\A_p\sem K_I(t)$ conformally onto $\A_{p-t}$ and maps $\TT$ onto $\TT$.
Moreover, $\til K_I(t)=(e^i)^{-1}(K_I(t))$, $\til g_I(t,\cdot)$ maps $\St_p\sem \til K_I(t)$ conformally
onto $\St_{p-t}$, maps $\R$ onto $\R$, satisfies $e^i\circ \til g_I(t,\cdot)=g_I(t,\cdot)\circ e^i$, and the equation
\BGE \pa_t{\til g}_I(t,z)=\HA_I(p-t,\til g_I(t,z)-\xi(t)),\quad \til g(0,z)=z.  \label{annulus-eqn-inver-covering}\EDE
We call $K_I(t)$ and $ g_I(t,\cdot)$ (resp.\ $\til K_I(t)$ and $\til  g_I(t,\cdot)$) the inverted annulus (resp.\ inverted covering annulus) Loewner hulls and maps of modulus $p$ driven by $\xi$.

\subsection{Disc Loewner equation}
Let $T\in(-\infty,0]$ and $\xi\in C((-\infty,T))$. Let $g_I(t,\cdot)$, $-\infty<t<T$, be the solution of
\BGE\pa_t  g_I(t,z)=g_I(t,z)\SA_I(-t, g_I(t,z)/e^{i\xi(t)}),\quad -\infty<t<T;\label{dot-g-I-deg}\EDE
$$\lim_{t\to-\infty}  g_I(t,z)=z,\quad \forall z\in \lin{\D}\sem\{0\}.$$
For each $t\in(-\infty,T)$, let $K_I(t)$ be the set of $z\in\D$ at which $g_I(t,\cdot)$ is not defined. Then $K_I(t)$ and $g_I(t,\cdot)$, $-\infty<t<T$, are called the disc Loewner hulls and maps driven by $\xi$.

\vskip 4mm
\no{\bf Remark.} From Proposition 4.1 and 4.2 in \cite{Zhan} we know that $K_I(t)$ and $g_I(t,\cdot)$ exist and are determined by $e^{i\xi(s)}$, $-\infty<s\le t$. Moreover, each $g_I(t,\cdot)$ maps $\D\sem K_I(t)$ conformally onto $\A_{-t}$ and maps $\TT$ onto $\TT$.
\vskip 4mm

We say that $\xi$ generates a disc Loewner trace $\beta$ if $$\beta_I(t):=\lim_{\A_{-t}\ni z\to e^{t+i\xi(t)}} g_I(t,\cdot)^{-1}(z)$$
exists for every $t\in(-\infty,T)$, and if $\beta_I(t)$, $-\infty\le t<T$, is a continuous curve in $\D$ with $\beta_I(-\infty)=0$. The trace is called simple if it has no self-intersection. If $\xi$ generates a disc Loewner trace $\beta_I$, then for each $t$, $\C\sem K_I(t)$ is the unbounded component of $\C\sem \beta_I([-\infty,t])$. In particular, if $\beta_I$ is simple, then $K_I(t)=\beta_I([-\infty,t])$ for each $t$.

Let $\beta_I(t)$, $-\infty\le t<T$, be a simple curve in $\D$ with $\beta_I(-\infty)=0$. If it is parameterized by capacity in $\D$, i.e., $\modd(\D\sem\beta_I([-\infty,t]))=-t$ for each $t$, then $\beta_I$ is a disc Loewner trace. In the general case, let $u(t)=-\modd(\D\sem\beta_I([-\infty,t]))$, then $\beta_I(u^{-1}(t))$ is a disc Loewner trace.

We will need the following inverted disc Loewner process, which grows from $\infty$.
For $-\infty<t<T$, let $K(t)=I_0(K_I(t))$ and $g(t,\cdot)=I_{-t}\circ g(t,\cdot)\circ I_0$.
Then each $g(t,\cdot)$ maps $\ha\C\sem\lin{\D}\sem K(t)$ conformally onto $\A_{-t}$ and maps $\TT$ onto $\TT_{-t}$. Moreover, $g(t,\cdot)$ satisfies (\ref{dot-g-I-deg}) with $\SA_I$ replaced by $\SA$.
We call $K(t)$ and $g(t,\cdot)$, $-\infty<t<T$, the inverted disc Loewner hulls and maps driven by $\xi$. If $\xi$ generates a disc Loewner trace $\beta_I$, then $\beta:=I_0\circ\beta_I$ is called the inverted disc Loewner trace driven by $\xi$.

The covering disc Loewner hulls and maps are defined as follows.
Let $\til K_I(t)=(e^i)^{-1}(K_I(t))$, $-\infty<t<T$. There is a unique family $\til g_I(t,\cdot)$, $-\infty<t<T$, which satisfy that, for each $t$, $\til g_I(t,\cdot)$ maps $\HH\sem \til K_I(t)$ conformally onto $\St_{-t}$ and maps $\R$ onto $\R$,
 $e^i\circ\til g_I(t,\cdot)=g_I(t,\cdot)\circ e^i$, and the following hold:
\BGE \pa_t{\til g}_I(t,z)=\HA_I(-t,\til g_I(t,z)-\xi(t)); \label{covering-disc-diff}\EDE
\BGE \lim_{t\to -\infty} \til g_I(t,z)=z.\label{covering-disc-end}\EDE
We call $\til K_I(t)$ and $\til g_I(t,\cdot)$ the covering disc Loewner hulls and maps driven by $\xi$.
Let $\til K(t)=\til I_0(\til K_I(t))$ and $\til g(t,\cdot)=\til I_{-t}\circ \til g_I(t,\cdot)\circ \til I_0$. Then $\til g(t,\cdot)$ maps
$-\HH\sem \til K(t)$ conformally onto $\St_{-t}$, maps $\R$ onto $\R_{-t}$, $e^i\circ\til g(t,\cdot)=g(t,\cdot)\circ e^i$, and satisfies
$\pa_t{\til g}(t,z)=\HA(-t,\til g(t,z)-\xi(t))$. 
We call $\til K(t)$ and $\til g(t,\cdot)$  the inverted covering disc Loewner hulls and maps driven by $\xi$.

\vskip 4mm
\no{\bf Remark.} Now we summarize the conformal maps that appear in the this section so far. The relations between a (inverted) whole-plane, annulus, or disc Loewner map $g(t,\cdot)$ or $g_I(t,\cdot)$ and its corresponding covering map $\til g(t,\cdot)$ or $\til g_I(t,\cdot)$ are $g(t,\cdot)\circ e^i=e^i\circ \til g(t,\cdot)$ and $g_I(t,\cdot)\circ e^i=e^i\circ \til g_I(t,\cdot)$. The relation between the inverted pair $\til g(t,\cdot)$ and $\til g_I(t,\cdot)$ depends on the three cases. For the whole-plane Loewner maps,
$$\til g_I(t,\cdot):\C\sem\til K_I(t)\conf -\HH,\quad \til g(t,\cdot):\C\sem\til K(t)\conf \HH,\quad \til I_0\circ \til g(t,\cdot)=\til g_I(t,\cdot)\circ \til I_0.$$
For the annulus Loewner maps of modulus $p$,
$$\til g(t,\cdot):(\St_p\sem\til K(t);\R_p)\conf (\St_{p-t};\R_{p-t}),\quad \til g_I(t,\cdot):(\St_p\sem \til K_I(t);\R)\conf (\St_{p-t};\R),$$ $$\til I_{p-t}\circ \til g_I(t,\cdot)=\til g(t,\cdot)\circ \til I_p,\quad t\in[0,p).$$
For the disc Loewner maps,
$$\til g_I(t,\cdot):(\HH\sem \til K_I(t);\R)\conf (\St_{-t};\R),\quad \til g(t,\cdot):(-\HH\sem \til K(t);\R)\conf (\St_{-t};\R_{-t}),$$
$$ \til I_{-t}\circ \til g(t,\cdot)=\til g_I(t,\cdot)\circ \til I_0,\quad t\in(-\infty,0).$$
The relation between $g(t,\cdot)$ and $g_I(t,\cdot)$ depends on the three cases in a similar way.

\subsection{SLE with Marked Points}\label{section-reference}
\begin{Definition} A covering crossing annulus drift function is a real valued $C^{0,1}$ differentiable function defined on  $(0,\infty)\times\R$. A covering crossing annulus drift function with period $2\pi$ is called a crossing annulus drift function.
\label{drift-function}
\end{Definition}

\begin{Definition}
Suppose $\Lambda$ is a covering crossing annulus drift function. Let $\kappa> 0$, $p>0$, and $x_0,y_0\in\R$. Let $\xi(t)$, $0\le t< p$, be the maximal solution to the SDE
\BGE d \xi(t)= \sqrt\kappa dB(t)+\Lambda(p-t,\xi(t)-\Ree \til g(t,y_0+ pi))dt,\quad  \xi(0)=x_0,\label{xi-crossing}\EDE
where $\til g(t,\cdot)$, $0\le t<p$, are the covering annulus Loewner maps of modulus $p$ driven by $\xi$. Then the covering annulus Loewner trace of modulus $p$ driven by $\xi$ is called the covering (crossing) annulus SLE$(\kappa,\Lambda)$ trace in $\St_p$ started from $x_0$ with marked point $y_0+p i$. \label{covering-ann-crossing}
\end{Definition}

\begin{Definition} Suppose $\Lambda$ is a crossing annulus drift function. Let $\kappa> 0$, $p>0$, $a\in\TT$
and $b\in\TT_p$. Choose  $x_0,y_0\in\R$ such that $a=e^{ix_0}$ and $b=e^{-p+iy_0}$.   Let $\xi(t)$, $0\le t< p$, be the maximal solution to (\ref{xi-crossing}). The annulus Loewner trace $\beta$ driven by $\xi$ is called the (crossing) annulus SLE$(\kappa,\Lambda)$ trace in $\A_p$ started from $a$ with marked point $b$. 
 \label{crossing-ann}
\end{Definition}

The above definition does not depend on the choices of $x_0$ and $y_0$ because $\Lambda$ has period $2\pi$, and for any $n\in\Z$, the annulus Loewner objects driven by $\xi(t)+2n\pi$ agree with those driven by $\xi(t)$.

A covering chordal-type annulus drift function is a real valued $C^{0,1}$ differentiable function defined on  $(0,\infty)\times(\R\sem 2\pi\Z)$.
The word ``covering'' is omitted if the function has period $2\pi$. If $\Lambda$ is a chordal-type annulus drift function, using the same idea, we may define the annulus SLE$(\kappa,\Lambda)$ processes, where the initial point $a=e^{ix_0}$ and marked point $b=e^{iy_0}$ both lie on $\TT$ and are distinct. The driving function $\xi$ is the solution to (\ref{xi-crossing}) with $\Ree \til g(t,y_0+ pi)$ replaced by $\til g(t,y_0)$.

Via conformal maps, we can then define annulus SLE$(\kappa,\Lambda)$ process and trace in any doubly connected domain started from one boundary point with another boundary point being marked. Here $\Lambda$ is a chordal-type or crossing annulus drift function depending on whether or not the initial point and the marked point lie on the same boundary component. Let $\Lambda_I(p,x)=-\Lambda(p,-x)$, then $\Lambda_I$ is called the dual of $\Lambda$. If $W$ is a conjugate conformal map of $\A_p$, and $\Lambda_I$ is the
 dual of $\Lambda$, then  $(W(K(t)))$ is an annulus SLE$(\kappa,\Lambda_I)$ process in $W(\A_p)$
 started from $W(a)$ with marked point $W(b)$.

\begin{Definition}
  Let $\kappa\ge 0$, $b\in\TT$, and $\Lambda$ be a crossing annulus drift function. Choose $y_0\in\R$ such that $e^{iy_0}=b$. Let $B_*^{(\kappa)}(t)$, $t\in\R$, be a pre-$(\TT;\kappa)$-Brownian motion. Suppose $\xi(t)$, $-\infty<t<0$,  satisfies the following SDE with the meaning in Definition \ref{SDE-T}:
$$ d\xi(t)=dB_*^{(\kappa)}(t)+\Lambda(-t,\xi(t)-\til g_I(t,y_0))dt, \quad -\infty<t<0,$$ 
where $\til g_I(t,\cdot)$ are the disc Loewner maps driven by $\xi$. Then we call the disc Loewner trace driven by $\xi$ the disc SLE$(\kappa,\Lambda)$ trace in $\D$ started from $0$ with marked point $b$. \label{def-disc-refer}
\end{Definition}

Via conformal maps, we can define disc SLE$(\kappa,\Lambda)$ trace in any simply connected domain started from an interior point with a marked boundary point.

\section{Coupling of Two Annulus SLE traces} \label{coupling}
The goal of this section is to prove Theorem \ref{coupling-thm} below, which says that when certain PDE is satisfied, we may couple two annulus SLE$(\kappa;\Lambda)$ processes such that they commute with each other. Although this result will not be used directly in the proof of the whole-plane reversibility, we prove this theorem because on the one hand, the result may be used in the future, and on the other hand, the proof will serve as a reference for a more complicated proof of the theorem about coupling two whole-plane SLE processes.

After some preparation in Section \ref{Trans}, the construction formally starts from Section \ref{ensemble-1}, which resembles Section 3 of \cite{reversibility}. The extra complexity comes from the appearance of covering maps and inverted maps. Then we construct a two-dimensional local martingale $M$ in Section \ref{martingale}, which resembles Section 4 of \cite{reversibility}. In the same subsection, we derive the boundedness of $M$ when the two processes are stopped at some exiting time. In Section \ref{local}, we first construct local commutation couplings using $M$, then construct the global commutation coupling using the coupling technique, and finishes the proof. 

\begin{Theorem}
Let $\kappa>0$ and $s_0\in\R$. Suppose $\Gamma$ is a positive $C^{1,2}$ differentiable function on $(0,\infty)\times\R$ that satisfies \BGE \pa_t\Gamma=\frac\kappa 2\Gamma''+\HA_I\Gamma'+\Big(\frac 3\kappa-\frac 12\Big)\HA_I'\Gamma;\label{dot-Gamma=}\EDE
\BGE \Gamma(t,x+2\pi)=e^{\frac{2\pi s_0}\kappa}\Gamma(t,x),\quad t>0,\,x\in\R.\label{Gamma-period-s}\EDE
Let $\Lambda=\kappa  \frac{\Gamma'}{\Gamma}$.
Then $\Lambda$ is a crossing annulus drift function. Let $\Lambda_1=\Lambda$ and $\Lambda_2$ be the dual of $\Lambda$.
Then for any $p>0$, $a_{1},a_{2}\in\TT$,  there is a coupling of two curves: $\beta_{1}(t)$, $0\le t<p$, and $\beta_{2}(t)$, $0\le t<p$,
such that for $j\ne k\in\{1,2\}$, the following hold.
\begin{enumerate}[(i)]
\item    $\beta_j$ is an annulus SLE$(\kappa,\Lambda_j)$ trace in $\A_p$
 started from $a_{j}$ with marked point $a_{I,k}:=I_p(a_{k})$.
\item    If $t_k<p$ is a stopping time w.r.t.\ $(\beta_{k}(t))$,  then conditioned on $\beta_k(t)$, $0\le t\le t_k$, after a time-change,
$\beta_j(t)$, $0\le t<T_j(t_k)$ is the annulus SLE$(\kappa,\Lambda_j)$ process in a connected component of $\A_p\sem I_p(\beta_k((0,t_k]))$ started from $a_j$ with marked point $I_p(\beta_k(t_k))$, where $T_j(t_k)$ is the first time that $\beta_j$ visits $I_p\circ \beta_k((0,t_k])$, which is set to be $p$ if such time does not exist.
\end{enumerate}
\label{coupling-thm}
\end{Theorem}

\no{\bf Remarks.}
\begin{enumerate}
  \item The $\Lambda$ in the theorem satisfies
  \BGE \pa_t \Lambda=\frac\kappa 2 \Lambda''+\Big(3-\frac\kappa 2\Big)\HA_I''+\Lambda\HA_I'+\HA_I\Lambda'+\Lambda\Lambda'.\label{dot-Lambda=}\EDE
 On the other hand, if $\Lambda$ satisfies (\ref{dot-Lambda=}), then there is $\Gamma$ such that $\Lambda=\kappa\frac{\Gamma'}{\Gamma}$ and (\ref{dot-Gamma=}) holds. 
  \item The theorem also holds for $\kappa=0$ if $\Lambda$ satisfies (\ref{dot-Lambda=}) with $\kappa=0$.
\end{enumerate}

\subsection{Transformations of PDE}\label{Trans}

\begin{Lemma} Let $\sigma,s_0\in\R$. Suppose $\Gamma$, $\Psi$, and $\Psi_{s_0}$ are functions defined on $(0,\infty)\times\R$, which satisfy $\Psi=\Gamma\Theta_I^{\frac 2\kappa}$, $\Psi_s=\Gamma_s\Theta_I^{\frac 2\kappa}$, and $\Psi_{s_0}(t,x)=e^{-\frac{s_0x}\kappa-\frac{s_0^2t}{2\kappa}}\Psi(t,x)$.
Then the following PDEs are equivalent:
\BGE \pa_t\Gamma=\frac\kappa 2\Gamma''+\HA_I\Gamma'+\Big(\sigma-\frac 1\kappa+\frac 12\Big)\HA_I'\Gamma;\label{dot-Gamma-c=}\EDE
\BGE \pa_t \Psi = \frac \kappa 2 \Psi'' + \sigma\HA_I' \Psi;\label{PDE-Psi*} \EDE
 \BGE \pa_t \Psi_{s_0} = \frac \kappa 2 \Psi_{s_0}'' +s_0\Psi_{s_0}'+ \sigma\HA_I' \Psi_{s_0}.\label{PDE-Psi-c} \EDE
\label{Gamma-Psi}
\end{Lemma}
\vskip -5 mm
{\bf Proof.} This follows from  (\ref{HA-Theta}),  (\ref{heat}), and some straightforward computations. $\Box$

\vskip 3mm

\no {\bf Remark.} When $\sigma=\frac 4\kappa-1$, (\ref{dot-Gamma-c=}) agrees with (\ref{dot-Gamma=}).

\begin{Lemma} Let $\sigma,s_0\in\R$. Suppose $\Psi_{s_0}$ is positive, has period $2\pi$, and solves (\ref{PDE-Psi-c}) in $(0,\infty)\times\R$. Then $\Psi_{s_0}(t,x)\to C$ as $t\to\infty$ for some constant $C>0$, uniformly in $x\in\R$.
\label{Psi-uniform}
\end{Lemma}
{\bf Proof.} Fix $t_0>0$ and $x_0\in\R$. For $0\le t<t_0$, let $X_{x_0}(t)=x_0+\sqrt\kappa B(t)+st$ and
$$M(t)= \Psi_{s_0}(t_0-t,X_{x_0}(t))\exp\Big(\sigma\int_0^t \HA_I'(t_0-r,X_{x_0}(r))dr\Big).$$
From (\ref{PDE-Psi-c}) and It\^o's formula, $M(t)$, $0\le t<t_0$, is a local martingale. Since $\Psi_{s_0}$ and $\HA_I'$ are continuous on $(0,\infty)\times\R$ and have period $2\pi$, we see that, for any $t_1\in(0,t_0]$, $M(t)$, $0\le t\le t_0-t_1$, is uniformly bounded, so it is a bounded martingale. Thus,
\BGE \Psi_{s_0}(t_0,x_0)= M(0)=\EE\Big[\Psi_{s_0}(t_1,X_{x_0}(t_0-t_1))\exp\Big(\sigma\int_0^{t_0-t_1} \HA_I'(t_0-r,X_{x_0}(r))dr\Big)\Big].\label{expectation*}\EDE
Now suppose $t_0>t_1\ge 3$. From Lemma \ref{estimation}, we see that,
\BGE \int_0^{t_0-t_1}|\HA_I'(t_0-r,X_{x_0}(r))|dr\le   \int_0^{t_0-t_1} 15 e^{r-t_0}dr\le 15  e^{-t_1}.\label{compare*}\EDE

Let $\eps>0$. Choose $t_1\ge 3$ such that $15\sigma e^{-t_1}<\eps/3$. For $t\in[t_1,\infty)$ and $x\in\R$, define $$\Psi_{s_0,t_1}(t,x)=\EE [\Psi_{s_0}(t_1,X_x(t-t_1))].$$
As $t\to \infty$, the distribution of $e^i(X_x(t-t_1))$ tends to the uniform distribution on $\TT$. Since $\Psi_{s_0}$ is positive, continuous, and has period $2\pi$, we see that $\Psi_{s_0,t_1}(t,x)\to  \frac 1{2\pi} \int_0^{2\pi} \Psi_{s_0}(t_1,x)dx>0$ as $t\to \infty$, uniformly in $x\in\R$. Thus, $\lim_{t\to\infty}\ln(\Psi_{s_0,t_1})$ converges uniformly in $x\in\R$.
So there is $t_2>t_1$ such that if $t_a,t_b\ge t_2$ and $x_a,x_b\in\R$, then $|\ln(\Psi_{t_1}(t_a,x_a))-\ln(\Psi_{t_1}(t_b,x_b))|<\eps/3$.
From (\ref{expectation*}) and (\ref{compare*}) we see that $$ |\ln(\Psi_{s_0}(t,x))-\ln(\Psi_{s_0,t_1}(t,x))|\le 15\sigma e^{-t_1}<\eps/3,\quad t\ge t_1,\,x\in\R.$$ Thus, $|\ln(\Psi_{s_0}(t_a,x_a))-\ln(\Psi_{s_0}(t_b,x_b))|<\eps$ if $t_a,t_b\ge t_2$ and $x_a,x_b\in\R$. So $\lim_{t\to \infty}\ln(\Psi_{s_0})$ converges uniformly in $x\in\R$, which implies the conclusion of the lemma. $\Box$

\begin{Lemma} Let $s_0\in\R$. Suppose $\Gamma$ is positive, satisfies (\ref{Gamma-period-s}), and solves (\ref{dot-Gamma-c=}). Then there is $C>0$ such that 
$\Gamma_{s_0}(t,x):=C^{-1}e^{-\frac{s_0x}\kappa-\frac{s_0^2t}{2\kappa}}\Gamma(t,x)$
has period $2\pi$ and satisfies $\lim_{t\to\infty}\Gamma_{s_0}(t,x)= 1$, uniformly in $x\in\R$.
\label{Gamma-uniform}
\end{Lemma}
 {\bf Proof.} Let $\Psi_{s_0}$ be given by Lemma \ref{Gamma-Psi}. Since $\Theta_I>0$, $\Psi_{s_0}$ is positive and solves (\ref{PDE-Psi-c}). Since $\Gamma$ satisfies (\ref{Gamma-period-s}) and $\Theta_I$ has period $2\pi$, $\Psi_{s_0}$ also has period $2\pi$. From Lemma \ref{Psi-uniform}, there is $C>0$ such that $\Psi_{s_0}\to C$ as $t\to \infty$, uniformly in $x\in\R$. Let $\Gamma_{s_0}(t,x):=C^{-1}e^{-\frac{s_0x}\kappa-\frac{s_0^2t}{2\kappa}}\Gamma(t,x)$. Then $\Gamma_{s_0}=C^{-1}\Psi_{s_0}\Theta_I(t,x)^{-\frac 2\kappa}$. From (\ref{Theta4}), $\Theta_I\to 1$ as $t\to \infty$, uniformly in $x\in\R$. Since $\Theta_I$ has period $2\pi$, we get the desired conclusion. $\Box$

\subsection{Ensemble} \label{ensemble-1}
Let $p>0$ and $\xi_1,\xi_2\in C([0,p))$. For $j=1,2$, let $g_j(t,\cdot)$ (resp.\ $g_{I,j}(t,\cdot)$), $0\le t<p$, be the annulus (resp.\ inverted annulus) Loewner  maps of modulus $p$ driven by $\xi_j$. Let  $\til g_j(t,\cdot)$ and $\til g_{I,j}(t,\cdot)$, $0\le t<p$, $j=1,2$,
 be the corresponding  covering Loewner  maps. Suppose $\xi_j$ generates a simple annulus Loewner trace of modulus $p$: $\beta_j$, $j=1,2$. Let $\beta_{I,j}=I_p\circ \beta_j$, $j=1,2$, be the inverted trace. Define \BGE {\cal D}=\{(t_1,t_2):\beta_{1}((0,t_1])\cap \beta_{I,2}((0,t_2])=\emptyset\}=\{(t_1,t_2):
 \beta_{I,1}((0,t_1])\cap \beta_{2}((0,t_2]) =\emptyset\}.\label{cal-D}\EDE
 For $(t_1,t_2)\in\cal D$, we  define
 \BGE \mA(t_1,t_2)=\modd(\A_p\sem \beta_1([0,t_1])\sem \beta_{I,2}([0,t_2]))=\modd(\A_p\sem \beta_{I,1}([0,t_1])\sem \beta_{2}([0,t_2])).\label{m}\EDE

Fix  any $j\ne k\in\{1,2\}$ and $t_k\in[0,p)$. Let $T_j(t_k)$ be the maximal number such that for any $t_j<T_j(t_k)$, we have $(t_1,t_2)\in\cal D$.  As $t_j\to T_j(t_k)^-$, the spherical distance between $\beta_j((0,t_j])$ and $\beta_{I,k}((0,t_k])$ tends to $0$, so $\mA(t_1,t_2)\to 0$.
For $0\le t_j<T_j(t_k)$, let $\beta_{j,t_k}(t_j)= g_{I,k}(t_k,\beta_j(t_j))$. Then $\beta_{j,t_k}(t_j)$, $0\le t_j<T_j(t_k)$, is a simple curve that starts from $g_{I,k}(t_k,e^{i\xi_j(t_j)})\in\TT$, and stays inside $\A_p$ for $t_j>0$. Let
\BGE v_{j,t_k}(t_j)=\ccap_{\A_{p-t_k}}(\beta_{j,t_k}((0,t_j]))=p-t_k-\mA(t_1,t_2).\label{v=}\EDE
Then $v_{j,t_k}$ is continuous and increasing and  maps $[0,T_j(t_k))$ onto $[0,S_{j,t_k})$ for some $S_{j,t_k}\in(0,p-t_k]$.
Since $\mA\to 0$ as $t_j\to T_j(t_k)$, $S_{j,t_k}=p-t_k$.
Then $\gamma_{j,t_k}(t):= \beta_{j,t_k}(v_{j,t_k}^{-1}(t))$, $0\le t<p-t_k$, are the annulus Loewner trace of modulus $p-t_k$
driven by some $\zeta_{j,t_k}\in C([0,p-t_k))$. Let $\gamma_{I,j,t_k}(t)$ be the corresponding inverted
annulus Loewner trace. Let $ h_{j,t_k}(t,\cdot)$ and $h_{I,j,t_k}(t,\cdot)$
be the corresponding annulus and inverted annulus Loewner maps. Let $\til h_{j,t_k}(t,\cdot)$,
and $\til h_{I,j,t_k}(t,\cdot)$ be the corresponding covering maps.

For $0\le t_j<T_j(t_k)$, let $\xi_{j,t_k}(t_j)$, $\beta_{I,j,t_k}(t_j)$, $g_{j,t_k}(t_j,\cdot)$,
$g_{I,j,t_k}(t_j,\cdot)$, $\til g_{j,t_k}(t_j,\cdot)$, and
$\til g_{I,j,t_k}(t_j,\cdot)$ be the time-changes of $\zeta_{j,t_k}(t)$, $\gamma_{I,j,t_k}(t)$,
$ h_{j,t_k}(t,\cdot)$, $h_{I,j,t_k}(t,\cdot)$, $\til h_{j,t_k}(t,\cdot)$, and $\til h_{I,j,t_k}(t,\cdot)$, respectively,
via the map $v_{j,t_k}$. For example, this means that $\xi_{j,t_k}(t_j)=\zeta_{j,t_k}(v_{j,t_k}(t_j))$ and
$g_{j,t_k}(t_j,\cdot)=h_{j,t_k}(v_{j,t_k}(t_j),\cdot)$.

For $0\le t_j<T_j(t_k)$, let
\BGE G_{I,k,t_k}(t_j,\cdot)= g_{j,t_k}(t_j,\cdot)\circ g_{I,k}(t_k,\cdot)\circ g_j(t_j,\cdot)^{-1},\label{G}\EDE
\BGE \til G_{I,k,t_k}(t_j,\cdot)=\til g_{j,t_k}(t_j,\cdot)\circ\til g_{I,k}(t_k,\cdot)\circ\til g_j(t_j,\cdot)^{-1}.\label{til-G}\EDE
Then $G_{I,k,t_k}(t_j,\cdot)$ maps $\A_{p-t_j}\sem  g_j(t_j,\beta_{I,k}((0,t_k]))$
 conformally onto $ \A_{\mA(t_1,t_2)}$ and maps $\TT$ onto $\TT$; $e^i\circ \til G_{I,k,t_k}(t_j,\cdot)=G_{I,k,t_k}(t_j,\cdot)\circ e^i$; and $\til G_{I,k,t_k}(t_j,\cdot)$ maps $\R$ onto $\R$. Since $\gamma_{j,t_k}(t)= \beta_{j,t_k}(v_{j,t_k}^{-1}(t))$, from (\ref{trace-2}) and a similar formula for $\gamma$, we find that $e^{i\xi_{j,t_k}(t_j)}=G_{I,k,t_k}(t_j,e^{i\xi_j(t_j)})$
 for $0\le t_j<T_j(t_k)$. So there is $n\in\Z$ such that $\til G_{I,k,t_k}(t_j,\xi_j(t_j))=\xi_{j,t_k}(t_j)+2n\pi$ for $0\le t_j<T_j(t_k)$.
Since $\zeta_{j,t_k}+2n\pi$ generates the same annulus Loewner hulls as $\zeta_{j,t_k}$, we may choose $\zeta_{j,t_k}$ such that for $0\le t_j<T_j(t_k)$,
\BGE \xi_{j,t_k}(t_j)=\til G_{I,k,t_k}(t_j,\xi_j(t_j)).\label{eta=W(xi)}\EDE

For $0\le t_j<T_j(t_k)$, let
\BGE A_{j,h}(t_1,t_2)=\til G_{I,k,t_j}^{(h)}(t_k,\xi_j(t_j)),\quad h=1,2,3,
\label{A}\EDE
\BGE A_{j,S}(t_1,t_2)=\frac{A_{j,3}(t_1,t_2)}{A_{j,1}(t_1,t_2)}-\frac 32\Big( \frac{A_{j,2}(t_1,t_2)}{A_{j,1}(t_1,t_2)}\Big)^2.
\label{AS}\EDE
Then $A_{j,S}(t_1,t_2)$ is the Schwarzian derivative of $\til G_{I,k,t_j}(t_k,\cdot)$ at $\xi_j(t_j)$.
A standard argument using Lemma 2.1 in \cite{Zhan} shows that, for  $0\le t_j<T_j(t_k)$,  \BGE v_{j,t_k}'(t_j)=|G_{I,k,t_k}'(t_j,\xi_j(t_j))|^2=\til G_{I,k,t_k}'(t_j,\xi_j(t_j))^2=A_{j,1}(t_1,t_2)^2,\label{v'}\EDE
so from (\ref{v=}) we have
\BGE \pa_j \mA=-A_{j,1}^2.\label{pa1m}\EDE
Moreover, for $0\le t_j<T_j(t_k)$,
\BGE\pa_t {\til g}_{j,t_k}(t_j,z)=A_{j,1}(t_1,t_2)^2\HA(\mA(t_1,t_2),\til g_{j,t_k}(t_j,z)-\xi_{j,t_k}(t_j));\label{dot-til-g-1}\EDE
\BGE\pa_t {\til g}_{I,j,t_k}(t_j,z)=A_{j,1}(t_1,t_2)^2\HA_I(\mA(t_1,t_2),\til g_{I,j,t_k}(t_j,z)-\xi_{j,t_k}(t_j)).\label{dot-til-g-1-I'}\EDE

From (\ref{til-G}) 
we have
\BGE \til G_{I,k,t_k}(t_j,\cdot)\circ \til g_j(t_j,z)=\til g_{j,t_k}(t_j,\cdot)\circ\til g_{I,k}(t_k,z).\label{circ=1}\EDE
Differentiate (\ref{circ=1}) w.r.t.\ $t_j$. Let $w=\til g_j(t_j,z)\to \xi_j(t_j)$. From (\ref{annulus-eqn-covering}),  (\ref{eta=W(xi)}), (\ref{dot-til-g-1}), and (\ref{Taylor}) we get
\BGE \pa_t{\til G}_{I,k,t_k}(t_j,\xi_j(t_j))=-3\til G_{I,k,t_k}''(t_j,\xi_j(t_j))=-3A_{j,2}(t_1,t_2).\label{dot=-3''W}\EDE
Differentiate (\ref{circ=1}) w.r.t.\ $t_j$ and $z$, and let $w=\til g_j(t_j,z)\to \xi_j(t_j)$. Then we get
\BGE\frac{\pa_t{\til G}'_{I,k,t_k}(t_j,\xi_j(t_j))}{{\til G}'_{I,k,t_k}(t_j,\xi_j(t_j))}=
\frac 12\cdot\Big(\frac{A_{j,2}}{A_{j,1}}\Big)^2
-\frac 43\cdot \frac{A_{j,3}}{A_{j,1}} +A_{j,1}^2 \rA(\mA )-\rA(p-t_j).\label{dot-W'-1}\EDE


Note that both $G_{I,k,t_k}(t_j,\cdot)$ and $g_{I,k,t_j}(t_k,\cdot)$ map $\A_{p-t_j}\sem  \beta_{I,k,t_j}((0,t_k])$
conformally onto $\A_{\mA(t_1,t_2)}$ and maps $\TT$ onto $\TT$. So they differ by a multiplicative constant of modulus $1$.
Thus, there is $C_k(t_1,t_2)\in\R$ such that
\BGE\til G_{I,k,t_k}(t_j,\cdot)=\til g_{I,k,t_j}(t_k,\cdot)+C_k(t_1,t_2).\label{CkI}\EDE
Interchanging $j$ and $k$ in (\ref{CkI}), we see that there is $C_j(t_1,t_2)\in\R$ such that
\BGE\til G_{I,j,t_j}(t_k,\cdot)=\til g_{I,j,t_k}(t_j,\cdot)+C_j(t_1,t_2).\label{CjI}\EDE
From (\ref{til-G}) we have
\BGE \til g_{I,j,t_k}(t_j,\cdot)\circ \til g_k(t_k,\cdot)+C_j=\til g_{k,t_j}(t_k,\cdot)\circ \til g_{I,j}(t_j,\cdot),\label{+1}\EDE
\BGE \til g_{I,k,t_j}(t_k,\cdot)\circ \til g_j(t_j,\cdot)+C_k=\til g_{j,t_k}(t_j,\cdot)\circ \til g_{I,k}(t_k,\cdot).\label{+2}\EDE
From the definition of inverted annulus Loewner maps, we have
$$ \til g_{j,t_k}(t_j,\cdot)=\til I_{\mA(t_1,t_2)}\circ \til g_{I,j,t_k}(t_j,\cdot)\circ \til I_{p-t_k},\quad \til g_{j}(t_j,\cdot)=\til I_{p-t_j}\circ \til g_{I,j}(t_j,\cdot)\circ \til I_p;$$ 
$$\til g_{I,k,t_j}(t_k,\cdot)=\til I_{\mA(t_1,t_2)}\circ \til g_{k,t_j}(t_k,\cdot)\circ \til I_{p-t_j},\quad \til g_{I,k}(t_k,\cdot)=\til I_{p-t_k}\circ \til g_k(t_k,\cdot)\circ \til I_p.$$ 
From (\ref{+2}) and the above formulas, we get $\til g_{k,t_j}(t_k,\cdot)\circ \til g_{I,j}(t_j,\cdot)+C_k=\til g_{I,j,t_k}(t_j,\cdot)\circ \til g_k(t_k,\cdot)$. Comparing this formula with (\ref{+1}), we see that $ C_1+C_2\equiv 0$. Now we define $X_1$ and $X_2$ on $\cal D$ by
\BGE X_j(t_1,t_2)=\xi_{j,t_k}(t_j)-\til g_{I,j,t_k}(t_j,\xi_k(t_k))=\til G_{I,k,t_k}(t_j,\xi_j(t_j))-\til g_{I,j,t_k}(t_j,\xi_k(t_k)).\label{Xj}\EDE
From (\ref{CkI}), (\ref{CjI}), and $C_1+C_2\equiv 0$, we have
\BGE X_1+X_2\equiv 0.\label{Xjk}\EDE
Since $\HA_I'''$ is even, we may define $Q$ on $\cal D$ by
\BGE Q=\HA_I'''(\mA,X_1)=\HA_I'''(\mA,X_2).\label{Q}\EDE

Differentiate (\ref{dot-til-g-1-I'}) w.r.t.\ $z$ twice. We get
\BGE\frac{\pa_t{\til g'}_{I,j,t_k}(t_j,z)}{\til g_{I,j,t_k}'(t_j,z)}=A_{j,1}^2
\HA_I'(\mA,\til g_{I,j,t_k}(t_j,z)-\xi_{j,t_k}(t_j)). \label{dot-til-g-1-diff}
\EDE
\BGE \pa_{t}\Big(\frac{\til g_{I,j,t_k}''(t_j,z) }{\til g_{I,j,t_k}'(t_j,z) }\Big)=
A_{j,1}^2\HA_I''(\mA,\til g_{I,j,t_k}(t_j,z)-\xi_{j,t_k}(t_j))\til g_{I,j,t_k}'(t_j,z).\label{dot-til-g-1-diff-2}
\EDE

Let $z=\xi_k(t_k)$ in (\ref{dot-til-g-1-I'}), (\ref{dot-til-g-1-diff}), and (\ref{dot-til-g-1-diff-2}).
Since $\HA_I$ and $\HA_I''$ are odd and $\HA_I'$ is even, from (\ref{CjI}) and (\ref{Xj}) we have
\BGE \pa_j \til g_{I,j,t_k}(t_j,\xi_k(t_k))
=-A_{j,1}^2\HA_I(\mA,X_j).\label{d-til-g=X}\EDE
\BGE \frac{\pa_j A_{k,1}}{ A_{k,1}}
=A_{j,1}^2 \HA_I'(\mA,X_j).\label{pa-A-k-alpha}\EDE
\BGE \pa_j\Big( \frac{A_{k,2}}{A_{k,1}}\Big) =-A_{j,1}^2\HA_I''(\mA,X_j)A_{k,1}.\label{d-til-g''=X}\EDE

Differentiate (\ref{dot-til-g-1-diff-2}) w.r.t.\ $z$ again, and let $z=\xi_k(t_k)$. Since  $\HA_I'''$ is
even, we get
$$\pa_{j}\Big(\frac{A_{k,3}}{A_{k,1} }
-\Big(\frac{A_{k,2} }{A_{k,1} }\Big)^2\Big)
 =A_{j,1}^2[\HA_I'''(\mA,X_j)A_{k,1} ^2
-\HA_I''(\mA,X_j)A_{k,2}], $$ 
which together with (\ref{Q}) and (\ref{d-til-g''=X}) implies that
\BGE \pa_j A_{k,S}= A_{j,1}^2 A_{k,1}^2  Q.  \label{pa-AS}\EDE
Define $F$ on $\cal D$ by  \BGE F(t_1,t_2)=\exp\Big(\int_0^{t_2}\!\int_0^{t_1}
 A_{1,1}(s_1,s_2)^2  A_{2,1}(s_1,s_2)^2  Q(s_1,s_2)ds_1ds_2\Big),\label{F}\EDE
Since $\til g_{I,j,t_k}(0,\cdot)=\til h_{I,j,t_k}(0,\cdot)=\id$, when $t_j=0$, we have $A_{k,1}=1$,
$A_{k,2}=A_{k,3}=0$, hence $A_{k,S}=0$. From (\ref{pa-AS}),   we see that
\BGE \frac{\pa_{j} F}{F }=A_{j,S}.\label{paF}\EDE

\vskip4mm
\no{\bf Remark.} There is an explanation of $F$ in terms of Brownian loop measure. If $R$ is a function on $(0,\infty)$ that satisfies $R'(t)=\rA(t)+\frac 1t$, then
$$-\frac 13\ln F(t_1,t_2)-R(t_1,t_2)+R(t_1,0)+R(0,t_2)-R(0,0)$$
is the Brownian loop measure of the loops in $\A_p$ that intersect both $\beta_1([0,t_1])$ and $\beta_{I,2}([0,t_2])$.

\subsection{Martingales in two time variables} \label{martingale}
 Let $a_{1},a_{2}\in\TT$ be as in Theorem \ref{coupling-thm}. Let $a_{I,j}=I_p(a_{j})\in\TT_p$, $j=1,2$.
 Choose $x_1,x_2\in\R$ such that $a_j=e^{ix_j}$, $j=1,2$.
Let $B_1(t)$ and $B_2(t)$ be two independent Brownian motions. For $j=1,2$, let $(\F^j_t)$ be the complete filtration generated by $(B_j(t))$. Let $\Gamma$, $\Lambda$, $\Lambda_1$, and $\Lambda_2$ be as in Theorem \ref{coupling-thm}.
Since $\Gamma$ satisfies (\ref{Gamma-period-s}), $\Lambda_j$, $j=1,2$, has period $2\pi$, which implies that they are annulus drift functions. For $j=1,2$, let  $\xi_j(t_j)$, $0\le t_j<p$, be the solution to the SDE:
\BGE d \xi_j(t_j)=\sqrt\kappa d B_j(t_j)+\Lambda_{j}(p-t_j,\xi_j(t_j)-\til g_{I,j}(t_j,x_{3-j}))d t_j, \quad
\xi_j(0)=x_j.\label{xi-j}\EDE Then $(\xi_1)$ and $(\xi_2)$ are independent.
For simplicity, suppose $\kappa\in(0,4]$ (for the case $\kappa>4$, we may work on Loewner chains and apply Proposition 2.1 in \cite{Zhan}). Then for $j=1,2$, a.s.\ $(\xi_j)$ generates a simple annulus Loewner trace $\beta_j$, which is an annulus SLE$(\kappa,\Lambda_j)$ trace $\beta_j$ in $\A_p$ started from $a_{j}$ with marked point $a_{I,3-j}$. We may apply the results in the prior subsection.

As the annulus Loewner objects driven by $\xi_j$,   $\beta_j$, $\beta_{I,j}=I_p\circ\beta_j$, $(g_{I,j}(t_j,\cdot))$,
 $(\til g_j(t_j,\cdot))$,  and $(\til g_{I,j}(t_j,\cdot))$ are all $(\F^j_{t_j})$-adapted. Fix $j\ne k\in\{1,2\}$.
Since  $\beta_j$ is $(\F^j_{t_j})$-adapted and $(g_{I,k}(t_k,\cdot))$  is $(\F^k_{t_k})$-adapted, we see that
$(t_1,t_2)\mapsto\beta_{j,t_k}(t_j)=g_{I,k}(t_k,\beta_j(t_j))$ defined on $\cal D$
 is $(\F^1_{t_1}\times \F^2_{t_2})$-adapted. Since $\til g_{j,t_k}(t_j,\cdot)$ and
 $\til g_{I,j,t_k}(t_j,\cdot)$ are determined by $\beta_{j,t_k}(s_j)$, $0\le s_j\le t_j$,
 they are $(\F^1_{t_1}\times \F^2_{t_2})$-adapted.
 From (\ref{til-G}), $(\til G_{I, k,t_k}(t_j,\cdot))$ is $(\F^1_{t_1}\times \F^2_{t_2})$-adapted.
 From (\ref{eta=W(xi)}), $(\xi_{j,t_k}(t_j))$ is also $(\F^1_{t_1}\times \F^2_{t_2})$-adapted.
 From (\ref{m}), (\ref{Xj}),  (\ref{A}), and (\ref{AS}),  we see that $(\mA)$,
 $(X_j)$, $(A_{j,h})$, $h=1,2,3$, and $(A_{j,S})$ are all $(\F^1_{t_1}\times \F^2_{t_2})$-adapted.

Fix $j\ne k\in\{1,2\}$ and any $(\F^k_t)$-stopping time $t_k\in[0,p)$.
Let $\F^{j,t_k}_{t_j}=\F^j_{t_j}\times \F^k_{t_k}$, $0\le t_j<p$. Then $(\F^{j,t_k}_{t_j})_{0\le t_j<p}$ is
a filtration. Since $(B_j(t_j))$ is independent of $\F^k_{t_k}$, it is also an  $(\F^{j,t_k}_{t_j})$-Brownian motion.
Thus, (\ref{xi-j}) is an  $(\F^{j,t_k}_{t_j})$-adapted SDE.
From now on, we will apply It\^o's formula repeatedly, all SDE will be $(\F^{j,t_k}_{t_j})$-adapted,
and $t_j$ ranges in $[0,T_j(t_k))$.

From (\ref{dot=-3''W}), (\ref{Xj}), (\ref{A}), and (\ref{d-til-g=X}), we see that $X_j$ satisfies
\BGE \pa_j X_j =A_{j,1}\pa \xi_j(t_j)+\Big(\frac\kappa 2-3\Big)A_{j,2}\pa t_j+A_{j,1}^2 \HA_I(\mA,X_j)\pa t_j.\label{dX}\EDE
Let $\Gamma_{1}=\Gamma$ and $\Gamma_{2}(t,x)=\Gamma(t,-x)$. Then for $j=1,2$, $\Lambda_j=\frac{\Gamma_j'}{\Gamma_j}$ and $\Gamma_{j}$ satisfies (\ref{dot-Gamma=}). From (\ref{Xjk}), we may define $Y$ on $\cal D$ by
\BGE Y=\Gamma_{1}(\mA,X_1)=\Gamma_{2}(\mA,X_2).\label{Y}\EDE
From  (\ref{dot-Gamma=}), (\ref{pa1m}),  (\ref{dX}), and  (\ref{Y}),  we have
\BGE\frac{\pa_j Y}{Y} =\frac1\kappa \Lambda_{j} (\mA,X_j) A_{j,1}{\pa \xi_j(t_j)}-\Big(\frac 3\kappa-\frac 12\Big)\Big(
 A_{j,1}^2\HA_I'(\mA,X_j)+\Lambda_{j}(\mA,X_j)A_{j,2}\Big)\pa t_j.\label{dY}\EDE

 From (\ref{dot-W'-1})  we have
$$\frac{\pa_j A_{j,1} }{ A_{j,1} }=\frac{A_{j,2}}{A_{j,1}}\cdot\pa \xi_j(t_j)
+ \Big(\frac 12\cdot\Big(\frac{A_{j,2} }{A_{j,1}  }\Big)^2+\Big(\frac\kappa 2-\frac43\Big)\cdot\frac{A_{j,3} }{A_{j,1}}\Big)\pa t_j +   A_{j,1}^2  \rA(\mA)\pa t_j-  \rA(p-t_j)\pa t_j.$$
Let $$ \alpha=\frac{6-\kappa}{2\kappa},\qquad c=\frac{(8-3\kappa)(\kappa-6)}{2\kappa}.$$ 
Actually,  $c$ is the central charge for SLE$_\kappa$.  Then we compute
\BGE\frac{\pa_j A_{j,1}^{\alpha}}{ A_{j,1}^{\alpha}}=\alpha\cdot\frac{A_{j,2}}{A_{j,1}}\cdot\pa \xi_j(t_j)
+\frac c 6 A_{j,S}\pa t_j
 +\alpha   A_{j,1}^2  \rA(\mA)\pa t_j-\alpha \rA(p-t_j)\pa t_j.\label{pa-A-j-alpha}\EDE

Recall the  $\RA$ defined in Section \ref{special}. Define $\ha M$ on $\cal D$ by
\BGE\ha M=A_{1,1}^\alpha A_{2,1}^\alpha F^{-\frac c6}Y\exp(\alpha \RA(\mA)).\label{haM}\EDE
Then $\ha M$ is positive. From (\ref{RA}), (\ref{pa1m}), (\ref{pa-A-k-alpha}), (\ref{paF}), (\ref{dY}), and (\ref{pa-A-j-alpha}), we have
\BGE \frac{\pa_j\ha M}{\ha M}=\alpha\frac{A_{j,2}}{A_{j,1}}\pa \xi_j(t_j)+
\frac{A_{j,1}}\kappa \Lambda_{j}(\mA,X_j) {\pa \xi_j(t_j)}-\alpha\rA(p-t_j)\pa t_j+\alpha A_{j,1}^2\rA(\infty)\pa t_j.\label{pahaM}\EDE
When $t_k=0$, we have $A_{j,1}=1$, $A_{j,2}=0$, $\mA=p-t_j$, and $X_j=\xi_j(t_j)-\til g_{I,j}(t_j,x_k)$, so the RHS of (\ref{pahaM}) becomes
$\frac 1\kappa \Lambda_{j}(p-t_j,\xi_j(t_j)-\til g_{I,j}(t_j,x_k))\pa \xi_j(t_j)$.
Define $M$ on  $\cal D$ by
\BGE M(t_1,t_2)=\frac{\ha M(t_1,t_2) \ha M(0,0)}{\ha M(t_1,0)\ha M(0,t_2)}.\label{M}\EDE
Then $M$ is also positive, and $M(\cdot,0)\equiv M(0,\cdot)\equiv 1$. From (\ref{xi-j}) and (\ref{pahaM}) we have
\BGE\frac{\pa_j M}{ M}=\Big[ \Big(3-\frac\kappa2\Big) \frac{A_{j,2}}{A_{j,1}}+\Lambda_{j} (\mA,X_j) A_{j,1}
-\Lambda_{j}(p-t_j,\xi_j(t_j)-\til g_{I,j}(t_j,x_k)) \Big] \frac{\pa B_j(t_j)}{\sqrt\kappa}.\label{paM}\EDE
 So when $t_k\in [0,p)$ is a fixed $(\F^k_t)$-stopping time, $M$ is a local martingale in $t_j$.

\vskip 3mm

Let $\cal J$ denote the set of Jordan curves in $\A_p$ that separate $\TT$ and $\TT_p$. For $J\in\cal J$ and $j=1,2$,
let $T_j(J)$ be the first time that $\beta_j$ visits $J$. It is also the first time that $\beta_{I,j}$ visits $I_p(J)$.
Let $\JP$ denote the set of pairs $(J_1,J_2)\in{\cal J}^2$ such that $I_p(J_1)\cap J_2=\emptyset$ and $I_p(J_1)$ is surrounded by  $J_2$.
This is equivalent to that $I_p(J_2)\cap J_1=\emptyset$ and $I_p(J_2)$ is surrounded by  $J_1$.
Then for every $(J_1,J_2)\in\JP$, $\beta_{I,1}((0,t_1])\cap \beta_{2}((0,t_2])=\emptyset$ when $t_1\le T_{1}(J_1)$ and $t_2\le T_2(J_2)$, so
$[0,T_{1}(J_1)]\times[0,T_{2}(J_2)]\subset \cal D$.

\begin{Lemma} There are positive continuous functions $N_L(p)$ and $N_S(p)$ defined on $(0,\infty)$ that satisfies
$N_L(p),N_S(p)=O(pe^{-p})$ as $p\to \infty$ and the following properties.
Suppose $K$ is an interior hull in $\D$ containing $0$,
$g$ maps $\D\sem K$ conformally onto $\A_p$ for some $p\in(0,\infty)$ and maps $\TT$ onto $\TT$, and $\til g$
is an analytic function that satisfies $e^i\circ \til g=g\circ e^i$. Then for any $x\in\R$,
$|\ln (\til g'(x))|\le N_L(p)$ and $|S\til g(x)|\le N_S(p)$, where $S\til g(x)$ is the Schwarzian derivative of $\til g$ at
$x$, i.e., $S\til g(x)=\til g'''(x)/\til g'(x)-\frac 32 (\til g''(x)/\til g'(x))^2$.
\label{modulus-bound}
\end{Lemma}
{\bf Proof.} Let $f=g^{-1}$ and $\til f=\til g^{-1}$. Then $e^i\circ\til f=f\circ e^i$. Since $\til f'(\til g(x))=1/\til g'(x)$ and $S\til f(\til g(x))=-S\til g(x)/\til g'(x)^2$, we suffice to prove the lemma for $\til f$. Let $P(p,z)=-\Ree \SA_I(p,z)-\ln|z|/p$ and $\til P(p,z)=P(p,e^{iz})=\Imm \HA_I(p,z)+\Imm z/p$. Then $P(p,\cdot)$ vanishes on $\TT$ and $\TT_p\sem\{e^{-p}\}$ and is harmonic inside $\A_p$.
Moreover, when $z\in\A_p$ is near $e^{-p}$, $P(p,z)$ behaves like $-\Ree(\frac{e^{-p}+z}{e^{-p}-z})+O(1)$.
Thus, $-P(p,\cdot)$ is a renormalized  Poisson kernel in $\A_p$ with the pole at $e^{-p}$.
Since $\ln|f|$ is negative and harmonic in $\A_p$ and vanishes on $\TT$,
there is a positive measure $\mu_K$ on $[0,2\pi)$ such that
$$\ln|f(z)|=-\int  P(p,z/e^{i\xi})d\mu_K(\xi),\quad z\in\A_p,$$
which implies that
$$\Imm \til f(z)=\int  P(p,e^{iz}/e^{i\xi})d\mu_K(\xi)=\int \til P(p,z-\xi)d\mu_K(\xi),\quad z\in\St_p$$
So for any $x\in\R$ and $h=1,2,3$, $\til f^{(h)}(x)=\int  \frac{\pa^{h}}{\pa x^{h-1}\pa y}\til P(p,x-\xi)d\mu_K(\xi)$. 
Let $$m_p=\inf_{x\in\R} \frac{\pa}{\pa y}\til P(p,x),\quad M_p=\sup_{x\in\R} \frac{\pa}{\pa y}\til P(p,x)
,\quad M_p^{(h)}=\sup_{x\in\R} |\frac{\pa^{h}}{\pa_x^{h-1}\pa y}\til P(p,x)|, \quad h=2,3.$$
We have $0<m_p<M_p<\infty$ and  $m_p|\mu_K|\le \til f'\le M_p |\mu_K|$ on $\R$. Since $\til f(2\pi)=\til f(0)+2\pi$, we get $1/M_p\le |\mu_K|\le 1/m_p$. Thus, $m_p/M_p\le  \til f'\le M_p/m_p$ and $|\til f^{(h)}|\le M^{(h)}_p/m_p$, $h=2,3$, from which follows that $|S\til f| \le \frac{M_p^{(3)}M_p}{m_p^2}+\frac 32(\frac {M_p^{(2)}M_p}{m_p^2})^2$ on $\R$. Since $\til P(p,z)=\Imm \HA_I(p,z)+\Imm z/p$, we see that $ \frac{\pa}{\pa y}\til P(p,x)=\HA_I'(p,x)+\frac 1p$ and  $\frac{\pa^h}{\pa x^{h-1}\pa y}\til P(p,x)=\HA_I^{(h)}(p,x)$, $h=2,3$. From Lemma \ref{estimation}, $M_p,m_p=\frac 1p+O(e^{-p})$ and $M_p^{(h)}=O(e^{-p})$, $h=2,3$, as $p\to \infty$. So we have $\ln(M_p/m_p)=O(pe^{-p})$ and $\frac{M_p^{(3)}M_p}{m_p^2}+\frac 32(\frac {M_p^{(2)}M_p}{m_p^2})^2=O(pe^{-p})$.  $\Box$

\begin{Proposition}  (Boundedness) Fix $(J_1,J_2)\in\JP$. Then
$|\ln(M)|$ is bounded on $[0,T_1(J_1)]\times[0,T_2(J_2)]$ by a constant
depending only on $J_1$ and $J_2$. \label{bounded}
\end{Proposition}
{\bf Proof.} In this proof, we say a function is uniformly bounded if its values on $[0,T_1(J_1)]\times[0,T_2(J_2)]$
are bounded in absolute value by a constant depending only on $p$, $J_1$, and $J_2$.
If there is no ambiguity, let $\Omega(A,B)$ denote the domain bounded by sets $A$ and $B$, and let $\modd(A,B)$ denote the modulus of this domain if it is doubly connected. Let $J_{I,2}=I_0(J_2)$. Let $p_0=\modd(J_1,J_{I,2})>0$. If $t_1\le T_1(J_1)$ and $t_2\le T_2(J_2)$, since $\Omega(J_1,J_{I,2})$ disconnects $K_1(t_1)$ and $K_{I,2}(t_2)$ in $\A_p$,  $\mA(t_1,t_2)\ge p_0$. Since $\mA\le p$ always holds,
 $\mA\in [p_0,p]$ on $[0,T_{1}(J_1)]\times [0,T_{2}(J_2)]$.
Since $\RA$ is continuous on $(0,\infty)$, $\RA(\mA)$ is uniformly bounded. Since $Q=\HA_I'''(\mA,X_1)$ and $\HA_I'''$ is continuous and has period $2\pi$, $Q$ is uniformly bounded. From Lemma \ref{modulus-bound}, for $j=1,2$, $|\ln(A_{j,1})|\le N_L(\mA)$, so it is uniformly bounded. From (\ref{paF}), $\ln(F)$ is uniformly bounded. Let $s_0\in\R$ be as in Theorem \ref{coupling-thm}. Let $\Gamma_{s_0}>0$ be defined by Lemma \ref{Gamma-uniform}, and $Y_{s_0}=\Gamma_{s_0}(X_1)$. Then $\Gamma_{s_0}$ has period $2\pi$. So $\ln(Y_{s_0})$ is uniformly bounded.  Define $\ha M_{s_0}$ and $M_{s_0}$ using (\ref{haM}) and (\ref{M}) with $Y$ and $\ha M$ replaced by $Y_{s_0}$ and $\ha M_{s_0}$, respectively. Then $\ln(\ha M_{s_0})$ and $\ln(M_{s_0})$ are uniformly bounded because their factors are. Now it suffices to show that $\ln(M)-\ln(M_{s_0})$ is uniformly bounded.
We have
$$\ln(M(t_1,t_2))-\ln(M_{s_0}(t_1,t_2))=\frac {s_0}\kappa(X_1(t_1,t_2)-X_1(t_1,0)-X_1(0,t_2)+X_1(0,0))$$
$$+\frac{s_0^2}{2\kappa}(\mA(t_1,t_2)-\mA(t_1,0)-\mA(0,t_2)+\mA(0,0)).$$
The second term on the RHS of the above formula is uniformly bounded because $\mA\in [p_0,p]$. So it suffices to show that $ X_1(t_1,t_2)-X_1(t_1,0)-X_1(0,t_2)+X_1(0,0)$ is uniformly bounded. Let
$$\til G(t_1,t_2)=\til G_{I,2,t_2}(t_1,\xi_1(t_1)),\quad \til g(t_1,t_2)=\til g_{I,1,t_2}(t_1,\xi_2(t_2)).$$ From (\ref{Xj}) we have $X_1=\til G-\til g$. So it suffices to show that $ \til G(t_1,t_2)-\til G(t_1,0)-\til G(0,t_2)+\til G(0,0)$ and $ \til g(t_1,t_2)-\til g(t_1,0)-\til g(0,t_2)+\til g(0,0)$ are both uniformly bounded. From (\ref{dot-til-g-1-I'})  we have
$\pa_1 \til g(t_1,t_2)=A_{1,1}^2\HA_I(\mA(t_1,t_2),\til g(t_1,t_2)-\xi_{1,t_2}(t_1))$.
Since $A_{1,1}^2$ is uniformly bounded, $\mA\in[p_0,p]$, and $\HA_I$ is continuous and has period $2\pi$, $\til g(t_1,t_2)-\til g(0,t_2)$ is uniformly bounded. Thus, $ \til g(t_1,t_2)-\til g(t_1,0)-\til g(0,t_2)+\til g(0,0)$ is uniformly bounded. Let $\til G_d(t_1,t_2)=\til G(t_1,t_2)-\xi_1(t_1)$. Then $ \til G(t_1,t_2)-\til G(t_1,0)-\til G(0,t_2)+\til G(0,0)= \til G_d(t_1,t_2)-\til G_d(t_1,0)-\til G_d(0,t_2)+\til G_d(0,0)$. To finish the proof it suffices to show that $\til G_d$ is uniformly bounded.

Let $J$ be a Jordan curve which is disjoint from $J_1$ and $I_p(J_2)$, and separates these two curves. Let $\til J=(e^i)^{-1}(J)$.
Since $\til G_d(t_1,t_2)=\til G_{I,2,t_2}(t_1,\xi_1(t_1))-\xi_1(t_1)$, from the Maximum principle, we suffice to show that $\sup_{z\in\til g_1(t_1,\til J)} (\til G_{I,2,t_2}(t_1,z)-z)$ is uniformly bounded. Recall from (\ref{til-G}) that $\til G_{I,2,t_1}(t_1,\cdot)=\til g_{1,t_2}(t_1,\cdot)\circ \til g_{I,2}(t_2,\cdot)\circ \til g_1(t_1,\cdot)^{-1}$. So we suffice to show that the following three quantities are uniformly bounded:
$$\sup_{z\in\til J}|\til g_1(t_1,z)-z|,\quad \sup_{z\in\til J}|\til g_{I,2}(t_2,z)-z|,\quad
\sup_{z\in\til g_{I,2}(t_2,\til J)} |\til g_{1,t_2}(t_1,z)-z|.$$

The uniformly boundedness of these quantities follow from similar arguments. We only work on the last one since it is the hardest. From (\ref{dot-til-g-1}) we have
$$\til g_{1,t_2}(t_1,z)-z=\int_0^{t_1} A_{1,1}(s,t_2)^2\HA(\mA(s,t_2),\til g_{1,t_2}(s,z)-\xi_{1,t_2}(s))ds.$$
Since $\int_0^{t_1}A_{1,1}(s,t_2)^2ds=\mA(0,t_2)-\mA(t_1,t_2)$ is uniformly bounded, we suffice to show that
$$\sup_{z\in \til g_{I,2}(t_2,\til J)}|\HA(\mA(t_1,t_2),\til g_{1,t_2}(t_1,z)-\xi_{1,t_2}(t_1))|$$
is uniformly bounded. From the properties of $\HA$, we suffice to show that there is a constant $h>0$ such that $\Imm \til g_{1,t_2}(t_1,\cdot)\circ\til g_{I,2}(t_2,z)\ge h$ for any $z\in\til J$. This is equivalent to that $|g_{1,t_2}(t_1,\cdot)\circ g_{I,2}(t_2,z)|\le e^{-h}$ for any $z\in J$. We suffice to show that the extremal distance (c.f.\ \cite{Ahl}) between $\TT$ and $g_{1,t_2}(t_1,\cdot)\circ g_{I,2}(t_2,J)$ is bounded below by a positive constant depending only on $p$, $J$, $J_1$ and $J_2$. From conformal invariance, that is equal to the extremal distance between $J$ and $\TT_p\cup \beta_I((0,t_2])$, which is not smaller than the extremal distance between $J$ and $I_p(J_2)$ since $I_p(J_2)$ separates $J$ from $\TT_p\cup \beta_I((0,t_2])$. So we are done. $\Box$

\subsection{Local couplings and global coupling}\label{local}
Let $\mu_j$ denote the distribution of $(\xi_j)$, $j=1,2$. Let $\mu=\mu_1\times\mu_2$. Then
$\mu$ is the joint distribution of $(\xi_1)$ and $(\xi_2)$, since $\xi_1$ and $\xi_2$ are independent.
 Fix $(J_1,J_2)\in\JP$. From  the local martingale property  of $M$ and  Proposition \ref{bounded}, we have
$ \EE_\mu[M(T_1(J_1),T_2(J_2))]=M(0,0)=1$. Define $\nu_{J_1,J_2}$ by
$d\nu_{J_1,J_2}/d\mu=M(T_1(J_1),T_2(J_2))$. Then $\nu_{J_1,J_2}$ is a probability measure.
Let $\nu_1$ and $\nu_2$ be the two marginal measures of $\nu_{J_1,J_2}$. Then
$d\nu_1/d\mu_1=M(T_1(J_1),0)=1$ and $d\nu_2/d\mu_2=M(0,T_2(J_2))=1$, so $\nu_j=\mu_j$, $j=1,2$.
Suppose temporarily that the joint distribution of $(\xi_1)$ and $(\xi_2)$ is $\nu_{J_1,J_2}$ instead of $\mu$.
Then the distribution of each $(\xi_j)$ is still $\mu_j$.

Fix an $(\F^2_t)$-stopping time $t_2\le T_2(J_2)$. From  (\ref{xi-j}), (\ref{paM}), and Girsanov theorem
(c.f.\ \cite{RY}), under the probability measure $\nu_{J_1,J_2}$,
there is an $(\F^1_{t_1}\times\F^2_{t_2})_{t_1\ge 0}$-Brownian
motion $\til B_{1,t_2}(t_1)$  such that $\xi_1(t_1)$,  $0\le t_1\le T_1(J_1)$, satisfies the
$(\F^1_{t_1}\times\F^2_{ t_2})_{t_1\ge 0}$-adapted SDE: $$ d\xi_1(t_1)=\sqrt\kappa d \til B_{1,t_2}(t_1)+\Big(3-\frac\kappa 2\Big)
\frac{A_{1,2} }{A_{1,1} }dt_1+\Lambda_{1}(\mA ,X_1 )A_{1,1}  dt_1,$$ 
which together with (\ref{eta=W(xi)}), (\ref{dot=-3''W}), and It\^o's formula implies that
$$ d\xi_{1,t_2}(t_1)=A_{1,1}  \sqrt\kappa d\til B_{1,t_2}(t_1) +A_{1,1} ^2\Lambda_{1}(\mA ,\xi_{1,t_2}(t_1)
-\til g_{I,1,t_2}(t_1,\xi_2(t_2))) dt_1.$$ 
Recall that $\zeta_{1,t_2}(s_1)=\xi_{1,t_2}(v_{1,t_2}^{-1}(s_1))$
and $\til h_{I,1,t_2}(s_1,\cdot)=\til g_{I,1,t_2}(v_{1,t_2}^{-1}(s_1),\cdot)$. So from (\ref{v=}) and (\ref{v'}),
there is another Brownian motion $\ha B_{1,t_2}(s_1)$ such that for $0\le s_1\le v_{1,t_2}(T_1(J_1))$,
$$ d\zeta_{1,t_2}(s_1)=\sqrt\kappa d\ha B_{1,t_2}(s_1)+\Lambda_{1}(p-t_2-s_1,\zeta_{1,t_2}(s_1)
-\til h_{I,1,t_2}(s_1,\xi_2(t_2)))ds_1.$$ 
Moreover, the initial values is
$ \zeta_{1,t_2}(0)=\xi_{1,t_2}(0)=\til G_{I,2,t_2}(0,x_1)=\til g_{I,2}(t_2,x_1)$.
Thus, after a time-change, $g_{I,2}(t_2,\beta_{1}(t_1))$, $0\le t_1\le  T_1(J_1)$, is
a partial annulus SLE$(\kappa,\Lambda_1)$ trace in $\A_{p-t_2}$ started from $g_{I,2}(t_2,a_{1})$ with marked point $I_{p-t_2}\circ e^i(\xi_2(t_2))$. This means that, conditioning on $\F^2_{t_2}$, after a time-change, $\beta_1(t_1)$, $0\le t_1\le T_1(J_1)$, is a partial annulus SLE$(\kappa,\Lambda_1)$ trace in $\A_p\sem \beta_{I,2}((0,t_2])$ started from $a_1$ with marked point $\beta_{I,2}(t_2)$.
Similarly, the above statement holds true if the subscripts ``$1$'' and ``$2$'' are exchanged.

The joint distribution $\nu_{J_1,J_2}$ is a local coupling such that the desired properties in the statement of Theorem \ref{coupling-thm} holds true up to the stopping time $T_1(J_1)$ and $T_2(J_2)$. Then we can apply the coupling technique developed in Section 7 of \cite{reversibility} to construct a global coupling using the local couplings for different pairs $(J_1,J_2)$.

The coupling technique is composed of several steps. First, let $\{(J_1^k,J_2^k):k\in\N\}$ denote the set of all pairs in $\JP$ such that $J_j^k$, $k\in\N$, $j=1,2$, are polygonal curves, whose vertices have rational coordinates. Second, for every $n\in\N$, one may find a coupling of $\beta_1$ and $\beta_2$ such that, for every $1\le k\le n$, if $\beta_1$ is stopped at $\tau_{J_1^k}$, and $\beta_2$ is stopped at $\tau_{J_2^k}$, then the joint distribution is $\nu_{J_1^k,J_2^k}$. To construct such coupling, we work on the two-dimensional random process $M$. One may prove that there is a process $M_n$ defined on $[0,p]^2$, which satisfies the following properties:
 \begin{enumerate}
   \item $M_n$ is a martingale in one variable, when the other variable is fixed;
\item $M_n=1$ when either variable is $0$;
\item $M_n=M$ on $[0,\tau_{J_1^k}]\times [0,\tau_{J_2^k}]$, $1\le k\le n$.
 \end{enumerate}
The reader is referred to Theorem 6.1 in \cite{reversibility} for the construction of $M_n$. The $\nu_n$ is then defined by $d\nu_n/d\mu=M_n(p,p)$. Finally, the global coupling measure $\nu$ is any subsequential weal limit of the sequence $(\nu_n)$ in some suitable topology.

\subsection{Other results}
Here we state without proofs some other results which can be proved using the idea in the proof of Theorem \ref{coupling-thm}.

\begin{Theorem}
Let $\kappa>0$. Suppose $\Gamma$ is a $C^{1,2}$ differentiable function on $(0,\infty)\times(\R\sem 2\pi\Z)$ that satisfies
 \BGE {\pa_t{ \Gamma}}=\frac{\kappa}2  { \Gamma''}
+\HA{ \Gamma'}+\Big(\frac 3\kappa-\frac 12\Big)\HA' \Gamma.\label{dot-Gamma=2}\EDE
Let $\Lambda=\kappa\frac{\Gamma'}{\Gamma}$, $\Lambda_1=\Lambda$, and $\Lambda_2$ be the dual of $\Lambda$.
 Then for any $p>0$ and $a_1\ne a_2\in\TT$, there is a coupling of two curves:
 $\beta_1(t)$, $0\le t<T_1$, and $\beta_{2}(t)$, $0\le t<T_2$, such that  for $j\ne k\in\{1,2\}$ the following hold.
\begin{enumerate}[(i)]
\item    $\beta_j $  is an annulus SLE$(\kappa,\Lambda_j)$ trace in $\A_p$
 started from $a_j$ with marked point $a_k$.
\item  If $t_k\in[0,T_k)$ is a stopping time w.r.t.\ $(K_k(t))$, then
conditioned on $\beta_{k}(t)$, $0\le t\le t_k$,  after a time-change, $\beta_j(t)$, $0\le t<T_j(t_k)$, is a partial annulus
SLE$(\kappa,\Lambda_j)$ process in a component of $\A_p\sem \beta_k((0,t_k])$ started from
$a_j$ with marked point $\beta_k(t_k)$, where $T_j(t_k)$ is the first time that $\beta_j$ hits $\beta_k([0,t_k])$, and is set to be $T_j$ if such time does not exist. If $\kappa\in(0,4]$, the word ``partial'' could be removed.
\end{enumerate}
\label{coupling-thm-2}
\end{Theorem}

\no{\bf Remark.} \begin{enumerate}
  \item The $\Lambda$ in the theorem satisfies the following partial differential equation:
\BGE \pa_t \Lambda=\frac\kappa 2 \Lambda''+\Big(3-\frac\kappa 2\Big)\HA''+\Lambda\HA'+\HA\Lambda'+\Lambda\Lambda'.\label{dot-Lambda=2}\EDE
On the other hand, if $\Lambda$ satisfies (\ref{dot-Lambda=2}), then there is $\Gamma$, which satisfies $\Lambda=\kappa\frac{\Gamma'}{\Gamma}$ and (\ref{dot-Gamma=2}).
\item Theorem \ref{coupling-thm-2} also holds for $\kappa=0$ if $\Lambda$ solves (\ref{dot-Gamma=2}).
\item We may also derive similar results for radial SLE$(\kappa,\Lambda)$ process and strip SLE$(\kappa,\Lambda)$ process. In these two cases, $\Gamma$ and $\Lambda$ are functions of a single variable, and $\Lambda=\kappa\frac{\Gamma'}{\Gamma}$. If $\Lambda=\frac\rho2\cot_2$ or $\Lambda=\frac\rho2\coth_2$, respectively, in these two cases, then we get the radial SLE$(\kappa,\rho)$ and strip SLE$(\kappa,\rho)$ processes, respectively (c.f.\ \cite{SW}). For the radial SLE$(\kappa,\Lambda)$ process, to have the commutation relation, we need that $\Gamma$ solves the ODE
    \BGE 0=\frac{\kappa}2  { \Gamma''}+\cot_2{ \Gamma'}+\Big(\frac 3\kappa-\frac 12\Big)\cot_2' \Gamma+C\Gamma,\label{radial}\EDE
    where $C$ is a constant. For the strip SLE$(\kappa,\Lambda)$ process, $\Gamma$ must solves (\ref{radial}) with $\cot_2$ replaced by $\coth_2$ to guarantee the existence of the commutation coupling.
\end{enumerate}

\section{Coupling of Whole-Plane SLE}\label{coupling-deg}
The goal of this section is to prove Theorem \ref{coupling-thm-int-skew} below, which is about commutation couplings of two whole-plane SLE processes. This result will later be used to prove the whole-plane reversibility. Since the proof is similar to the proof of Theorem \ref{coupling-thm}, we will frequently quote the arguments in the previous section.

\begin{Theorem} Let $\kappa>0$ and $s_0\in\R$. Suppose $\Gamma$ is a positive $C^{1,2}$ differentiable function on $(0,\infty)\times\R$ that satisfies (\ref{dot-Gamma=}) and (\ref{Gamma-period-s}). Let $\Lambda=\kappa\frac{\Gamma'}{\Gamma}$, $\Lambda_1=\Lambda$, and $\Lambda_2$ be the dual of $\Lambda_1$.
Let $s_1=s_0$ and $s_2=-s_0$. Then there is a coupling of two curves $\beta_{I,1}(t)$, $-\infty<t<\infty$, and $\beta_{I,2}(t)$, $-\infty<t<\infty$,
such that for $j\ne k\in\{1,2\}$, the following hold.
\begin{enumerate}[(i)]
\item $\beta_{I,j}$ is a whole-plane SLE$(\kappa,s_j)$ trace in $\ha\C$ from $0$ to $\infty$;
\item Let $t_k$ be a finite stopping time w.r.t.\ $(K_{I,k}(t))$.
Then conditioned on $\beta_{I,k}(s)$, $-\infty<s\le  t_k$, after a time-change, the curve $\beta_{I,j}(t_j)$, $-\infty<t_j<T_j( t_k)$, is a disc SLE$(\kappa,\Lambda_j)$ process in a component of $\ha\C\sem I_0(\beta_{I,j}([-\infty,t_j]))$ started from $0$ with marked point $I_0(\beta_{I,j}(t_j))$, where $T_j(t_k)$ is the first time that $\beta_j$ hits $\beta_k([-\infty,t_k])$, or $\infty$ if such time does not exist.
\end{enumerate}
\label{coupling-thm-int-skew}
\end{Theorem}

\subsection{Estimations on Loewner maps}
Let $\til g(t,\cdot)$, $t\in\R$, be the inverted covering {\it whole-plane} Loewner maps driven by some $\xi\in C(\R)$. Let $z\in\C$ and $h(t)=\Imm \til g(t,z)>0$ for $t\in (-\infty,\tau_z)$, the interval on which $\til g(t,z)$ is defined.
From (\ref{radial-eqn-covering}) we have $-\tanh_2(h(t))\ge h'(t)\ge -\coth_2(h(t))$, and
\BGE |\pa_t\til g(t,z)+i|\le \frac{2}{e^{\Imm \til g(t,z)}-1}=\frac{2}{e^{h(t)}-1},\quad t\in(-\infty,\tau_z).\label{|pat-til-whole|}\EDE
So $h(t)$ decreases, and $\frac d{dt} \ln(\cosh_2(h(t)))\ge -1/2$, which
together with (\ref{invert-covering-whole-limit}) and integration implies that $ \cosh_2(h(t))\ge \frac 12 e^{\frac{\Imm z}2-\frac{t}2}$. Then we have $e^{h(t)}\ge e^{\Imm z-t}-3$. From (\ref{|pat-til-whole|}) we see that, if $t<\Imm z-\ln(8)$, then  $|\pa_t\til g(t,z)+i|\le \frac{2}{e^{\Imm z-t}-4}\le 4e^{t-\Imm z}$. From (\ref{invert-covering-whole-limit}) and integration we have
\BGE |\til g(t,z)+it-z|\le 4 e^{t-\Imm z}\le 1/2,\quad\mbox{if } t\le \Imm z-\ln(8).\label{invert-whole>=}\EDE
If $\til g(t,\cdot)$ are the covering whole-plane Loewner maps, then from $\til g(t,\cdot)=\til I_0\circ \til g_I(t,\cdot)\circ \til I_0$, we have
\BGE |\til g_I(t,z)-it-z|\le 4 e^{t+\Imm z}\le 1/2,\quad\mbox{if } t\le -\Imm z- \ln(8).\label{whole>=}\EDE

 Let $\til g_I(t,\cdot)$, $-\infty<t<0$, be the covering {\it disc} Loewner maps driven by some $\xi\in(-\infty,0)$.
Let $z\in\HH$ and $h(t)=\Imm \til g_I(t,z)>0$ for $t\in (-\infty,\tau_z)$. From Lemma \ref{estimation} and  (\ref{covering-disc-diff})  we see that if $-t\ge h(t)+2$ then $|h'(t)|\le 5.5e^{h(t)+t}$, so $|\frac d{dt} e^{-h(t)}|\le 5.5 e^{t}$. From (\ref{covering-disc-end}) we see that $-t\ge h(t)+2$ when $t$ is close to $-\infty$. Suppose that $-t\ge h(t)+2$ does not hold for all $t\in(-\infty,\tau_z)$, and let $t_0$ be the first $t$ such that $-t=h(t)+2$. Then $|\frac d{dt} e^{-h(t)}|\le 5.5 e^{t}$ on $(-\infty,t_0]$. From (\ref{covering-disc-end}) and integration we have $e^{t_0+2}=e^{-h(t_0)}\ge e^{-\Imm z}-5.5e^{t_0}$, which implies that $e^{-\Imm z}\le (e^2+5.5)e^{t_0}<13e^{t_0}$. Thus, if $t\le -\Imm z-\ln(13)$ then  $-t\ge h(t)+2$, so  $|h'(t)|\le 5.5e^{h(t)+t}$. From (\ref{covering-disc-end}) and integration, we see that, if $t\le -\Imm z-\ln(13)$, then  $e^{-\Imm \til g_I(t,z)} \ge \frac {7.5}{13}e^{-\Imm z}$,
which implies that
$\Imm \til g_I(t,z)\le \Imm z+\ln(13/7.5)<-t-2$,
which, together with Lemma \ref{estimation}, implies that
$|\HA_I(-t,\til g_I(t,z)-\xi(t))| \le 5.5\frac{13}{7.5} e^{\Imm z+t}<10e^{\Imm z+t}$.
From (\ref{covering-disc-diff}), (\ref{covering-disc-end}) and integration we have
$ |\til g_I(t,z)-z|\le 10 e^{\Imm z+t}$, if $t\le -\Imm z- \ln(13)$.
If $\til g(t,\cdot)$ are the inverted covering disc Loewner maps, then from $\til g(t,\cdot)=\til I_{-t}\circ \til g_I(t,\cdot)\circ \til I_0$, we have
\BGE |\til g(t,z)+it-z|\le 10 e^{-\Imm z+t}\le 10/13,\quad \mbox{if } t\le\Imm z- \ln(13).\label{invert-disc>=}\EDE

\subsection{Ensemble} \label{ensemble-2}
The argument in this subsection is parallel to that in Section \ref{ensemble-1}.
Let $\xi_1,\xi_2\in C(\R)$. For $j=1,2$, let $g_{I,j}(t,\cdot)$ (resp.\ $g_{j}(t,\cdot)$), $t\in\R$, be the whole-plane (resp.\ inverted whole-plane) Loewner  maps driven by $\xi_j$. Let  $\til g_{I,j}(t,\cdot)$ and $\til g_{j}(t,\cdot)$, $t\in\R$, $j=1,2$,
 be the corresponding  covering Loewner  maps. Suppose $\xi_j$ generates a simple whole-plane Loewner trace: $\beta_{I,j}$, $j=1,2$. Let $\beta_{I,j}=I_0\circ \beta_j$, $j=1,2$, be the inverted trace. Let $K_j(t)$ and $K_{I,j}(t)$ be the corresponding hulls.
Define $\cal D$ and $\mA$ using (\ref{cal-D}) and (\ref{m}) with $0$ replaced by $-\infty$ and $\A_p$ replaced by $\ha\C$.
Fix  any $j\ne k\in\{1,2\}$ and $t_k\in\R$. Let $T_j(t_k)$ be as defined as before.
Then for any $t_j<T_j(t_k)$, we have $(t_1,t_2)\in\cal D$. Moreover, as $t_j\to T_j(t_k)^-$, $\mA(t_1,t_2)\to 0$.

For $-\infty\le t_j<T_j(t_k)$, let $\beta_{I,j,t_k}(t_j)= g_{k}(t_k,\beta_{I,j}(t_j))$. Then $\beta_{j,t_k}$ is a simple curve in $\D$ starts from $0$. For $-\infty<t_j<T_j(t_k)$, let
$v_{j,t_k}(t_j)=-\modd(\D\sem \beta_{I,j, t_k}([-\infty,t_j]))=-\mA(t_1,t_2)$. 
Then $v_{j, t_k}$ is continuous and increasing and maps $(-\infty,T_j(t_k))$ onto $(-\infty,0)$.
 Let $\gamma_{I,j,t_k}(t)=\beta_{I,j,t_k}(v_{j,t_k}^{-1}(t))$, $-\infty\le t<0$. Then $\gamma_{I,j,t_k}$ is the disc Loewner trace driven by some $\zeta_{j,t_k}\in C((-\infty,0))$. Let $\gamma_{j,t_k}$ be the corresponding inverted disc Loewner trace. Let  $h_{I,j,t_k}(t,\cdot)$ and
$h_{j,t_k}(t,\cdot)$  be the corresponding disc and inverted disc Loewner maps. Let
$\til h_{I,j,t_k}(t,\cdot)$ and $\til h_{j,t_k}(t,\cdot)$ be the corresponding covering Loewner maps.

 For $-\infty<t_j<T_j(t_k)$, let $\xi_{j,t_k}(t_j)$, $\beta_{j,t_k}(t_j)$, $g_{I,j,t_k}(t_j,\cdot)$, $g_{j,t_k}(t_j,\cdot)$, $\til g_{I,j,t_k}(t_j,\cdot)$, and $\til g_{j,t_k}(t_j,\cdot)$ be the time-changes of $\zeta_{j,t_k}(t)$, $\gamma_{j,t_k}(t)$, $h_{I,j,t_k}(t,\cdot)$,
$ h_{j,t_k}(t,\cdot)$,  $ \til h_{I,j,t_k}(t,\cdot)$, and $ \til h_{j,t_k}(t,\cdot)$, respectively, via $v_{j,t_k}$.

Define $G_{I,k,t_k}(t_j,\cdot)$ and $\til G_{I,k,t_k}(t_j,\cdot)$ by (\ref{G}) and (\ref{til-G}). Then we could choose the driving function $\zeta_{j,t_k}$ such that (\ref{eta=W(xi)}) holds. Define $A_{j,h}$ and $A_{j,S}$ using (\ref{A}) and (\ref{AS}). A standard argument using Lemma 2.1 in \cite{Zhan} shows (\ref{v'}) and (\ref{pa1m}) hold here. From the definition of $\til G_{I,k,t_k}(t_j,\cdot)$, we get (\ref{circ=1}), which can be differentiated to conclude that (\ref{dot=-3''W}) holds here, and (\ref{dot-W'-1}) holds with $p-t_j$ replaced by $\infty$.
Let $X_j$  be defined by (\ref{Xj}). Then (\ref{Xjk}) holds. Let $Q$ be defined by (\ref{Q}). Then (\ref{dot-til-g-1-diff})-(\ref{pa-AS}) still hold.

From Lemma \ref{estimation},  we have
\BGE Q=O(e^{-\mA}),\quad \mbox{as }\mA\to \infty.\label{Q-int-bound}\EDE
From Lemma \ref{modulus-bound}, we see that, for $j=1,2$,
\BGE \ln(A_{j,1}),A_{j,S}=O(\mA e^{-\mA}),
\quad \mbox{as } \mA\to \infty.\label{A-j-bound}\EDE

Since $e^{t_j}$ is the capacity of $K_{I,j}(t_j)$, which contains $0$, we have $K_{I,j}(t_j)\subset\{|z|\le  4e^{t_j}\}$. This then implies that $K_j(t_j)\subset \{|z|\ge e^{-t_j}/4\}$, $\til K_{I,j}(t_j)\subset\{\Imm z\ge -t_j-\ln(4)\}$, and $\til K_j(t_j)\subset\{\Imm z\le \ln(4)+t_j\}$. Thus,
\BGE \{(t_1,t_2)\in\R^2:t_1+t_2 <-\ln(16)\}\subset \cal D,\label{Dsupset}\EDE
\BGE \mA(t_1,t_2)\ge -t_1-t_2-\ln(16),\quad \mbox{if }(t_1,t_2)\in\cal D.\label{ma>}\EDE
From (\ref{Q-int-bound})-(\ref{ma>}), we see that
$ A_{1,1}^2  A_{2,1}^2  Q=O(e^{t_1+t_2})$ as $t_1+t_2\to -\infty$. 
Define $F $ on $\cal D$  using (\ref{F}) with the lower bounds $0$ replaced by $-\infty$.
From Lemma \ref{modulus-bound},  $A_{k,S}\to 0$ as $t_j\to -\infty$.
Thus, (\ref{paF}) still holds here, and $\ln (F(t_1,t_2))=\int_{-\infty}^{t_1}\frac{A_{1,S}(s_1,t_2)}{A_{1,1}(s_1,t_2)^2}
\cdot A_{1,1}(s_1,t_2)^2ds_1$. Changing variable with $x(s_1)=\mA(s_1,t_2)$, and using (\ref{pa1m}) and (\ref{A-j-bound}),we conclude that
\BGE \ln(F)=O(\mA e^{-\mA}),\quad\mbox{as } \mA\to \infty.\label{F-bound}\EDE

\subsection{Martingales in two time variables} \label{martingale-2}
 The argument in this subsection is parallel to that in Section \ref{martingale}.
Let $(B^{(\kappa)}_1(t))$ and $(B^{(\kappa)}_2(t))$  be two independent pre-$(\TT;\kappa)$-Brownian motion.
Let $\xi_j(t)=B^{(\kappa)}_j(t)+s_jt$, $t\in\R$, $j=1,2$. For simplicity, suppose $\kappa\in(0,4]$. Then for $j=1,2$, a.s.\ $(\xi_j)$ generates a simple whole-plane Loewner trace $\beta_{I,j}$, which is a whole-plane SLE$(\kappa,s_0)$ trace in $\ha \C$ from $0$ to $\infty$. We may apply the results in the prior subsection. For $j=1,2$, let $(\F^j_t)_{t\in\R}$ be the complete filtration generated by $e^i(\xi_j)$.
The whole-plane Loewner objects driven by $\xi_j$ are all $(\F^j_t)$-adapted, because they are all determined by $(e^i(\xi_j(t)))$. It is easy to check that for $j\ne k\in\{1,2\}$, the processes $(\beta_{I,j,t_k})$, $(\til g_{I,j,t_k}(t_j,\cdot))$, $(A_{j,h})$, $h=1,2,3$, $(A_{j,S})$, $(G_{I,j,t_j}(t_k,\cdot))$, $(\til G_{I,j,t_j}(t_k,\cdot))$, $(e^i(\xi_{j,t_k}))$, $(e^i(X_j))$, $(\mA)$, $(\HA_I^{(h)}(\mA,X_1))$,  $(\Gamma_{j}(\mA,X_j))$,  $(\Lambda_{j}(\mA,X_j))$, $(Q)$ and $(F)$ defined on $\cal D$ are all $(\F^1_{t_1}\times \F^2_{t_2})$-adapted. This is not true for $(\xi_{j,t_k}(t_j))$ and $(X_j)$, but is true for their images under the map $e^i$. Define $Y$ using (\ref{Y}). Then $(Y)$ is also $(\F^1_{t_1}\times \F^2_{t_2})$-adapted since $\Gamma_j$ has period $2\pi$.

In this section we work on SDE with the meaning as in Definition \ref{SDE-T}: the stochastic part contains pre-$(\TT;\kappa)$-Brownian motions, and the time intervals start from $-\infty$. The traditional It\^o's formula works only for time intervals that start from $0$ or a finite number. To derive the results in this section, we may truncate the interval ``$(-\infty,T)$'' by an arbitrary real number $c$ (and we work on the interval $[c,T)$), which is close to $-\infty$. Fix $j\ne k\in\{1,2\}$ and any $(\F^k_t)$-stopping time $t_k\in\R$.
Let $\F^{j,t_k}_{t_j}=\F^j_{t_j}\times \F^k_{t_k}$, $t_j\in\R$. From now on,  all SDE will be $(\F^{j,t_k}_{t_j})$-adapted (with the meaning as in Definition \ref{SDE-T}), and $t_j$ ranges in $[0,T_j(t_k))$.

First, we find that (\ref{dX}) still holds here, which then implies (\ref{dY}). From the modified (\ref{dot-W'-1}),  we see that (\ref{pa-A-j-alpha}) holds here with $p-t_j$ replaced by $\infty$. Let $\ha M$ be defined by (\ref{haM}). Then (\ref{pahaM}) holds with $p-t_j$ replaced by $\infty$.

Define $M$ on $\cal D$ by
\BGE M=\ha M  \exp\Big(\alpha\rA(\infty) (\mA+t_1+t_2)+\sum_{j=1,2}-\frac{s_j}\kappa\xi_j(t_j)+\frac{s_j^2}{2\kappa}t_j\Big).\label{M-int}\EDE
Then $M$ is $(\F^1_{t_1}\times \F^2_{t_2})$-adapted.
From the modified (\ref{pahaM}) and that $\xi_j(t_j)=B^{(\kappa)}_j(t_j)+s_jt_j$, we compute
\BGE \frac{\pa_j M}{ M}=\Big[\Big(3-\frac\kappa 2\Big)\frac{A_{j,2}}{A_{j,1}}+
{A_{j,1}} \Lambda_{j}(\mA,X_j)-s_j\Big]\frac{\pa B^{(\kappa)}_j(t_j)}\kappa.\label{paM-int}\EDE
So $M$ is a local martingale in $t_j$ when $t_k$ is a finite stopping time.

 \vskip 3mm

Let ${\cal J}$ denote the set of Jordan curves in $\C\sem\{0\}$ that surround $0$.
For $J\in\cal J$ and $j=1,2$, let $T_j(J)$ denote the first time that $\beta_j$ hits $J$.
Then $T_j(J)$ is also the first time that $\beta_{I,j}$ hits $I_0(J)$.
 Let $H_J$ denote the closure of the domain bounded by $I_0(J)$, and let $C_J$ denote the capacity of $H_J$. If $K_{I,j}(t)\subset H_J$, then $e^t\le C_J$. So we have $T_j(J)\le \ln(C_J)$.

 Let $\JP$ denote the set of pairs $(J_1,J_2)\in{\cal J}^2$ such that $I_0(J_1)\cap J_2=\emptyset$ and $I_0(J_1)$ is surrounded by  $J_2$.
This is equivalent to that $I_0(J_2)\cap J_1=\emptyset$ and  $I_0(J_2)$ is surrounded by  $J_1$.
Then for every $(J_1,J_2)\in\JP$, $\beta_{I,1}(t_1)\ne\beta_{2}(t_2)$ when $t_1\le T_{1}(J_1)$ and $t_2\le T_2(J_2)$, so
$(-\infty,T_{1}(J_1)]\times(-\infty,T_{2}(J_2)]\subset \cal D$.

\begin{Proposition}  (Boundedness) Fix $(J_1,J_2)\in\JP$. (i) $|\ln(M)|$ is bounded on $(-\infty,T_{1}(J_1)]\times(-\infty,T_2(J_2)]$
by a constant depending only on $J_1$ and $J_2$. (ii) Fix any $j\ne k\in\{1,2\}$. Then $\ln(M)\to 0$ as $t_j\to -\infty$
uniformly in $t_k\in(-\infty,T_k(J_k)]$.
\label{bounded-int} \end{Proposition}
{\bf Proof.} Let $\Gamma_{s_0}$ be given by Lemma \ref{Gamma-uniform}. Let $\Gamma_{s_0,1}=\Gamma_{s_0}$, and $\Gamma_{s_0,2}(t,x)=\Gamma_{s_0}(t,-x)$. Define $Y_{s_0}$ on $\cal D$ by $Y_{s_0}=\Gamma_{s_0,1}(\mA,X_1)=\Gamma_{s_0,2}(\mA,X_2)$. 
Then $Y_{s_0}=Y\exp(-\frac {s_0}\kappa X_1-\frac{s_0^2\mA}{2\kappa})$.
From Lemma \ref{Gamma-uniform}, \BGE \ln(Y_{s_0})=o(\mA )\quad \mbox{as }\mA \to \infty.\label{Y-bound}\EDE
 Define $\ha M_{s_0}$ using (\ref{haM}) with $Y$ replaced by $Y_{s_0}$. From (\ref{M-int}) we have
 \BGE M=M_{s_0}\exp\Big((\alpha\rA(\infty)+\frac{s_0^2}{2\kappa})(\mA+t_1+t_2)+\frac{s_0}\kappa (X_1-\xi_1(t_1)+\xi_2(t_2))\Big).
 \label{M-Ms}\EDE

 From (\ref{A-j-bound}), (\ref{F-bound}), (\ref{Y-bound}), and
that $\RA(p)=O(e^{-p})$ as $p\to\infty$,  we see that there is a positive continuous function $f$ on $(0,\infty)$
with $\lim_{x\to\infty} f(x)= 0$ such that
\BGE|\ln(M_{s_0}(t_1,t_2))|\le f(\mA(t_1,t_2)).\label{M-f}\EDE

 Let $\Omega(I_0(J_1),J_2)$ denote the doubly connected domain bounded by $I_0(J_1)$ and $J_2$. Let $p_0>0$ denote its modulus.
 For $(t_1,t_2)\in (-\infty,T_{1}(J_1)]\times(-\infty,T_{2}(J_2)]$, since $\Omega(I_0(J_1),J_2)$ disconnects $K_{I,1}(t_1)$ from $K_2(t_2)$,  we have $\mA(t_1,t_2)\ge p_0$. On the other hand, $\mA\le p$. From (\ref{M-f}) we see that $\ln(M_{s_0})$ is uniformly bounded.
From (\ref{Dsupset}),  (\ref{ma>}),  (\ref{M-f}), and that $T_k(J_k)\le C_{J_k}<\infty$, we see that  $\ln(M)\to 0$ as $t_j\to -\infty$
uniformly in $t_k\in(-\infty,T_k(J_k)]$. The rest of the proof follows from (\ref{M-Ms}) and the following proposition.  $\Box$

\begin{Proposition} Fix $(J_1,J_2)\in\JP$. (i) $|X_1-\xi_1+\xi_2|$ and $|\mA+t_1+t_2|$ are bounded on $(-\infty,T_1(J_1)]\times (-\infty,T_2(J_2)]$ by constants depending only on $J_1$ and $J_2$. (ii) For any $j\ne k\in\{1,2\}$, $X_1-\xi_1+\xi_2\to 0$ and $\mA+t_1+t_2\to 0$ as $t_j\to -\infty$, uniformly in  $t_k\in(-\infty,T_k(J_k)]$.
\label{bounded-skew1}
\end{Proposition}
{\bf Proof.} Recall that $T_j(J_j)\le C_{J_j}<\infty$ for $j=1,2$, and $\mA\ge p_0>0$ on $(-\infty,T_1(J_1)]\times (-\infty,T_2(J_2)]$. If there is no ambiguity, let $\Omega(A,B)$ denote the domain bounded by sets $A$ and $B$, and let $\modd(A,B)$ denote the modulus of this domain if it is doubly connected.

From (\ref{Xj}) we have $X_1(t_1,t_2)=\til G_{I,2,t_2}(t_1,\xi_1(t_1))-\til g_{I,1,t_2}(t_1,\xi_2(t_2))$. So
\BGE |X_1(t_1,t_2)-\xi_1(t_1)+\xi_2(t_2)|\le |\til g_{I,1,t_2}(t_1,\xi_2(t_2))-\xi_2(t_2)|+|\til G_{I,2,t_2}(t_1,\xi_1(t_1))-\xi_1(t_1)|.\label{X-+<}\EDE
From (\ref{covering-disc-end}) we have $\lim_{t_1\to -\infty} \til g_{I,1,t_2}(t_1,\xi_2(t_2))=\xi_2(t_2)$. From (\ref{pa1m}), (\ref{dot-til-g-1-I'}), and Lemma \ref{estimation}, we see that there is a deterministic positive decreasing function $f(x)$ with $\lim_{x\to\infty}f(x)= 0$  such that $|\til g_{I,1,t_2}(t_1,\xi_2(t_2))-\xi_2(t_2)|\le f(\mA(t_1,t_2))$. Since $\mA\ge p_0$ on $(-\infty,T_{1}(J_1)]\times(-\infty,T_2(J_2)]$, $|\til g_{I,1,t_2}(t_1,\xi_2(t_2))-\xi_2(t_2)|$ is uniformly bounded by $f(p_0)$. From (\ref{ma>}) and that $T_2(J_2)\le C_{J_2}$, we see that $\til g_{I,1,t_2}(t_1,\xi_2(t_2))-\xi_2(t_2)\to 0$ as $t_1\to-\infty$, uniformly in $t_2\in(-\infty,T_2(J_2)]$.

Let $J$ be a Jordan curve separating $J_1$ and $J_{I,2}$. Let $p_1=\modd(J,J_1)$ and $p_2=\modd(J,J_{I,2})$. Let $\til J=(e^i)^{-1}(J)$. Let $h_m=\inf \{\Imm z:z\in\til J\}$ and $h_M=\sup\{\Imm z:z\in\til J\}$. Then both $h_m$ and $h_M$ are finite. For $j=1,2$, there is $h_j>0$ depending only on $p_j$, such that, if $K\subset \D$ is an interior hull with $0\in K$ and $\modd(\D\sem K)\ge p_j$, then $K\subset\{|z|\le e^{-h_j}\}$. If $t_1\le T_1(J_1)$, then $J_1$ disconnects $J$ from $K_1(t_1)$, so $\modd(J,K_1(t_1))\ge p_1$. Since $\Omega(J,K_1(t_1))$ is mapped by $g_1(t_1,\cdot)$ conformally onto $\Omega(g_1(t_1,J),\TT)\subset\D$, $\modd(g_1(t_1,J),\TT)\ge p_1$. Since $g_1(t_1,J)$ surrounds $0$, $g_1(t_1,J)\subset\{|z|\le e^{-h_1}\}$. Since $\til g_1(t_1,\til J)=(e^i)^{-1}(g_1(t_1,J))$, $\til g_1(t_1,\til J)\subset\{\Imm z\ge h_1\}$. Similarly, if $t_2\le T_2(J_2)$, then $\til g_{I,2}(t_2,\til J)\subset\{\Imm z\le -h_2\}$. If $t_1\le T_1(J_1)$ and $t_2\le T_2(J_2)$, then $g_{1,t_2}(t_1,\cdot)\circ  g_{I,2}(t_2,\cdot)$ maps $\C\sem K_1(t_1)\sem K_{I,2}(t_2)$ conformally onto $\A_{\mA}$. A similar argument shows that the image of $J$ under this map lies in $\{e^{-\mA+h_2}\le |z|\le e^{-h_1}\}$. Thus, $\til g_{1,t_2}(t_1,\til g_{I,2}(t_2,\til J))\subset\{h_1\le \Imm z\le \mA-h_2\}$, if $t_1\le T_1(J_1)$ and $t_2\le T_2(J_2)$.

Let $z_0\in \C\sem \til K_1(t_1)\sem \til K_{I,2}(t_2)$, $w_1=\til g_1(t_1,z_0)$, $w_2=\til g_{I,2}(t_2,z_0)$, and $w_3=\til g_{1,t_2}(t_1,w_2)$. From (\ref{invert-whole>=}), (\ref{whole>=}), and (\ref{invert-disc>=}) we see that
\BGE |w_1-(z_0-it_1)|\le 4e^{t_1-\Imm z_0}\le 1/2,\quad \mbox{if } t_1\le \Imm z_0-\ln(8);\label{>1}\EDE
\BGE |w_2-(z_0+it_2)|\le 4e^{t_2+\Imm z_0}\le 1/2,\quad \mbox{if } t_2\le -\Imm z_0-\ln(8);\label{>2}\EDE
\BGE |w_3-(w_2+i\mA)|\le 10e^{-\mA-\Imm w_2}<1,\quad \mbox{if } \Imm w_2+\mA\ge \ln(13).\label{>3}\EDE
Now let $z_0\in \til J$. From the prior paragraph, $\Imm \til g_1(s,z_0)\ge h_1$ for $s\le t_1$, $\Imm\til g_{I,2}(s,z_0)\le -h_2$ for $s\le t_2$, and $\mA(s,t_2)-h_2\ge \til g_{1,t_2}(s,w_2)\ge h_1$ for $s\le t_1$. From (\ref{|pat-til-whole|})  we have $|\pa_t\til g_1(s,z_0)+i|\le \frac{2}{e^{h_1}-1}$, for $s\le t_1$. Similarly,  $|\pa_t\til g_{I,2}(s,z_0)-i|\le \frac{2}{e^{h_2}-1}$ for $s\le t_2$.
If $t_1\le \Imm z_0-\ln(8)$, then from (\ref{>1}), $|w_1-(z_0-it_1)|\le 1/2$. If $t_1>\Imm z_0-\ln(8)$, we let $t_1'=\Imm z_0-\ln(8)$, and $w_1'=\til g_1(t_1',z_0)$. Then we have $|w_1'-(z_0-it_1')|\le 1/2$. From the bound of $|\pa_t\til g_1(s,z_0)+i|$, we see that $$|(w_1+it_1)-(w_1'+it_1')|\le \frac{2(t_1-t_1')}{e^{h_1}-1}\le \frac{2C_{J_1}-2(\Imm z_0-\ln(8))}{e^{h_1}-1}\le \frac{2C_{J_1}+2\ln(8)-2h_m}{e^{h_1}-1}. $$ Let $A_1=\frac 12 +\max \{0,\frac{2C_{J_1}+2\ln(8)-2h_m}{e^{h_1}-1} \}$. Then in all cases we have \BGE |w_1-(z_0-it_1)|\le A_1.\label{w1-z0}\EDE
Similarly, let $A_2=\frac 12 +\max \{0,\frac{2C_{J_2}+2\ln(8)+2h_M}{e^{h_2}-1} \}$. Then we always have
\BGE |w_2-(z_0+it_2)|\le A_2.\label{w2-z0}\EDE
Since $\Imm z_0\ge h_m$, we have
\BGE t_2-\Imm w_2\le A_2-\Imm z_0\le A_2-h_m.\label{w2>=}\EDE

If $\Imm w_2+\mA(t_1,t_2)\ge \ln(13)$, from (\ref{>3}), we have $|w_3-(w_2+i\mA(t_1,t_2))|<1$. Now suppose that $\Imm w_2+\mA(t_1,t_2)< \ln(13)$. We may choose $t_1'<t_1$ such that $\Imm w_2+\mA(t_1',t_2)=\ln(13)$. Let $w_3'=\til g_{1,t_2}(t_1',w_2)$. Then we have $|w_3'-(w_2+i\mA(t_1',t_2))|<1$. For $s\le t_1$, since $h_1\le \Imm \til g_{1,t_2}(s,w_2)\le \mA(s,t_2)-h_2$, from Lemma \ref{estimation} we have $$|\HA_I(\mA(s,t_2),i\mA(s,t_2)-\til g_{1,t_2}(s,w_2)+\xi_{1,t_2}(s))|\le \frac{4e^{-h_1}}{(1-e^{-h_1})^3}.$$ Since $\HA_I(\mA,z)+i=\HA_I(\mA,z-i\mA)=-\HA_I(\mA,i\mA-z)$, we have
$$|\HA(\mA(s,t_2),\til g_{1,t_2}(s,w_2)-\xi_{1,t_2}(s))+i|\le \frac{4e^{-h_1}}{(1-e^{-h_1})^3},\quad\mbox{if } s\le t_1.$$
Let $C_0=\frac{4e^{-h_1}}{(1-e^{-h_1})^3}$. From (\ref{pa1m}), (\ref{dot-til-g-1-I'}), (\ref{ma>}), (\ref{w2>=}), and the above inequality, we have $$|(w_3-i\mA(t_1,t_2))-(w_3'-i\mA(t_1',t_2))|\le C_0(\mA(t_1',t_2)-\mA(t_1,t_2))$$ $$\le
 C_0(\ln(13)-\Imm w_2+t_1+t_2+\ln(16))\le C_0(\ln(13)+\ln(16)+C_{J_1}+A_2-h_m). $$
Let $A_3=1+\max\{0,C_0(\ln(13)+\ln(16)+C_{J_1}+A_2-h_m)\}$. Then
$ |w_3-(w_2+i\mA)|\le A_3 $ 
always holds, which together with
(\ref{w1-z0}) and (\ref{w2-z0}) implies that, for any $t_1\le T_1(J_1)$ and $t_2\le T_2(J_2)$,  \BGE |\til G_{I,2,t_2}(t_1,w_1)-w_1-i(\mA+t_1+t_2)|\le A_1+A_2+A_3,\quad w_1\in\til g_1(t_1,\til J).\label{<1-bound}\EDE
Now $\til g_1(t_1,\til J)$ is a curve with period $2\pi$ above $\R$, the function $w\mapsto\til G_{I,2,t_2}(t_1,w)-w$ has period $2\pi$, is analytic in $\Omega(\til g_1(t_1,\til J),\R)$, and its imaginary part vanishes on $\R$. 
Applying the maximum principle to the real part of this function, and using (\ref{<1-bound}), we conclude that
$$ |\til G_{I,2,t_2}(t_1,\xi_1(t_1))-\xi_1(t_1)|\le A_1+A_2+A_3,\quad\mbox{if }t_1\le T_1(J_1)\mbox{ and } t_2\le T_2(J_2).$$ 
This together with (\ref{X-+<}) and the estimation of  $|\til g_{I,1,t_2}(t_1,\xi_2(t_2))-\xi_2(t_2)|$ implies that $|X_1-\xi_1+\xi_2|$ is uniformly bounded on $(-\infty,T_1(J_1)]\times (-\infty,T_2(J_2)]$.


Since $G_{I,2,t_2}(t_1,\cdot)$ maps $\TT$ onto $\TT$, and is conformal in the domain that contains the region between $g_1(t_1,J)$ and $\TT$, there must exist $z_1\in g_1(t_1,J)$ such that $|G_{I,2,t_2}(t_1,z_1)|=|z_1|$. Choose $w_1\in \til g_1(t_1,\til J)$ such that $e^i(w_1)=z_1$. Then $\Imm \til G_{I,2,t_2}(t_1,w_1)=\Imm w_1$. From (\ref{<1-bound})   we get
 $|\mA+t_1+t_2|\le A_1+A_2+A_3$, if $t_1\le T_1(J_1)$ and $t_2\le T_2(J_2)$, which
finishes the proof of (i).

\vskip 3mm

Now suppose $t_1+t_2\le -1-2\ln(13)-2\ln(16)$ and $\Imm z_0=\frac{t_1-t_2}2$. Then
$$\Imm z_0-t_1=-\Imm z_0-t_2=-\frac{t_1+t_2}2\ge \frac12+\ln(13)+\ln(16)\ge \ln(8).$$
Since $\til K_1(t_1)\subset\{\Imm z\le \ln(4)+\ln(t_1)\}$ and $\til K_{I,2}(t_2)\subset\{\Imm z\ge -\ln(t_2)-\ln(4)\}$, we have $z_0\in \C\sem \til K_1(t_1)\sem \til K_{I,2}(t_2)$. From (\ref{>1}) and (\ref{>2}) we have \BGE |w_1-(z_0-it_1)|,|w_2-(z_0+it_2)|\le 4e^{\frac{t_1+t_2}2}\le 1/2.\label{>4}\EDE
From (\ref{ma>}), (\ref{>4}), and the upper bound of $t_1+t_2$, we have
$$\Imm w_2+\mA\ge \Imm z_0+t_2-\frac 12-t_1-t_2-\ln(16)=-\frac{t_1+t_2+1}2-\ln(16)\ge \ln(13).$$
Thus, from (\ref{>3}) and the above inequality we have
\BGE |\til g_{1,t_2}(t_1,w_2)-(w_2+i\mA)|\le 10e^{-\mA-\Imm w_2}\le 264e^{\frac{t_1+t_2}2}.\label{>5}\EDE
From (\ref{til-G}), (\ref{>4}), and (\ref{>5}) we see that if $t_1+t_2\le -1-2\ln(13)-2\ln(16)$, then
$$ |\til G_{I,2,t_2}(t_1,w_1)-w_1-i(\mA+t_1+t_2)|\le 272 e^{\frac{t_1+t_2}2},\quad w_1\in\til g_1(t_1,\R_{({t_1-t_2})/2}).$$
The argument between (\ref{<1-bound}) and the end of part (i) can be used here to show that, if $t_1+t_2\le -1-2\ln(13)-2\ln(16)$, then $|\til G_{I,2,t_2}(t_1,\xi_1(t_1))-\xi_1(t_1)|\le 272 e^{\frac{t_1+t_2}2}$ and $|\mA+t_1+t_2|\le 272 e^{\frac{t_1+t_2}2}$. These inequalities together with the uniform limit of $\til g_{I,1,t_2}(t_1,\xi_2(t_2))-\xi_2(t_2)$ and the fact that $T_2(J_2)\le C_{J_2}$ imply that (ii) hold for $j=1$ and $k=2$.
Interchanging $t_1$ and $t_2$, we find that $\mA+t_1+t_2\to 0$ and $X_2-\xi_2+\xi_1\to 0$ as $t_2\to 0$, uniformly in $t_1\in(-\infty,T_1(J_1)]$. From (\ref{Xjk}) we see that $X_2-\xi_2+\xi_1=-(X_1-\xi_1+\xi_2)$, so we have $X_1-\xi_1+\xi_2\to 0$ as $t_2\to 0$, uniformly in $t_1\in(-\infty,T_1(J_1)]$. This completes the proof of part (ii). $\Box$

\vskip 3mm

Let $\ha{\cal D}={\cal D}\cup\{(t_1,-\infty):t_1\in [-\infty,\infty)\}\cup\{(-\infty,t_2):t_2\in [-\infty,\infty)\}$,
 and extend $M$ to $\ha{\cal D}$ such that $M=1$ if $t_1$ or $t_2$ equals $-\infty$. From Proposition \ref{bounded-int},
 we see that $M$ is positive and continuous on $\ha{\cal D}$. So  for any fixed $j\ne k\in\{1,2\}$
 and any $(\F^k_t)$-stopping time $t_k$ which is uniformly bounded above, $M$ is a local martingale in $t_j\in [-\infty,T_j(t_k))$.

\subsection{Local coupling and global coupling} \label{lcaoandglobal}

 Let  $\mu_j$ denote the distribution of $(\xi_j)$, $j=1,2$. Let $\mu=\mu_1\times\mu_2$. Then $\mu$ is
the joint distribution of $(\xi_1)$ and $(\xi_2)$, since $\xi_1$ and $\xi_2$ are independent.
Fix $(J_1,J_2)\in\JP$. From the local martingale property of $M$ and Proposition \ref{bounded-int}, we have
$\EE_\mu[M(T_1(J_1),T_{2}(J_2))]=M(-\infty,-\infty)=1$.
 Define $\nu_{J_1,J_2}$ by $d\nu_{J_1,J_2}=M(T_1(J_1),T_{2}(J_2)) d\mu$. Then $\nu_{J_1,J_2}$ is a probability measure.
Let $\nu_1$ and $\nu_2$ be the two marginal measures of $\nu_{J_1,J_2}$. Then
$d\nu_1/d\mu_1=M(T_1(J_1),-\infty)=1$ and $d\nu_2/d\mu_2=M(-\infty,T_2(J_2))=1$, so $\nu_j=\mu_j$, $j=1,2$.
Suppose temporarily that the distribution of $(\xi_1,\xi_2)$ is $\nu_{J_1,J_2}$ instead of $\mu$. Then the
distribution of each $(\xi_j)$ is still $\mu_j$.

We may now use the argument in Section \ref{local} with a few changes. Here $M(t_1,t_2)$ satisfies (\ref{paM-int}) instead of (\ref{paM}); $\xi_j(t_j)$ does not satisfy (\ref{xi-j}), but is a pre-$(\TT;\kappa)$-Brownian motion with drift $s_j\cdot t$. The traditional Girsanov theorem needs to be modified to work for the current setting.
Eventually, we can conclude that, under the probability measure $\nu_{J_1,J_2}$, for any $j\ne k\in\{1,2\}$, if $t_k$ is a fixed $(\F^k_t)$-stopping time with $t_k\le T_k(J_k)$, and $g_k(t,\cdot)$, $-\infty<t<\infty$, are the inverted whole-plane Loewner maps driven by $\xi_k$, then conditioned on $\F^k_{t_k}$,  after a time-change, $g_k(t_k,K_{I,j}(t_j))$, $-\infty< t_j\le T_j(J_j)$, is a partial disc SLE$(\kappa,\Lambda_{j})$ process in $\D$ started from $0$ with marked point $e^i(\xi_k(t_k))$.

The proof of Theorem \ref{coupling-thm-int-skew} can be now completed by applying the coupling technique.

\section{Partial Differential Equations} \label{Section-PDE}
With Theorem \ref{coupling-thm-int-skew} at hand, to prove the main theorem we need to find particular solutions to (\ref{dot-Gamma=}) that satisfy certain properties. The main part of this section serves this purpose.
From Lemma \ref{Gamma-Psi} we see that solving (\ref{dot-Gamma=}) is equivalent to solving (\ref{PDE-Psi*}) with $\sigma=\frac4\kappa-1$. Throughout this section, we assume that $\kappa>0$ and $\sigma\in[0,\frac 4\kappa)$, and will find solutions to (\ref{PDE-Psi*}) in these cases. In particular, we will obtain solutions to (\ref{dot-Gamma=}) when $\kappa\in(0,4]$.

The solutions to (\ref{PDE-Psi*}) is obtained by construction. We will transform (\ref{PDE-Psi*}) into a similar PDE (\ref{PDE-ha-Psi}), where $\HA_I$ is replaced by $\ha\HA_I$. We know that as $t\to\infty$, $\ha\HA_I(t,\cdot)\to \coth_2$, and PDE (\ref{PDE-ha-Psi}) tends to another PDE (\ref{PDE-ha-Psi-inf}), which has a simple solution $\ha\Psi_\infty$ given by (\ref{Def-Psi0}). Then we let $\ha\Psi_q=\ha\Psi/\ha\Psi_\infty$, and find that $\ha\Psi$ solves (\ref{PDE-ha-Psi}) if and only if $\ha\Psi_q$ solves PDE (\ref{PDE-ha-Psi-q}). A formal solution to (\ref{PDE-ha-Psi-q}) is expressed by a Feynman-Kac formula (\ref{ha-Psi-q-sol}), which involves diffusion processes. In Section \ref{regularity} we prove that the solution given by (\ref{ha-Psi-q-sol}) is smooth, and solves (\ref{PDE-ha-Psi-q}). So we obtain a solution $\Gamma$ to (\ref{dot-Gamma=}). However, such $\Gamma$ does not satisfy (\ref{Gamma-period-s}). For this purpose, note that (\ref{dot-Gamma=}) is a linear PDE, and $\HA_I$ has period $2\pi$, so any translation of $\Gamma$ by an integer multiple of $2\pi$ also solves (\ref{dot-Gamma=}). The solutions to (\ref{dot-Gamma=}) which also satisfy (\ref{Gamma-period-s}) will be obtained by summing over all translations of $\Gamma$ with suitable weights.

In the last subsection, we apply Theorem \ref{coupling-thm-int-skew} to the solutions $\Gamma$ obtained in the previous subsections and obtain the whole-plane reversibility. In fact, we not only prove Theorem \ref{Main-Thm}, but also prove a more general result: the reversibility of skew whole-plane SLE$_\kappa$ processes for $\kappa\in(0,4]$, which is Theorem \ref{Main-Thm-skew}.

For the proofs in this section, the following symbols will be used. For any $n,j\in\N$, we call an $j$-tuple $\lambda=(\lambda_1,\dots,\lambda_j)\in\N^j$ a partition of $n$ if $\lambda_1\ge\dots\ge \lambda_j$ and $\sum_{k=1}^j \lambda_k=n$. The length of such partition is denoted by $l(\lambda)=j$. Let ${\cal P}_n$ denote the set of all partitions of $n$. For example, $(n)$ is the only element in ${\cal P}_n$ with length  $1$. Let ${\cal P}_{\N}=\bigcup_{n\in\N}{\cal P}_n$ denote the set of all partitions.

\subsection{Diffusion processes}
Fix  $\tau\le 0$. For $x\in\R$, let $u(t,x)$, $t\ge 0$, be the solution to \BGE \pa_t u(t,x)=  \tau \tanh_2(u(t,x)+\sqrt\kappa B(t));\quad u(0,x)=x.\label{ODE-u}\EDE
Then $X_x(t):=u(t,x)+\sqrt\kappa B(t)$ satisfies the SDE
\BGE dX_x(t)=\sqrt \kappa dB(t)+ \tau \tanh_2(X_x(t))dt,\quad X_x(0)=x. \label{SDE-X}\EDE

\begin{Lemma}
For any $x\in\R$, we have a.s.\ $\int_0^\infty \tanh_2'(X_x(t))dt=\infty$ and
\BGE \limsup_{t\to\infty} X_x(t)=+\infty, \qquad \liminf_{t\to\infty}X_x(t)=-\infty.\label{Xsupinf*}\EDE
 \label{strip-compare}
\end{Lemma}
 \vskip -5mm
{\bf Proof.} Fix $x\in\R$. Let $X(t)=X_x(t)$. Define
$f(t)=\int_0^t \cosh_2(s)^{-\frac{4}\kappa \tau}ds,\quad t\in\R$. 
Then $f$ is a differentiable increasing odd function and satisfies
$\frac\kappa 2 f''+\tau\tanh_2f'=0$. Let $Y(t)=f(X(t))$. From (\ref{SDE-X})  and It\^o's formula, we have
$dY(t)=f'(X(t))\sqrt\kappa dB(t)$.
Define a time-change function $ u(t)=\int_0^t \kappa f'(X(s))^2  ds$. 
Since $\tau\le 0$, $f'(t)\ge 1$, $t\in\R$. Thus, $u(t)\ge t$ for all $t\in\R$. So $u$ maps $[0,\infty)$ onto $[0,\infty)$, and $Y(u^{-1}(t))$, $0\le t<\infty$, has the distribution of a Brownian motion. Thus, (\ref{Xsupinf*}) holds with $X$ replaced by $Y$, which then implies (\ref{Xsupinf*}). Since $X$ is recurrent, and $\tanh_2'>0$ on $\R$, we immediately have a.s.\ $\int_0^\infty \tanh_2'(X_x(t))dt=\infty$. $\Box$


\begin{Lemma} For any $b,c>0$ and $x\in\R$,
  \BGE \PP[\exists t\ge 0,|X_x(t)|>ct+b]\le 2e^{\frac{2c}\kappa(|x|-b)}.\label{inq}\EDE\label{Lemma-inq}
\end{Lemma}
\vskip -5mm
{\bf Proof.} First, it is well known that (\ref{inq}) holds with $X_x(t)$ replaced by $x+\sqrt\kappa B(t)$.
So it suffices to show that $(|X_x(t)|)$ is bounded above by a process that has the distribution of $(|x+\sqrt\kappa B(t)|)$. This can be proved by using Theorem 4.1 in \cite{BC}. Here we give a direct proof.

Let $Y(t)=|X_x(t)|$ From (\ref{SDE-X}) and Tanaka-It\^o's formula, we have
\BGE Y(t)=|x|+\sqrt\kappa B_0(t)+\frac\kappa 2 \tau \int_0^t \tanh_2(Y(s))ds+L(t),\quad t\ge 0,\label{SDE-Y}\EDE
where $B_0(t)$ is a Brownian motion, $L(t)$ is a non-decreasing function, which satisfies $L(0)=0$ and is constant on every interval of $\{Y(t)>0\}$.

Fix $t_0\ge 0$. There is $t_0'\in[0,t_0]$ such that $L(t)$ is constant on $[t_0',t_0]$. We may assume $t_0'$ is the smallest such number. There are two cases. Case 1: $t_0'=0$. Then $L(t_0)=L(t_0')=L(0)=0$. Since $\tau\le 0$, from (\ref{SDE-Y}), $Y(t_0)\le |x|+\sqrt\kappa B_0(t)$. Case 2: $t_0'>0$. Then $Y(t_0')=0$. Since $\tau\le 0$, from (\ref{SDE-Y}),
$$Y(t_0)-|x|-\sqrt\kappa B_0(t_0)\le Y(t_0')-|x|-\sqrt\kappa B_0(t_0')=-|x|-\sqrt\kappa B_0(t_0').$$
Thus, in either case, we have
$$Y(t_0)\le |x|+\sqrt\kappa B_0(t_0)+\max\{0,\sup_{0\le s\le t_0}\{-|x|-\sqrt\kappa B_0(s)\}\}.$$
The RHS of the above inequality defines a process that has the distribution of $|x+\sqrt\kappa B(t_0)|$, $t_0\ge 0$ (see \cite{RY}), so the proof is completed. $\Box$

\begin{Lemma}  There are $C_n>0$, $n\in\N$, with $C_1=1$,
 such that $$|\tanh_2^{(n)}(x)|\le C_n\tanh_2'(x)\le \frac{C_n}2,\quad x\in\R,\, n\in\N.$$
\label{estoftanh2}
\end{Lemma}
\vskip -5mm
{\bf Proof.} Note that $\tanh_2'(x)=\frac 12\cosh_2^{-2}(x)\in(0,1/2]$. So the second ``$\le$'' holds. By induction, one can prove that for every $n$, there are $a^{(n)}_j\in\R$, $0\le j\le n-1$, such that
$$\tanh_2^{(n)}(x)=\sum_{j=0}^{n-1}a^{(n)}_j  \cosh_2^{-2-j}(x)\sinh_2^{j}(x)=
\sum_{j=0}^{n-1} a^{(n)}_j \cosh_2^{-2}(x)\tanh_2^{j}(x).$$
Since $ |\tanh_2^j(x)|\le 1$ and $\cosh_2^{-2}=2\tanh_2'$, we may choose $C_n=2\sum_j |a^{(n)}_j|$. $\Box$

\begin{Lemma} For every $m\in\N$, there is a polynomial $P_m(t)$ of degree $m-1$ such that for any $t>0$ and $x\in\R$,
$|\frac{\pa^m}{\pa x^m} X_x(t)|\le P_m(t)$. 
  \label{lemma-est-pax-X}
\end{Lemma}
{\bf Proof.} Since $X_x(t)=u(t,x)+\sqrt\kappa B(t)$, $\frac{\pa^n}{\pa x^n} X_x(t)=u^{(n)}(t,x)$. It suffices to show that for every $m\in\N$, there is some polynomial $P_m(t)$ of degree $m-1$, such that \BGE |u^{(m)}(t,x)|\le P_m(t), \quad t>0,\,x\in\R.\label{Pnt2}\EDE Let $f_x(t)= \tau \tanh_2'(X_x(t))$. Since $\tau\le 0$ and $\tanh_2'>0$, $f_x(t)\le 0$. Differentiating (\ref{ODE-u}) w.r.t.\ $x$ and using $u'(0,x)=1$,
we get $ u'(t,x)=\exp(\int_0^t f_x(s)ds)\in(0,1]$. Thus, (\ref{Pnt2}) holds in the case $n=1$ with $P_1(t)\equiv 1$.

Let $n\in\N$, $n\ge 2$. Suppose that (\ref{Pnt2}) holds for any $m\le n-1$. Differentiating (\ref{ODE-u}) $n$ times, by induction we find that there are $b_n(\lambda)\in\R$ for $\lambda\in{\cal P}_n$ with $b_n((n))=\tau$ such that \BGE \pa_t u^{(n)}(t,x)=\sum_{\lambda\in{\cal P}_n} b_n(\lambda)\tanh_2^{(l(\lambda))}(X_x(t)) \prod_{k=1}^{l(\lambda)} u^{(\lambda_k)}(t,x),\quad u^{(n)}(0,x)=0. \label{panu}\EDE
Observe that the term $u^{(n)}(t,x)$ appears only once in (\ref{panu}), i.e, in the case $\lambda=((n))$, and the coefficient is $\tau\tanh_2'(X_x(t))=f_x(t)$. From  Lemma \ref{estoftanh2} and induction hypothesis, there is a polynomial $g_x(t)$ of degree $n-2$ such that
$$\pa_t u^{(n)}(t,x)=f_x(t)u^{(n)}(t,x)+g_x(t),\quad u^{(n)}(0,x)=0.$$
Solving this inequality using the fact that $f_x(t)\le 0$, we can conclude that (\ref{Pnt2}) holds in the case $m=n$, which finishes the proof. $\Box$

\subsection{Some estimations}
We will need some estimations about the limits of $\ha\HA_{I}-\tanh_2$ as $t\to\infty$.  Let \BGE \ha\HA_{I,q}(t,z)=\ha\HA_I(t,z)-\tanh_2(z).\label{ha-HA-iq}\EDE
From (\ref{ha-HA-I}) we have
\BGE\ha\HA_{I,q}'(t,x)=\sum_{2|n\ne 0}\tanh_2'(x-nt)=\sum_{2|n\ne 0} \frac 12\cosh_2^{-2}(x-nt)>0.\label{expofhaHAIq'}\EDE

\begin{Lemma} Let $C_n$, $n\in\N$, be as in Lemma \ref{estoftanh2}. Note that $C_1=1$. Then
\BGE |\ha\HA_{I,q}(t,x)|\le \frac{|x|}t+3+\frac{2e^{-t}}{1-e^{-2t}},\quad t>0,\,x\in\R.\label{est-ha-HA1}\EDE
\BGE |\ha\HA_{I,q}^{(n)}(t,x)|\le C_n\Big(\frac 12 + \frac{4e^{-t}}{1-e^{-2t}}\Big),\quad t>0,\,x\in\R,\,n\in\N.\label{est-ha-HA'1}\EDE
Moreover, for any $c>0$,
\BGE |\ha\HA_{I,q}(t,x)|\le \frac{2e^{(c-2)t}}{1-e^{-2t}},\quad \mbox{if }t>0,\,x\in\R,\,|x|\le ct.\label{est-ha-HA2}\EDE
\BGE |\ha\HA_{I,q}^{(n)}(t,x)|\le C_n\frac{4e^{(c-2)t}}{1-e^{-2t}},\quad \mbox{if } t>0,\,x\in\R,\,|x|\le ct,\,n\in\N.\label{est-ha-HA'2}\EDE
\label{est-ha-HA'}\end{Lemma}
{\bf Proof.} We first show (\ref{est-ha-HA2}). From (\ref{ha-HA-I}) and (\ref{ha-HA-iq}) we have
$$\ha\HA_{I,q}(t,x)=\sum_{m=1}^{\infty}(\tanh_2(x-2mt)+\tanh_2(x+2mt))$$
$$= \sum_{m=1}^{\infty}\Big(-\frac{e^{2mt}-e^x}{e^{2mt}+e^x}+\frac{e^{2mt}-e^{-x}}{e^{2mt}+e^{-x}}\Big)= \sum_{m=1}^{\infty}\frac{2(e^x-e^{-x})}{e^{2mt}+e^{-2mt}+e^x+e^{-x}}.$$
Thus, \BGE |\ha\HA_{I,q}(t,x)|\le  \sum_{m=1}^{\infty}\frac{2e^{|x|}}{e^{2mt}}=\frac{2e^{|x| -2t}}{1-e^{-2t}}.\label{|ha-HA-Iq|}\EDE
Then (\ref{est-ha-HA2}) is a direct consequence of this inequality.

Secondly, we show (\ref{est-ha-HA1}). Since $|\tanh_2(x)|\le 1$, from (\ref{ha-HA-iq}) it suffices to show
\BGE |\ha\HA_{I}(t,x)|\le \frac{|x|}t+2+\frac{2e^{-t}}{1-e^{-2t}}.\label{est-ha-HA10}\EDE
We first consider the case $|x|\le t$. From (\ref{|ha-HA-Iq|}) we have
\BGE |\ha\HA_I(t,x)|\le |\tanh_2(x)|+|\ha\HA_{I,q}(t,x)|\le 1+\frac{2e^{|x| -2t}}{1-e^{-2t}}\le 1+\frac{2e^{-t}}{1-e^{-2t}}.\label{|ha-HA-I|}\EDE
Thus, (\ref{est-ha-HA10}) holds in this case.

Then we consider the case $|x|\ge t$. There exists $m\in\N$ such that $(2m-1)t\le |x|\le (2m+1)t$. Since $\ha\HA_{I}$ is odd, we only need to consider the case that $(2m-1)t\le x\le (2m+1)t$. Let $x_0=x-2mt$. Then $|x_0|\le t$. From (\ref{ha-HA-I-period}) we have
$\ha\HA_I(t,x)=2m+\ha\HA_I(t,x_0)$. From (\ref{|ha-HA-I|}) with $x=x_0$ we have
$$|\ha\HA_I(t,x)|\le 2m+|\ha\HA_I(t,x_0)|\le 2m+1+\frac{2e^{-t}}{1-e^{-2t}}\le \frac{|x|}t+2+\frac{2e^{-t}}{1-e^{-2t}},$$
where the last inequality uses $\frac{|x|}t\ge 2m-1$. So we have (\ref{est-ha-HA10}) and (\ref{est-ha-HA1}).

Now we prove (\ref{est-ha-HA'1}) and (\ref{est-ha-HA'2}). From (\ref{expofhaHAIq'}) we have
$$0<\ha\HA_{I,q}'(t,x)=\sum_{2|n\ne 0} \frac 12\cosh_2^{-2}(|nt-x|)\le \sum_{2|n\ne 0} \frac 12\cosh_2^{-2}(|n|t-|x|)$$
\BGE =\sum_{m=1}^\infty \cosh_2^{-2}(2mt-|x|)\le 4\sum_{m=1}^\infty e^{|x|-2mt}=\frac{4e^{|x|-2t}}{1-e^{-2t}}.\label{ha-HAq<}\EDE
which implies (\ref{est-ha-HA'2}) in the case $n=1$. From (\ref{ha-HAq<}) we have $\ha\HA_I'(t,x) <\frac 12+ \frac{4e^{-t}}{1-e^{-2t}}$ if $|x|\le t$. Since $\ha\HA_I'$ has period $2t$, this inequality holds for all $x\in\R$. Since $\ha\HA_{I,q}'<\ha\HA_I'$,  (\ref{est-ha-HA'1}) holds in the case $n=1$.
From (\ref{expofhaHAIq'}) and Lemma \ref{estoftanh2} we have $ |\ha\HA_{I,q}^{(n)}(t,x)|\le C_n \ha\HA_{I,q}'(t,x)$.
So (\ref{est-ha-HA'1}) and (\ref{est-ha-HA'2}) in the case $n\ge 2$ follow from those in the case $n=1$. $\Box$

\begin{Lemma}
  For every $n\in\N\cup\{0\}$, there is a constant $D_n>0$ such that for any $j\in\{1,2\}$, $t>0$, and $x\in\R$,
\BGE |\pa_t^j\ha\HA_{I,q}^{(n)}(t,x)|\le D_n\Big(\frac{|x|}t+3+\frac{2e^{-t}}{1-e^{-2t}}\Big)^j\Big(\frac 12 + \frac{4e^{-t}}{1-e^{-2t}}\Big). \label{est-pat-ha-HA1}\EDE
Moreover, for any $n\in\N\cup\{0\}$ and $c>0$, there is a constant $D_n>0$, such that
\BGE |\pa_t^j\ha\HA_{I,q}^{(n)}(t,x)|\le D_n\Big(\frac{2e^{(c-2)t}}{1-e^{-2t}}\Big)^{j+1},\quad \mbox{if }t>0,\,x\in\R,\,|x|\le ct; \label{est-pat-ha-HA3}\EDE
\label{lemma-est-pat-ha}
\end{Lemma}
\vskip -3mm
{\bf Proof.} Let $A(t,x)=\frac{|x|}t+3+\frac{2e^{-t}}{1-e^{-2t}}$, $B(t,x)=\frac 12 + \frac{4e^{-t}}{1-e^{-2t}}$, and $C_c(t,x)=\frac{2e^{(c-2)t}}{1-e^{-2t}}$. In this proof, by $X\lesssim Y$ we mean that there is a constant $C$ such that $X\le C Y$. Here $C$ may depend on $n$ if $X$ depends on $n$. From (\ref{est-ha-HA1}), (\ref{est-ha-HA'1}), (\ref{est-ha-HA2}), and (\ref{est-ha-HA'2}), we see that
\BGE |\HA_{I,q}(t,x)|\lesssim A(t,x),\qquad |\HA_{I,q}^{(n)}(t,x)|\lesssim B(t,x)\lesssim A(t,x),\quad x\in\R,\, n\in\N.\label{<=AB}\EDE
\BGE |\HA_{I,q}^{(n)}(t,x)|\lesssim C_c(t,x),\quad \mbox{ if }x\in\R,\,|x|\le ct,\, n\in\N\cup\{0\}.\label{<=C}\EDE

As $t\to\infty$, $\ha\HA_I\to\tanh_2$. Then (\ref{pat-ha-HA}) becomes
$ 0=\tanh_2''+\tanh_2'\tanh_2$, 
which can be proved directly. From  (\ref{pat-ha-HA}), (\ref{ha-HA-iq}), and the above equation, we get
\BGE \pa_t \ha\HA_{I,q}=\ha\HA_{I,q}'' +\ha\HA_{I,q}'\ha\HA_{I,q}+\tanh_2'\ha\HA_{I,q}+\ha\HA_{I,q}'\tanh_2.\label{a}\EDE
Then (\ref{est-pat-ha-HA1}) and (\ref{est-pat-ha-HA3}) in the case $j=1$ and $n=0$ follow from (\ref{<=AB}), (\ref{<=C}), (\ref{a}), and Lemma \ref{estoftanh2}.

Differentiating (\ref{pat-ha-HA}) w.r.t.\ $x$ twice, we get
\BGEN\pa_t \ha\HA_I'=\ha\HA_I'''+\ha\HA_I''\ha\HA_I+(\ha\HA_I')^2.
\EDEN
\BGEN\pa_t \ha\HA_I''=\ha\HA_I^{(4)}+\ha\HA_I'''\ha\HA_I+3\ha\HA_I''\ha\HA_I'.
\EDEN
Differentiating (\ref{pat-ha-HA}) w.r.t.\ $t$ and using the above two displayed formulas, we obtain
$$ \pa_t^2 \ha\HA_I=\ha\HA_I^{(4)}+2\ha\HA_I'''\ha\HA_I+ 4\ha\HA_I''\ha\HA_I'+\ha\HA_I''(\ha\HA_I)^2+2(\ha\HA_I')^2\ha\HA_I.$$
As $t\to\infty$, this equation tends to the following equation, which can also be checked directly.
$$0=\tanh_2^{(4)}+2\tanh_2'''\tanh_2+4\tanh_2''\tanh_2'+\tanh_2''\tanh_2^2+2(\tanh_2')^2\tanh_2.$$
From (\ref{ha-HA-iq}),  and the above two equations, we compute
$$\pa_t^2 \ha\HA_{I,q}=\ha\HA_{I,q}^{(4)} +2\ha\HA_{I,q}'''\ha\HA_{I,q}+2\tanh_2'''\ha\HA_{I,q} +2\ha\HA_{I,q}'''\ha\HA_{I,q}+4\ha\HA_{I,q}''\ha\HA_{I,q}' $$
$$+4\ha\HA_{I,q}''\tanh_2'+\tanh_2''(\ha\HA_{I,q})^2+2\ha\HA_{I,q}''\ha\HA_{I,q}\tanh_2+2\tanh_2''\ha\HA_{I,q}\tanh_2 $$
$$ +\ha\HA_{I,q}''(\tanh_2)^2+ \ha\HA_{I,q}''(\ha\HA_{I,q})^2+2(\ha\HA_{I,q}')^2\ha\HA_{I,q}+4\ha\HA_{I,q}'\tanh_2'\ha\HA_{I,q}$$ \BGE +2(\tanh_2')^2\ha\HA_{I,q} +2(\ha\HA_{I,q}')^2\tanh_2+4\ha\HA_{I,q}'\tanh_2'\tanh_2 .\label{b}\EDE
Then (\ref{est-pat-ha-HA1}) and (\ref{est-pat-ha-HA3}) in the case $j=2$ and $n=0$ follow from (\ref{<=AB}), (\ref{<=C}), (\ref{b}), and Lemma \ref{estoftanh2}.

Differentiate (\ref{a}) and (\ref{b}) $n$ times w.r.t.\ $x$. We see that $\pa_t \ha\HA_{I,q}^{(n)}$ can be expressed as a sum of finitely many terms, whose factors are $\HA_{I,q}^{(k)}$ or $\tanh_2^{(k)}$, $k\in\N\cup\{0\}$. In every term, the factors of the kind $\HA_{I,q}^{(k)}$ appear at most twice, and the factor $\HA_{I,q}$ appears at most once. So we derive (\ref{est-pat-ha-HA1}) and (\ref{est-pat-ha-HA3}) in the case $j=1$ and $n\in\N$ from (\ref{<=AB}), (\ref{<=C}),  and Lemma \ref{estoftanh2}.  We see that $\pa_t^2 \ha\HA_{I,q}^{(n)}$ can be expressed as a sum of finitely many terms, whose factors are constant, $\HA_{I,q}^{(k)}$, or $\tanh_2^{(k)}$. In every term, the factors of the kind $\HA_{I,q}^{(k)}$ appear at most three times, and the factor $\HA_{I,q}$ appears at most twice. So we derive (\ref{est-pat-ha-HA1}) and (\ref{est-pat-ha-HA3}) in the case $j=2$ and $n\in\N$ from (\ref{<=AB}), (\ref{<=C}), and Lemma \ref{estoftanh2}. $\Box$

\subsection{Feynman-Kac expression}\label{Feynman}
We begin with a lemma, which can be proved directly. Recall the definition of $\ha\HA_I$ in (\ref{ha-HA-0*}).

\begin{Lemma}
For  $\Psi$ and $\ha\Psi$ on $(0,\infty)\times\R$, the following are equivalent:  \BGE \ha\Psi(t,x)=e^{\frac{x^2}{2\kappa t}} \Big(\frac \pi t\Big)^{\sigma+\frac 12}\Psi\Big(\frac{\pi^2}t,\frac\pi t x\Big).\label{ha-Psi-1}\EDE
  \BGE\Psi(t,x)= e^{-\frac{x^2}{2\kappa t}} \Big(\frac \pi t\Big)^{\sigma+\frac 12}\ha\Psi\Big(\frac{\pi^2}t,\frac\pi t x\Big).\label{ha-Psi-2}\EDE
If the above two equalities hold, then $\Psi$ satisfies (\ref{PDE-Psi*}) if and only if $\ha\Psi$ satisfies
\BGE -\pa_t {\ha \Psi} = \frac \kappa 2 \ha\Psi'' + \sigma\ha \HA_I' \ha \Psi.\label{PDE-ha-Psi} \EDE
  \label{Lemma-transform}
\end{Lemma}
\vskip -5mm

As $t\to\infty$, $\ha\HA_I'\to\tanh_2'$, so (\ref{PDE-ha-Psi}) tends to
\BGE -\pa_t {\ha \Psi}_\infty  = \frac \kappa 2 \ha\Psi_\infty'' + \sigma \tanh_2'(x) \ha \Psi_\infty .\label{PDE-ha-Psi-inf}\EDE
Let $\tau$ be the non-positive root of the equation $\frac{\tau^2}{2\kappa}=\frac\tau 4 +\frac\sigma2$, i.e.,
$\tau=\kappa/4-\sqrt{{\kappa^2}/{16}+\kappa\sigma}$. 
Then  $\tau=\frac\kappa 2-2$ when $\sigma=\frac 4\kappa-1$. It is easy to check that (\ref{PDE-ha-Psi-inf}) has a simple solution:
\BGE \ha\Psi_\infty(t,x)=e^{-\frac{\tau^2t}{2\kappa}}\cosh_2^{\frac 2\kappa\tau}(x).\label{Def-Psi0}\EDE

Recall the $\ha\HA_{I,q}$ defined in (\ref{ha-HA-iq}).  The proof of the following lemma is straightforward.
\begin{Lemma}
  Let $\ha\Psi$ and $\ha\Psi_q$ be defined on $(0,\infty)\times\R$, and satisfy $\ha\Psi=\ha\Psi_\infty\ha\Psi_q$, where $\Psi_\infty$ is defined by (\ref{Def-Psi0}). Then $\ha\Psi$ satisfies (\ref{PDE-ha-Psi}) if and only if $\ha\Psi_q$ satisfies
\BGE -\pa_t{\ha\Psi}_q=\frac\kappa 2 \ha\Psi_q''+\tau\tanh_2 \ha\Psi_q'+\sigma\ha\HA_{I,q}'\ha\Psi_q.\label{PDE-ha-Psi-q}\EDE \label{Lemma-psi-q}
\end{Lemma}

\vskip -5mm

Suppose $\ha\Psi_q$ solves (\ref{PDE-ha-Psi-q}). Let $X_{x_0}(t)$ be as in (\ref{SDE-X}). Fix $t_0>0$ and $x_0\in\R$. Let
$$ M(t)= \ha\Psi_q(t_0+t,X_{x_0}(t))\exp\Big(\sigma \int_{0}^{t} \ha\HA_{I,q}'(t_0+s,X_{x_0}(s))ds\Big).$$ 
From (\ref{SDE-X}), (\ref{PDE-ha-Psi-q}), and It\^o's formula, we see that $M(t)$ is a local martingale. If $M(t)$ is a martingale on $[0,\infty]$, and $\ha\Psi_q\to 1$ as $t\to\infty$, then from $M_0=\ha\Psi_q(t_0,x_0)$ we have
\BGE \ha\Psi_q(t_0,x_0)=\EE\Big[\exp\Big(\sigma \int_{0}^{\infty} \ha\HA_{I,q}'(t_0+s,X_{x_0}(s))ds\Big)\Big].\label{ha-Psi-q-sol}\EDE
This Feynman-Kac formula holds under many additional assumptions. We do not try to prove it. Instead, we now define $\ha\Psi_q$ by (\ref{ha-Psi-q-sol}). We will prove that $\ha\Psi_q$ is finite and differentiable, and solves (\ref{PDE-ha-Psi-q}).

\subsection{Regularity}\label{regularity}
Fix $c_0\in(1+\frac\kappa 4\sigma ,2)$. This is possible because
$\sigma\in[0,\frac 4\kappa)$. Then we have
\BGE \exp\Big({\frac{\sigma}{2(c_0-1)}-\frac{2}\kappa}\Big)<1.\label{exp<1}\EDE
Throughout this subsection, we use $C$ to denote a positive constant, which depends only on $\kappa,\sigma,c_0$, and could change between lines. The symbol $X\lesssim Y$ means that $X\le CY$ for some $C$. Let $\alpha(t)=\frac 4{1-e^{-2t}}$. Then $t^{-1}+1\lesssim\alpha(t)\lesssim t^{-1}+1$. For  $m\in\N\cup\{0\}$, let ${\cal E}_m$ denote the event that $|X_x(s)|\le s+m$ for all $s\ge 0$. From (\ref{inq}) we have
\BGE \PP[{\cal E}_m^c]\le 2e^{\frac 2\kappa(|x|-m)},\quad m\in\N\cup\{0\}.\label{p(E)<}\EDE

\begin{Proposition}
$\ha\Psi_q$ is finite and satisfies
\BGE 1\le \ha\Psi_q(t,x)\le \exp\Big(C(t^{-1}+1) e^{(c_0-2)t}\Big)(1+Ce^{\frac 2\kappa|x|-\frac 2\kappa c_0t}).\label{psiq1}\EDE
\end{Proposition}
{\bf Proof.} Fix $t>0$ and $x\in\R$. Assume that ${\cal E}_m$ occurs for some $m\in\N\cup\{0\}$.  If $s\ge \frac{m-c_0t}{c_0-1}$ then $|X_x(s)|\le s+m\le c_0(s+t)$, so from (\ref{est-ha-HA'2}) with $C_1=1$ we have $$\ha\HA_{I,q}'(t+s,X_{x}(s)) \le \frac{4e^{(c_0-2)(s+t)}}{1-e^{-2(s+t)}}\le \alpha(t) e^{(c_0-2)(s+t)}.$$ If $0\le s\le \frac{m-c_0t}{c_0-1}$, from $-1\le c_0-2$ and (\ref{est-ha-HA'1}) with $C_1=1$, we have
$$\ha\HA_{I,q}'(t+s,X_{x}(s)) < \frac 12+\frac{4e^{-(s+t)}}{1-e^{-2(s+t)}}\le \frac 12+ \alpha(t) e^{(c_0-2)(s+t)},$$
Since $c_0-2<0$, at the event ${\cal E}_m$,
\BGE \int_{0}^{\infty} \ha\HA_{I,q}'(t+s,X_{x}(s))ds\le  \frac 12\cdot \frac{(m-c_0t)\vee 0}{c_0-1}+\frac{\alpha(t)e^{(c_0-2)t}}{2-c_0};\label{int-eq0}\EDE
Let $H(t)=\exp (\sigma \int_{0}^{\infty} \ha\HA_{I,q}'(t+s,X_{x}(s))ds )$. From (\ref{p(E)<}) and (\ref{int-eq0}) we have
$$\ha\Psi_q(t,x)= \EE [1_{{\cal E}_{\lfloor c_0t\rfloor}} H(t)  ]+\sum_{m=\lfloor c_0t\rfloor}^\infty \EE [1_{{\cal E}_{m+1}\sem{\cal E}_{m}}H(t) ]$$
\BGE \le \exp\Big(\frac{\sigma\alpha(t)e^{(c_0-2)t}}{2-c_0}\Big)+\sum_{m=\lfloor c_0t\rfloor}^\infty 2e^{\frac 2\kappa(|x|-m)}
 \exp\Big(\frac 12 \frac{\sigma(m+1-\lfloor c_0t\rfloor)}{c_0-1}+\frac{\sigma \alpha(t)e^{(c_0-2)t}}{2-c_0} \Big).\label{RHS1}\EDE
Change index using $m=l+\lfloor c_0t\rfloor$. The second term of the RHS of (\ref{RHS1}) equals
\BGE 2\exp\Big(\frac {2|x|}\kappa-\frac{2 \lfloor c_0t\rfloor}\kappa+\frac{\sigma}{2(c_0-1)}+\frac{\sigma \alpha(t)e^{(c_0-2)t}}{2-c_0} \Big)\sum_{l=0}^\infty
\exp\Big(\frac{\sigma}{2(c_0-1)}-\frac 2\kappa\Big)^l.\label{argument}\EDE
From (\ref{exp<1}), the infinite sum is finite. Thus, from $\ha\HA_{I,q}'>0$, $\sigma\ge 0$, and (\ref{RHS1}), we have
$$ 1\le \ha\Psi_q(t,x)\le\exp\Big(\frac{\sigma\alpha(t)e^{(c_0-2)t}}{2-c_0}\Big)\Big(1+Ce^{\frac {2|x|}\kappa - \frac{2c_0t}\kappa}\Big).$$ 
Then (\ref{psiq1}) follows from this formula and that $\alpha(t)\lesssim t^{-1}+1$. $\Box$

\vskip 4mm

Let $n\in\N$. Formally differentiate (\ref{ha-Psi-q-sol}) $n$ times w.r.t.\ $x$. If the differentiation commutes with the integration and expectation at every time, then we should have
\BGE \ha\Psi_q^{(n)}(t,x)=\EE\Big[\exp\Big(\sigma \int_{0}^{\infty} \ha\HA_{I,q}'(t+s,X_{x}(s))ds\Big) \cdot {\bf Q}_{0,n}(v_{0,k,\lambda}(t,x))\Big],\label{Qn}\EDE
where ${\bf Q}_{0,n}$ is a polynomial of degree $\le n$ without constant term in the following variables:
\BGE v_{0,k,\lambda}(t,x):=\int_0^\infty \ha\HA_{I,q}^{(k)}(t+s,X_x(s))\prod_{r=1}^{l(\lambda)} \frac{\pa^{\lambda_r}}{\pa x^{\lambda_r}} X_x(s)ds,\quad k\in\N,\,\lambda\in{\cal P}_{\N}.\label{v0}\EDE

With ${\bf Q}_{0,0}\equiv 1$, (\ref{Qn}) becomes (\ref{ha-Psi-q-sol}). Let $n\in\N\cup\{0\}$. Formally differentiate (\ref{Qn}) w.r.t.\ $t$. If the differentiation commutes with the integration and expectation, then we should have
 \BGE \pa_t{\ha\Psi}_q^{(n)}(t,x)= \EE\Big[\exp\Big(\sigma \int_{0}^{\infty} \ha\HA_{I,q}'(t+s,X_{x}(s))ds\Big) \cdot {\bf Q}_{1,n}(v_{0,k,\lambda},v_{1,k,\lambda})],
 \label{Q1n}\EDE
where ${\bf Q}_{1,n}$ is a polynomial of degree $\le n+1$ without constant term in the variables $v_{0,k,\lambda}$ defined by (\ref{v0}) and
\BGE v_{1,k,\lambda}(t,x):=\int_0^\infty \pa_t \ha\HA_{I,q}^{(k)}(t+s,X_x(s))\prod_{r=1}^{l(\lambda)} \frac{\pa^{\lambda_r}}{\pa x^{\lambda_r}} X_x(s)ds,\quad k\in\N,\lambda\in{\cal P}_{\N}\cup\{\N^0\}.\label{v1}\EDE
Here by $\lambda\in\N^0$ we mean that the factor $\prod \frac{\pa^{\lambda_r}}{\pa x^{\lambda_r}} X_x(s)$ disappears. Moreover, in every term of ${\bf Q}_{1,n}$, factors $v_{1,k,\lambda}$ appear at most once.

Formally differentiate (\ref{Q1n}) w.r.t.\ $t$. If the differentiation commutes with the integration and expectation, then we should have
\BGE \pa_t^2 {\ha\Psi}_q^{(n)}(t,x)= \EE\Big[\exp\Big(\sigma \int_{0}^{\infty} \ha\HA_{I,q}'(t+s,X_{x}(s))ds\Big) \cdot {\bf Q}_{2,n}(v_{0,k,\lambda},v_{1,k,\lambda},v_{2,k,\lambda})\Big],\label{Q2n}\EDE
where ${\bf Q}_{2,n}$ is a polynomial of degree $\le n+2$ without constant term in the variables $v_{0,k,\lambda}$ defined by (\ref{v0}), $v_{1,k,\lambda}$ defined by (\ref{v1}), and
 $$ v_{2,k,\lambda}(t,x):=\int_0^\infty \pa_t^2 \ha\HA_{I,q}^{(k)}(t+s,X_x(s))\prod_{j=1}^{l(\lambda)} \frac{\pa^{\lambda_j}}{\pa x^{\lambda_j}} X_x(s)ds,\quad k\in\N,\lambda\in{\cal P}_{\N}\cup\{\N^0\}.$$
 Moreover, in every term of ${\bf Q}_{2,n}$, factors $v_{2,k,\lambda}$ appears at most once; when a factor $v_{2,k,\lambda}$ appears, factors $v_{1,k,\lambda}$ disappear; and when factors $v_{2,k,\lambda}$ disappear, factors $v_{1,k,\lambda}$ appear at most twice.

Now we suppose ${\cal E}_m$ occurs for some $m\in\N\cup\{0\}$.
Using (\ref{est-ha-HA'1}), (\ref{est-ha-HA'2}), (\ref{v0}), Lemma \ref{lemma-est-pax-X}, and the argument in (\ref{int-eq0}), we conclude that, for any $k\in\N$ and $\lambda\in{\cal P}_{\N}$, there is a polynomial $P_{k,\lambda}$ with no constant term such that
\BGE |v_{0,k,\lambda}(t,x)|\le P_{k,\lambda}((m-c_0t)\vee 0)+C \alpha(t)e^{(c_0-2)t}.\label{v<}\EDE
Let $j\in\{1,2\}$ and $n\in\N\cup\{0\}$. If $s\ge \frac{m-c_0t}{c_0-1}$ then $|X_x(s)|\le s+m\le c_0(s+t)$, so from (\ref{est-pat-ha-HA3})
we have
$$|\pa_t^j\ha\HA_{I,q}^{(n)}(t+s,X_x(s))|\le D_n\Big(\frac{2e^{(c_0-2)(t+s)}}{1-e^{-2t}}\Big)^{j+1}\lesssim \alpha(t)^{j+1}e^{(c_0-2)(t+s)}.$$
If $m\ge c_0t$ and $0\le s\le \frac{m-c_0t}{c_0-1}$, from (\ref{est-pat-ha-HA1})
and the definition of ${\cal E}_m$, we see that for $j=1,2$,
$$|\pa_t^j\ha\HA_{I,q}^{(n)}(t,X_x(s))|\le D_n\Big(\frac{|X_x(s)|}{t+s}+3+\frac{2e^{-t}}{1-e^{-2t}}\Big)^j\Big(\frac 12 + \frac{4e^{-t}}{1-e^{-2t}}\Big)$$
$$\le D_n\Big(\frac{m-c_0t}{t}+c_0+4+\frac{2e^{-t}}{1-e^{-2t}}\Big)^j\Big(\frac 12 + \frac{4e^{-t}}{1-e^{-2t}}\Big)\lesssim ((m-c_0t)^j+1)\alpha(t)^{j+1}.$$
Thus, from Lemma \ref{lemma-est-pax-X}, for  $k\in\N$ and $\lambda\in{\cal P}_{\N}\cup\{\N^0\}$,
\BGE |v_{j,k,\lambda}(t,x)|\lesssim \alpha(t)^{j+1}( e^{(c_0-2)t}+P_{j,k,\lambda}((m-c_0t)\vee 0)),\quad j=1,2,\label{vj<}\EDE
where $P_{j,k,\lambda}$ is a polynomial with no constant term.

Let $(j,n)\in \{0,1,2\}\times(\N\cup\{0\})\sem\{(0,0)\}$. From (\ref{v<}), (\ref{vj<}), and the properties of ${\bf Q}_{0,n}$, $n\in\N$, ${\bf Q}_{1,n}$ and ${\bf Q}_{2,n}$, $n\in\N\cup\{0\}$, we see that, at the event ${\cal E}_m$,
\BGE|{\bf Q}_{j,n}|\lesssim \alpha(t)^{2j} [P_{j,n}((m-c_0t)\vee 0)+Q_{j,n}((m-c_0t)\vee 0)\alpha(t)^ne^{(c_0-2)t}],
\label{Qj<}\EDE
where $P_{j,n}$ and $Q_{j,n}$ are polynomials, and $P_{j,n}(0)=0$.

\begin{Proposition} For $(j,n)\in \{0,1,2\}\times(\N\cup\{0\})\sem\{(0,0)\}$, $$ \EE\Big[\exp\Big(\sigma \int_{0}^{\infty} \ha\HA_{I,q}'(t+s,X_{x}(s))ds\Big)\cdot|{\bf Q}_{j,n}|\Big]$$ 
\BGE \lesssim \exp\Big(C(t^{-1}+1)e^{(c_0-2)t}\Big)(t^{-n-2j}+1)(e^{(c_0-2)t}+e^{\frac{2|x|}\kappa-\frac{2c_0t}\kappa}),  \label{psiqjn}\EDE
\label{Psi-q<infty}
\end{Proposition}
\vskip -3mm
{\bf Proof.} Let $H_{j,n}(t)=\exp\Big(\sigma \int_{0}^{\infty} \ha\HA_{I,q}'(t+s,X_{x}(s))ds\Big)\cdot|{\bf Q}_{j,n}|$. Recall that (\ref{int-eq0}) and (\ref{Qj<}) hold at the event ${\cal E}_m$. Using (\ref{p(E)<}) and the argument in (\ref{RHS1}) and (\ref{argument}), we see that
$$\EE [H_{j,n}(t)]=\EE [1_{{\cal E}_{\lfloor c_0t\rfloor}}H_{j,n}(t)]+\sum_{m=\lfloor c_0t\rfloor}^\infty\EE[ 1_{{\cal E}_{m+1}\sem {\cal E}_m} H_{j,n}(t)] $$
$$\lesssim \exp\Big(C\alpha(t)e^{(c_0-2)t}\Big)\alpha(t)^{2j+n} e^{(c_0-2)t}+ \exp\Big(\frac{2|x|}\kappa-\frac{2c_0t}{\kappa}+C\alpha(t)e^{(c_0-2)t}\Big)\cdot$$
$$\cdot\sum_{l=0}^\infty \alpha(t)^{2j}(P_{j,n}(l+1)+Q_{j,n}(l+1)\alpha(t)^ne^{(c_0-2)t})\exp\Big(\frac{\sigma}{2(c_0-1)}-\frac 2\kappa\Big)^l.$$
Then (\ref{psiqjn}) follows from (\ref{exp<1}) and that $\alpha(t)\lesssim t^{-1}+1$ and $\alpha(t)^{2j+n} \lesssim t^{-n-2j}+1$.  $\Box$

\begin{Theorem}
  The function $\ha\Psi_q$ is $C^{\infty,\infty}$ differentiable and solves (\ref{PDE-ha-Psi-q}). Moreover, for $j\in\{0,1,2\}$, $n\in\N\cup\{0\}$, there is a positive continuous function $c_{j,n}(t)$ on $(0,\infty)$ such that for any $t\in(0,\infty)$ and $x\in\R$,
  $|\pa_t^j \ha\Psi_q^{(n)}(t,x)|\le c_{j,n}(t)e^{\frac2\kappa|x|}$. \label{Thm-regul}
\end{Theorem}
{\bf Proof.} For $n\in\N\cup\{0\}$, define $\ha\Psi_q^{[0,n]}$, $\ha\Psi_q^{[ 1,n]}(t,x)$ and $\ha\Psi_q^{[ 2,n]}(t,x)$ to be equal to the RHS of (\ref{Qn}), (\ref{Q1n}) and (\ref{Q2n}), respectively. From the above two propositions, these functions are well defined, and there are positive continuous functions $c_{j,n}(t)$ on $(0,\infty)$ such that
\BGE |\ha\Psi_q^{[j,n]}(t,x)|\le c_{j,n}(t)e^{\frac2\kappa|x|},\quad j=0,1,2,\quad n\in\N\cup\{0\}.\label{cnt}\EDE

Let $n\in\N\cup\{0\}$, $j\in\{0,1,2\}$, $t\in(0,\infty)$, and $x_1<x_2\in\R$. Since $|\ha\Psi_q^{[ j,n+1 ]}|$ satisfies (\ref{cnt}), from Fubini's Theorem, we have \BGE \int_{x_1}^{x_2} \ha\Psi_q^{[ j, n+1 ]} (t,x)dx=\ha\Psi_q^{[ j,n ]}(t,x_2)-\ha\Psi_q^{[ j,n ]}(t,x_1).\label{Lip}\EDE Thus, $\ha\Psi_q^{[ j, n ]}$ is absolutely continuous in $x$ when $t$ is fixed, and its partial derivative w.r.t.\ $x$ is a.s.\ equal to $\ha\Psi_q^{[ j,n+1 ]}$. Since $\ha\Psi_q^{[ j,n+1 ]}$ is continuous in $x$ for fixed $t$, we see that $\ha\Psi_q^{[ j, n ]}$ is continuously differentiable in $x$, and the partial derivative exactly equals $\ha\Psi_q^{[j, n ]}$. The above holds for any $n\in\N$, so $\ha\Psi_q^{[ j,0 ]} $ is $C^\infty$ differentiable in $x$ when $t$ is fixed, and $\ha\Psi_q^{[j, n]}$ is its $n$-th partial derivative w.r.t.\ $x$. Especially, since $\ha\Psi_q=\ha\Psi_q^{[0,0]}$, we see that $\ha\Psi_q$ is $C^\infty$ differentiable in $x$ when $t$ is fixed, and $\ha\Psi_q^{[0,n]}$ is its $n$-th partial derivative w.r.t.\ $x$.

A similar argument using Fubini's Theorem shows that, for any $n\in\N\cup\{0\}$, $j\in\{0,1\}$, $\ha\Psi_q^{[j, n ]}$ is absolutely continuous in $t$ when $x$ is fixed, and its partial derivative w.r.t.\ $t$ is a.s.\ equal to $\ha\Psi_q^{[j+1,n ]}$. So $\ha\Psi_q^{[ 0, n ]}$ is continuously differentiable in $t$ when $x$ is fixed, and the partial derivative exactly equals $\ha\Psi_q^{[1, n ]}$.  From (\ref{cnt}) and  (\ref{Lip}), we see that  $\ha\Psi_q^{[j,n]}$ is locally uniformly Lipschitz continuous in $x$. We have seen that $\ha\Psi_q^{[ j, n ]}$ is continuous in $t$ for every fixed $x$. So $\ha\Psi_q^{[ j, n ]}$ is continuous in both $t$ and $x$. Thus, $\ha\Psi_q=\ha\Psi_q^{[0,0]}$ is $C^{1,\infty}$ differentiable.

Fix $t_0\in(0,\infty)$  and $x_0\in\R$. Let $ M(t)=\EE\Big[\exp\Big(\sigma\int_{0}^\infty \ha\HA_{I,q}'(t_0+s,X_{x_0}(s))ds\Big)\Big|\F_t\Big]$, $t\ge 0$.
Then $M(t)$, $0\le t<\infty$, is a uniformly integrable martingale. From (\ref{ha-Psi-q-sol}) we have
\BGE M(t)=\ha\Psi_q(t_0+t,X_{x_0}(t))\exp\Big(\sigma\int_{0}^{t} \ha\HA_{I,q}'(t_0+s,X_{x_0}(s))ds\Big).\label{M=}\EDE
From (\ref{SDE-X}), It\^o's formula, and the differentiability of $\ha\Psi_q$, we see that $\ha\Psi_q$ solves (\ref{PDE-ha-Psi-q}) for $t\ge t_0$. Since this is true for any $t_0\in(0,\infty)$, $\ha\Psi_q$ solves (\ref{PDE-ha-Psi-q}).

Since $\ha\Psi_q$ is $C^{1,\infty}$ differentiable, the same is true for the RHS of (\ref{PDE-ha-Psi-q}). Thus, $\pa_t\ha\Psi_q$ is also $C^{1,\infty}$ differentiable. So $\ha\Psi_q$ is $C^{2,\infty}$ differentiable. Iterating this argument, we conclude that $\ha\Psi_q$ is $C^{\infty,\infty}$ differentiable. 
The previous argument shows that $\pa_t^j\ha\Psi_q^{(n)}=\ha\Psi_q^{[j,n]}$ for any $j\in\{0,1,2\}$ and $n\in\N\cup\{0\}$. The bounds of $|\pa_t^j\ha\Psi_q^{(n)}|$ then follow from (\ref{cnt}). $\Box$

\begin{Theorem}
  Let $\ha\Psi_0=\ha\Psi_q\cdot\ha\Psi_\infty$, where $\ha\Psi_\infty$ is defined by (\ref{Def-Psi0}). Then $\ha\Psi_0$ is a positive $C^{\infty,\infty}$ differentiable function on $(0,\infty)\times\R$ and solves (\ref{PDE-ha-Psi}). Moreover, $j\in\{0,1,2\}$, for $n\in\N\cup\{0\}$, there is a positive continuous function $c_{j,n}(t)$ on $(0,\infty)$ such that, for any $t\in(0,\infty)$ and $x\in\R$,
  $|\pa_t^j \ha\Psi_0^{(n)}(t,x)|\le c_{j,n}(t)e^{\frac2\kappa|x|}$. \label{Cor-regul}
\end{Theorem}
{\bf Proof.} Since $\ha\Psi_q$ and $\ha\Psi_\infty$ are both positive and $C^{\infty,\infty}$ differentiable, the same is true for $\ha\Psi_0=\ha\Psi_q\cdot\ha\Psi_\infty$. Since $\ha\Psi_q$ solves (\ref{PDE-ha-Psi-q}), from Lemma \ref{Lemma-psi-q}, $\ha\Psi_0$ solves (\ref{PDE-ha-Psi}). From Lemma \ref{estoftanh2}, (\ref{Def-Psi0}), and that $\tau\le 0$, we see that for any $j,n\in\N\cup\{0\}$, $|\pa_t^j\ha\Psi_\infty^{(n)}(t,x)|$ is bounded by a positive continuous function in $t$, which, together with Theorem \ref{Thm-regul}, implies the upper bounds of $|\pa_t^j \ha\Psi_0^{(n)}(t,x)|$. $\Box$

\begin{Theorem} Let $\Psi_0$ be the $\Psi$ in (\ref{ha-Psi-2}) with $\ha\Psi$ replaced by $\ha\Psi_0$ in the above theorem. Then $\Psi_0$ is a $C^{\infty,\infty}$ differentiable positive function on $(0,\infty)\times\R$ and solves (\ref{PDE-Psi*}). Moreover, for $j\in\{0,1,2\}$, $n\in\N\cup\{0\}$, there is a function $h_{j,n}(t,|x|)$, which is a polynomial in $|x|$ for any fixed $t$, and every coefficient is a positive continuous function in $t$, such that for any $t\in(0,\infty)$ and $x\in\R$,
  $|\pa_t^j  \Psi_0^{(n)}(t,x)|\le h_{j,n}(t,|x|)e^{-\frac{x^2}{2\kappa t}+\frac{2\pi|x|}{\kappa t}}$.
    \label{Cor-regul2}
\end{Theorem}
{\bf Proof.} Since $\ha\Psi_0>0$, $\Psi_0>0$ also. The differentiability of $\Psi_0$ is obvious. Since $\ha\Psi_0$ solves (\ref{PDE-ha-Psi-q}), from Lemma \ref{Lemma-transform}, $\Psi_0$ solves (\ref{PDE-Psi*}). Let $\Psi_a(t,x)=\ha\Psi_0(\frac{\pi^2}t,\frac\pi t x)$.
From Theorem \ref{Cor-regul}, it is straightforward to check that for every $j\in\{0,1,2\}$, $n\in\N\cup\{0\}$, there is a function $f_{j,n}(t,|x|)$, which is a polynomial in $|x|$ of degree $j$ when $t$ is fixed, and every coefficient is a positive continuous function in $t$, such that
\BGE\Big|\pa_t^j \Psi^{(n)}_0(t,x)\Big|\le f_{j,n}(t,|x|) e^{\frac 2\kappa \frac\pi t|x|},\quad t>0,x\in\R.\label{bound1}\EDE
It is easy to verify that for every $j\in\{0,1,2\}$, $n\in\N\cup\{0\}$, there is a function $g_{j,n}(t,|x|)$, which is a polynomial in $|x|$, and every coefficient is a positive continuous function in $t$, such that
\BGE\Big|\pa_t^j \pa_x^n (e^{-\frac{x^2}{2\kappa t}}(\frac \pi t)^{\sigma+\frac 12})\Big| \le g_{j,n}(t,|x|)  e^{-\frac{x^2}{2\kappa t}},\quad t>0, x\in\R.\label{bound2}\EDE
From (\ref{ha-Psi-2}), $\Psi_0(t,x)=e^{-\frac{x^2}{2\kappa t}}(\frac \pi t)^{\sigma+\frac 12}\Psi_a(t,x)$. So we get the upper bounds of $|\pa_t^j \Psi_0^{(n)}(t,x)|$ from (\ref{bound1}) and (\ref{bound2}). $\Box$

\begin{Theorem}
Let $\Psi_0$ be as in the above theorem. Let $\Gamma_0=\Psi_0\Theta_I^{-\frac2\kappa}$ and $\Gamma_m(t,x)=\Gamma_0(t,x-2m\pi)$, $m\in\Z$. For $s_0\in\R$, let $\Gamma_{\langle s_0\rangle}=\sum_{m\in\Z}e^{\frac{2\pi}\kappa m s_0}\Gamma_m$. Then $\Gamma_{\langle s_0\rangle}$ is a $C^{\infty,\infty}$ differentiable positive function on $(0,\infty)\times\R$, satisfies (\ref{Gamma-period-s}), and solves (\ref{dot-Gamma-c=}).
 \label{Thm-Psi-period*}
\end{Theorem}
{\bf Proof.} Let $\Psi_m(t,x)=\Psi_0(t,x-2m\pi)$ for $m\in\Z$ and $\Psi_{\langle s_0\rangle}=\sum_{m\in\Z}e^{\frac{2\pi}\kappa m s_0}\Psi_m$. Since $\Theta_I$ has period $2\pi$, we have $\Gamma_{\langle s_0\rangle}=\Psi_{\langle s_0\rangle}\Theta_I^{-\frac2\kappa}$. Since $\Theta_I$ is a $C^{\infty,\infty}$ differentiable positive function with period $2\pi$, from Lemma \ref{Gamma-Psi} we suffice to show that $\Psi_{\langle s_0\rangle}$ is a $C^{\infty,\infty}$ differentiable positive function, satisfies (\ref{Gamma-period-s}), and solves (\ref{PDE-Psi*}). It is clear from the definition that $\Psi_{\langle s_0\rangle}$ satisfies (\ref{Gamma-period-s}).
Since $\Psi_0$ is a $C^{\infty,\infty}$ differentiable positive function that solves (\ref{PDE-Psi*}), and $\HA_I$ has period $2\pi$, every $\Psi_m$ also satisfies these properties. So $\Psi_{\langle s_0\rangle}$ is positive. The  upper bounds of $|\pa_t^j  \Psi_0^{(n)}(t,x)|$ imply that $\Psi_{\langle s_0\rangle}$ is finite, and   the series $\sum_{m\in\Z}e^{\frac{2\pi}\kappa m s_0}\pa_t^j\Psi_m^{(n)}$ converges locally uniformly for every $j,n\ge 0$.  Fubini's Theorem implies that  $\Psi_{\langle s_0\rangle}$ is $C^{\infty,\infty}$ differentiable and $\pa_t^j\Psi_{\langle s_0\rangle}^{(n)}=\sum_{m\in\Z}e^{\frac{2\pi}\kappa m s_0}\pa_t^j\Psi_m^{(n)}$
 Thus,  $\Psi_{\langle s_0\rangle}$ also solves (\ref{PDE-Psi*}). $\Box$

\subsection{Distributions}

\begin{Proposition}  Let $p>0$, $s_0\in\R$, and $x_0,y_0\in\R$. Let $\Gamma_m$, $m\in\Z$, and $\Gamma_{\langle s_0\rangle}$ be as in Theorem \ref{Thm-Psi-period*}. Let $\Lambda_*=\kappa\frac{\Gamma_*'}{\Gamma_*}$ for $*\in\{m,\langle s_0\rangle\}$. For $m\in\Z$, let $\til\beta_m$ be the covering annulus SLE$(\kappa,\Lambda_0)$ trace in $\St_p$ started from $x_0$ with marked point $y_0+2m\pi+p i$. Let $\til\beta_{\langle s_0\rangle}$ be the covering annulus SLE$(\kappa,\Lambda_{\langle s_0\rangle})$ trace in $\St_p$ started from $x_0$ with marked point $y_0+p i$. Let $\PP_{\til\beta,m}$, $m\in\Z$, and $\PP_{\til\beta,\langle s_0\rangle}$ denote the distributions of $\til\beta_m$, $m\in\Z$, and $\til\beta_{\langle s_0\rangle}$, respectively. Then
\BGE \PP_{\til\beta,\langle s_0\rangle}=\sum_{m\in\Z} e^{\frac{2\pi}\kappa m s_0}\frac{\Gamma_m(p,x_0-y_0)}{\Gamma_{\langle s_0\rangle}(p,x_0-y_0)}\,\PP_{\til\beta,m}.\label{conv8}\EDE
\label{conv}
\end{Proposition}
{\bf Proof.} For $m\in\Z$, let $\xi_m(t)$, $0\le t<p$, be the solution to (\ref{xi-crossing}) with $\Lambda=\Lambda_0$ and $y_0$ replaced by $y_0+2m\pi$. Let $\xi_{\langle s_0\rangle}(t)$ be the solution to (\ref{xi-crossing}) with $\Lambda=\Lambda_{\langle s_0\rangle}$. Then the covering annulus Loewner traces of modulus $p$ driven by $\xi_m$, $m\in\Z$, and $\xi_{\langle s_0\rangle}$   have distributions $\PP_{\til\beta,m}$, $m\in\Z$, and $\PP_{\til\beta,\langle s_0\rangle}$, respectively.
Let $X_m(t)=\xi_m(t)-\Ree \til g^{\xi_m}(t,y_0+2m\pi+pi)+2m\pi$, $m\in\Z$, and $X_{\langle s_0\rangle}(t)=\xi_{\langle s_0\rangle}(t)-\Ree \til g^{\xi_{\langle s_0\rangle}}(t,y_0+pi)$. Since $\Gamma_m(t,x)=\Gamma_0(t,x-2m\pi)$, we have $\Lambda_m(t,x)=\Lambda_0(t,x-2m\pi)$.
Since $\Ree g(t,y+pi)=\til g_I(t,y)$ for $y\in\R$, and $\HA_I$ is odd and has period $2\pi$, from (\ref{annulus-eqn-inver-covering}), we find that, for $*\in\{m,\langle s_0\rangle\}$, with $\Phi_*:=\Lambda_*+\HA_I$, $X_*(t)$ satisfies
   $$ dX_*(t)=\sqrt\kappa dB(t)+\Phi_*(p-t,X_*(t))dt,\quad X_*(0)=x_0-y_0.\label{X-SDE-m*}$$
Let $\PP_{X,*}$ denote the distributions of $(X_*(t))$. Since $\xi_{*}(t)= X_{*}(t)+y_0-\int_0^t \HA_I( p-r, X_{*}(r))dr$, $0\le t<p$,
we suffice to show that (\ref{conv8}) holds with the subscripts ``$\til\beta$'' replaced by ``$X$''.
The rest of the proof is a standard application of Girsanov theorem. One may check that for every $m\in\Z$, $M_m(t):=e^{\frac{2\pi}\kappa m s_0}\frac{\Gamma_m(p-t,X_{\langle s_0\rangle}(t))}{\Gamma_{\langle s_0\rangle}(p-t,X_{\langle s_0\rangle}(t))}$ is a nonnegative martingale w.r.t.\ $\PP_{X,\langle s_0\rangle}$, and satisfies that $\frac{dM_m(t)}{M_m(t)}=(\Lambda_m-\Lambda_{\langle s_0\rangle})\frac{dB(t)}{\sqrt\kappa}$ and $\sum_{m\in\Z} M_m(t)=1$; and we have $\frac{d\PP_{X,m}}{d\PP_{X,\langle s_0\rangle}}=\frac{M_m(\infty)}{M_m(0)}$. $\Box$

\vskip 4mm
\no{\bf Remark.} Since $\Gamma_{\langle s_0\rangle}$ satisfies (\ref{Gamma-period-s}), $\Lambda_{\langle s_0\rangle}$ has period $2\pi$. So $\Lambda_{\langle s_0\rangle}$ is a crossing annulus drift function, and we could define the annulus SLE$(\kappa,\Lambda_{\langle s_0\rangle})$ process. However, each $\Lambda_m$ does not have period $2\pi$. It only makes sense to define the covering annulus SLE$(\kappa,\Lambda_m)$ processes.

\begin{Proposition}
  Let $p>0$ and $x_0,y_0\in\R$. Let $\Gamma_0$ be as in Theorem \ref{Thm-Psi-period*}, and  $\Lambda_0=\kappa\frac{\Gamma_0'}{\Gamma_0}$.
  Let $\til\beta(t)$, $0\le t<p$, be the covering annulus SLE$(\kappa,\Lambda_0)$ trace in $\St_p$ started from $x_0$ with marked point $y_0+pi$. Then a.s.\ $\dist(y_0+ pi,\til\beta([0,p))+2\pi\Z)=0$.
\end{Proposition}
{\bf Proof.} Let $\xi(t)$ be the driving function, and $\til g(t,\cdot)$, $0\le t<p$, be the covering Loewner maps. Then $\til g(t,\cdot)$ maps $\St_p\sem (\til\beta([0,p))+2\pi\Z)$ conformally onto $\St_{p-t}$, and maps $\R_p$ onto $\R_{p-t}$. From Koebe's $1/4$ Theorem, we suffice to show that a.s.\ $\til g'(t,y_0+ pi)\cdot\frac{p}{p-t}\to\infty$ as $t\to p$.

Let $X(t)=\xi(t)-\Ree \til g(t,y_0+p i)$ and $\Phi_0=\Lambda_0+\HA_I$. Then $X(t)$  satisfies the SDE:
$$dX(t)=\sqrt\kappa dB(t)+{\Phi_{0}(p-t, X(t))}\,dt,\quad 0\le t<p.$$
From (\ref{deriv2*}) we have $\ln(\til g'(t,y_0+pi)\cdot\frac p{p-t})=\int_0^t (\HA_I'(p-s,X(s))+\frac{1}{p-s})ds$.
Let $\ha\Phi_0=\kappa\frac{\ha\Psi_0'}{\ha\Psi_0}$.
Since $\Psi_0$ and $\ha\Psi_0$ satisfy (\ref{ha-Psi-1}), we have $\ha\Phi_0(s,z)=\frac\pi s\Phi_0(\frac{\pi^2}s,\frac\pi s z)+\frac{z}s$. Let $\ha p=\frac{\pi^2}p$ and $\ha X(t)=\frac{\ha p+t}{\pi}X(p-\frac{\pi^2}{\ha p+t})$, $0\le t<\infty$. Then $\ha X(0)=\frac{\ha p}\pi X(0)=\frac\pi p(x_0-y_0)$. Applying It\^o's formula and time-change of a semimartingale, we see that $\ha X(t)$ satisfies the SDE:
$$d\ha X(t)=\sqrt\kappa \ha B(t)+\ha\Phi_0(\ha p+t,\ha X(t))dt,\quad 0\le t<\infty,$$
for some standard Brownian motion $\ha B(t)$. Changing variables using $\ha s=\frac{\pi^2}{p-s}-\ha p$, we get
$$\int_0^t (\HA_I'(p-s,X(s))+\frac{1}{p-s})ds=\int_0^{\ha t}\Big(\HA_I'(\frac{\pi^2}{\ha p+\ha s},X(p-\frac{\pi^2}{\ha p+\ha s}))+\frac{\ha p+s}{\pi^2}\Big)\frac{\pi^2}{(\ha p+\ha s)^2}\,d\ha s$$
$$=\int_0^{\ha t}\Big(\frac{\pi^2}{(\ha p+\ha s)^2}\HA_I'(\frac{\pi^2}{\ha p+\ha s},\frac{\pi}{\ha p+\ha s}\ha X(\ha s))+\frac1{\ha p+s} \Big) d\ha s=\int_0^{\ha t} \ha\HA_I'(\ha p+\ha s,\ha X(\ha s))d\ha s,$$
where $\ha t=\frac{\pi^2}{p-t}-\ha p$, and the last equality follows from (\ref{ha-HA-0*}). So we have
$$\lim_{t\to p^-}\ln(\til g'(t,y_0+ pi)\cdot\frac{p}{p-t})=\int_0^\infty  \ha\HA_I'(\ha p+\ha s,\ha X(\ha s))d\ha s\ge \int_0^\infty \tanh_2'(\ha X(\ha s))d\ha s,$$
where the last inequality follows from (\ref{ha-HA-I}).

From Girsanov theorem and the fact that $\kappa\frac{\ha\Psi_0'}{\ha\Psi_0}=\kappa\frac{\ha\Psi_q'}{\ha\Psi_q}+\tau\tanh_2$, we find that the distribution of $(\ha X(t))$ is equivalent to that of $(X_{\frac\pi p(x_0-y_0)}(t))$ defined by (\ref{SDE-X}), and the Radon-Nikodym derivative is $M(\infty)/M(0)$, where $M(t)$ is defined by (\ref{M=}).  Since $(X_{\frac\pi p(x_0-y_0)}(t))$ is homogeneous and recurrent, we have a.s.\ $\int_0^\infty \tanh_2'(X_{\frac\pi p(x_0-y_0)}(t))dt=\infty$, which implies that a.s.\  $\int_0^\infty \tanh_2'(\ha X(\ha s))d\ha s=\infty$.
Thus, a.s.\ $\til g'(t,y_0+ pi)\cdot\frac{p}{p-t}\to \infty$ as $t\to p^-$. $\Box$

\begin{Corollary}
  Let $p>0$, $s_0\in\R$, and $x_0,y_0\in\R$. Let $\Gamma_{\langle s_0\rangle}$ be as in Theorem \ref{Thm-Psi-period*}, and  $\Lambda_{\langle s_0\rangle}=\kappa\frac{\Gamma_{\langle s_0\rangle}'}{\Gamma_{\langle s_0\rangle}}$. Let $\beta(t)$, $0\le t<p$, be the annulus SLE$(\kappa,\Lambda_{\langle s_0\rangle})$ trace in $\St_p$ started from $e^{ix_0}$ with marked point $e^{-p+iy_0}$. Then a.s.\ $\dist(e^{-p+iy_0},\beta([0,p)))=0$.\label{sub-lim-union*}
\end{Corollary}
{\bf Proof.} This follows immediately from the above two propositions. $\Box$

\subsection{Reversibility}\label{Reversibility}
Throughout this subsection we assume that $\kappa\in(0,4]$, $s_0\in\R$, $\sigma=\frac4\kappa-1$, and $\Lambda_{\langle s_0\rangle}$ is given by Proposition \ref{conv}. 
We will prove the theorem below which generalizes Theorem \ref{Main-Thm}.
\begin{Theorem} Let $\kappa\in(0,4]$ and $s_0\in\R$. If $\beta(t)$, $-\infty\le t<\infty$, is a whole-plane SLE$(\kappa,s_0)$ trace in $\ha\C$ from $a$ to $b$, then the reversal of $\beta$, up to a time-change, has the distribution of a whole-plane SLE$(\kappa,s_0)$ trace in $\ha\C$ from $b$ to $a$.
  \label{Main-Thm-skew}
\end{Theorem}
{\bf Proof.} We only need to consider the case $a=0$ and $b=\infty$. Let $\Gamma_{\langle s_0\rangle}$ be given by Theorem \ref{Thm-Psi-period*} with $\sigma=\frac4\kappa-1$. Then $\Gamma_{\langle s_0\rangle}$ solves (\ref{dot-Gamma=}) and satisfies (\ref{Gamma-period-s}). We now apply Theorem \ref{coupling-thm-int-skew} to $\Gamma=\Gamma_{\langle s_0\rangle}$. Let $\Lambda_j$, $s_j$ and $\beta_{I,j}(t)$, $j=1,2$, be given by Theorem \ref{coupling-thm-int-skew}. Then for $j=1,2$, $\beta_{I,j}$ is a whole-plane SLE$(\kappa,s_j)$ trace in $\ha\C$ from $0$ to $\infty$, and satisfies that, for any $t_2\in\Q$, conditioned on $\beta_{I,2}(s)$, $-\infty\le s\le  t_2$, after a time-change, the curve $\beta_{I,1}(t_1)$, $-\infty\le t_1<T_1( t_2)$, has the distribution of a disc SLE$(\kappa,\Lambda_1)$ trace in $\C\sem I_0(\beta_{I,2}([-\infty,t_2]))$ started from $0$ with marked point $\beta_{I,2}(t_2)$, where $T_1(t_2)$ is the maximal number in $(-\infty,+\infty]$ such that  $\beta_1(t)\cap \beta_{2}([-\infty,t_2])=\emptyset$ for $-\infty< t<T_1(t_2)$.

Let  $\xi_2$ be the driving function for $(\beta_{I,2}(t))$, and $g_2(t,\cdot)$, $-\infty<t<\infty$, be the inverted whole-plane Loewner maps driven by $\xi_2$. Then $g_2(t,\cdot)$ maps $\C\sem I_0(\beta_{I,2}([-\infty,t_2]))$ conformally onto $\D$, fixes $0$, and takes $\beta_{I,2}(t_2)$ to $e^i(\xi_x(t_2))$. Thus, conditioned on $\beta_{I,2}(s)$, $-\infty\le s\le  t_2$, $g_2(t,\beta_{I,1}(t_1))$, $0\le t_1<T_1(t_2)$, is a time-change of a disc  SLE$(\kappa,\Lambda_1)$ trace in $\D$ started from $0$ with marked point $e^i(\xi_x(t_2))$. Since $\Lambda_1=\Lambda=\kappa\frac{\Gamma_{\langle s_0\rangle}'}{\Gamma_{\langle s_0\rangle}}$,
from Corollary \ref{sub-lim-union*} and the relation between the disc SLE$(\kappa,\Lambda)$ process and the annulus SLE$(\kappa,\Lambda)$ process, we conclude that a.s.\ $e^i(\xi_2(t_2))$ is a subsequential limit of $g_2(t_2,\beta_{I,1}(t))$ as $t\to T_1(t_2)^-$. Thus, $\beta_2(t_2)$ is a subsequential limit of $\beta_{I,1}(t)$ as $t\to T_1(t_2)^-$. If $T_1(t_2)=\infty$, then $\lim_{t\to T_1(t_2)^-} \beta_{I,1}(t)=\infty=\beta_2(-\infty)\ne \beta_2(t_2)$, which is a.s.\ a contradiction. So $T_1(t_2)<\infty$ a.s., and we have $\beta_{I,1}(T_1(t_2))=\lim_{t\to T_1(t_2)^-} \beta_{I,1}(t)=\beta_2(t_2)$ a.s.  Since $\Q$ is countable, we conclude that, a.s.\ $\beta_{I,1}(T_1(t_2))=\beta_2(t_2)$ for every $t_2\in\Q$, which implies that a.s.\ $\beta_2(\R)\subset\beta_{I,1}(\R)$.
Since both $\beta_{I,1}$ and $\beta_2$ are simple, and the initial (resp.\ final) point of $\beta_{I,1}$ agrees with the final (resp.\ initial) point of $\beta_2$,  we see that $\beta_2$ is a reversal of $\beta_{I,1}$. Now $\beta_{I,1}$ is a whole-plane SLE$(\kappa,s_0)$ trace in $\ha\C$ from $0$ to $\infty$, and $\beta_{I,2}$ is a whole-plane SLE$(\kappa,-s_0)$ trace in $\ha\C$ started from $0$ to $\infty$. Since $I_0$ is conjugate conformal, $\beta_2=I_0(\beta_{I,2})$ is  a whole-plane SLE$(\kappa,s_0)$ trace in $\ha\C$ from $\infty$ to $0$. So we proved the theorem in the case $a=0$ and $b=\infty$.
$\Box$

\begin{Theorem} If $\beta(t)$, $0\le t<\infty$, is a radial SLE$(\kappa,-s_0)$ trace in a simply connected domain $D$ from $a$ to $b$, then a.s.\ $\lim_{t\to\infty}\beta(t)=b$, and after a time-change, the reversal of $\beta$ becomes a disc SLE$(\kappa,\Lambda_{\langle s_0\rangle})$ trace in $D$ started from $b$ with marked point $a$.
  \label{reversal-radial}
\end{Theorem}
{\bf Proof.} This follows from the property of the coupling in Theorem \ref{Main-Thm-skew} and the relation between whole-plane SLE$(\kappa,s_0)$ and radial SLE$(\kappa,-s_0)$. $\Box$

\begin{Theorem} Let $D$ be a doubly connected domain with two boundary points  $a,b$ lying on different boundary components. If $\beta(t)$, $0\le t<p$, is an annulus SLE$(\kappa,\Lambda_{\langle s_0\rangle})$ trace in $D$ started from $a$ with marked point $b$, then  $\lim_{t\to p}\beta(t)=b$, and after a time-change, the reversal of $\beta$ becomes an annulus SLE$(\kappa,\Lambda_{\langle s_0\rangle})$ trace in $D$ started from $b$ with marked point $a$.   \label{reversal-ann}
\end{Theorem}
{\bf Proof.} This follows from  the property of the coupling in Theorem \ref{Main-Thm-skew}, and the relation between disc SLE$(\kappa,\Lambda_{\langle s_0\rangle})$ and annulus SLE$(\kappa,\Lambda_{\langle s_0\rangle})$. $\Box$

\vskip 4mm
\no{\bf Remark.} For $\kappa\in(0,6)$ and $\sigma=\frac 12+\frac 1\kappa \in[0,\frac 4\kappa)$, the $\Lambda_{\langle 0\rangle}$ given by Proposition \ref{conv} can be used to decompose an annulus SLE$_\kappa$ process (without marked point). The statement is similar to Lemma 3.1 in \cite{duality2}.

\section{Some Particular Solutions}\label{solution}

In this section,  for $\kappa\in\{4,2,3,0,16/3\}$, we will find  solutions to  the PDE for $\Lambda$ ((\ref{dot-Lambda=}) and (\ref{dot-Lambda=2})) and the PDE for $\Gamma$ ((\ref{dot-Gamma=}) and (\ref{dot-Gamma=2})), which can be expressed in terms of $\HA$ and $\HA_I$. Since $\Lambda=\kappa\frac{\Gamma'}\Gamma$, multiplying a function in $t$ to $\Gamma$ does not change the value of $\Lambda$. So we may as well consider the following PDEs for $\Gamma$, where $C(t)$ is some real valued continuous depending only on $t$:
  \BGE {\pa_t{\Gamma}}=\frac{\kappa}2  {\Gamma''}
+\HA_I{\Gamma'}+\Big(\frac 3\kappa-\frac 12\Big)\HA_I'\Gamma+  C(t){\Gamma}.\label{C'-1}\EDE
 \BGE {\pa_t{ \Gamma}}=\frac{\kappa}2  { \Gamma''}
+\HA{ \Gamma'}+\Big(\frac 3\kappa-\frac 12\Big)\HA' \Gamma+  C(t){ \Gamma}.\label{C'-2}\EDE

\subsection{$\kappa=4$}
Let $\kappa =4$. From Lemma \ref{Gamma-Psi} we see that  if $\Psi$ solves.
\BGE \pa_t\Psi=2\Psi'',\label{heat-Psi}\EDE
then $\Gamma=\Psi\Theta_I^{-2/\kappa}$ solves (\ref{dot-Gamma=}). Similarly, $\Gamma=\Psi\Theta^{-2/\kappa}$ solves (\ref{dot-Gamma=2}) if $\Psi$ solves (\ref{heat-Psi}). The solutions to (\ref{heat-Psi}) are well-known. For example, we have the following solutions: $e^{2c^2t+cx}$, $\frac 1{\sqrt{8\pi t}}e^{-\frac{(x-c)^2}{8t}}$, $e^{-t/2}\sin_2(x-c)$, $\Theta(2t,x-c)$, and $\Theta_I(2t,x-c)$. The function $\Theta_I(2t,x-\pi)$ corresponds to the solution $\Gamma(t,x)=\Theta_I(2t,x-\pi)\Theta_I(t,x)^{-1/2}$ of (\ref{C'-1}), which agrees with the solution given by Section \ref{Feynman} for $\kappa=4$ and $\sigma=\frac 4\kappa-1=0$. Some of these solutions are related to the Gaussian free field (\cite{SS}) in doubly connected domains.

\subsection{$\kappa=2$} \label{Section-kappa2}
Let $\kappa=2$. In this case if $\Xi$ on $(0,\infty)\times\R$ solves
\BGE \pa_t \Xi=\Xi''+\Xi'\HA_I+C(t)\Xi\label{Theta}\EDE
then $\Gamma:=\Xi'$ solves (\ref{C'-1}).  Similarly, if $\Xi$ on
$(0,\infty)\times(\R\sem\{2n\pi:n\in\Z\})$ solves \BGE \pa_t \Xi=\Xi''+\Xi'\HA+C(t)\Xi.\label{Theta2}\EDE
then $\Gamma:=\Xi'$ solves  (\ref{C'-2}).

From (\ref{pat-HA}) we see that $\Xi_1=\HA_I$ solves (\ref{Theta}) and
$\Xi_2=\HA$ solves (\ref{Theta2}) with $C(t)= 0$. It is also easy to check that $\Xi_3(t,x)=t\HA_I(t,x)+x$ solves
(\ref{Theta}) and $\Xi_4(t,x)=t\HA(t,x)+x$ solves (\ref{Theta2}) with $C(t)=0$. The $\Xi_3$ corresponds to the solution $\Gamma(t,x)=t\HA_I'(t,x)+1$, which agrees with the solution given by Section \ref{Feynman} for $\kappa=2$ and $\sigma=\frac 4\kappa-1=1$. Such $\Gamma$ is also the density function of the distribution of the limit point of an annulus SLE$_2$ trace.


We now derive more solutions. Fix $t>0$.  Let $L_t=\{2n\pi+i2kt:n,k\in\Z\}$. Let $F_{1,t}$ denote the set of odd analytic functions $f$ on $\C\sem L_t$ such that each $z\in L_t$ is a simple pole of $f$, $2\pi$ is a period of $f$, and $i2t$ is an antiperiod of $f$, i.e.,
$f(z+i2t)=-f(z)$.  Let $F_{2,t}$ denote the set of odd analytic functions $f$ on $\C\sem L_t$
such that each $z\in L_t$ is a simple pole of $f$, $2\pi$ is an antiperiod of $f$, and $i2t$ is a  period of $f$.
Let $F_{3,t}$ denote the set of odd analytic functions $f$ on $\C\sem L_t$ such that each $z\in L_t$ is a simple pole of $f$,
and both $2\pi$ and $i2t$ are antiperiods of $f$.
Define
$$ \Xi_1(t,z)=\HA(2t,z)-\HA_I(2t,z), \quad
\Xi_2(t,z)=\frac 12\HA(\frac t2,\frac z2)-\frac 12\HA(\frac t2,\frac z2+\pi),$$ 
$$ \Xi_3(t,z)=\frac 12\HA(t,\frac z2)-\frac 12\HA_I(t,\frac z2)-\frac 12\HA(t,\frac z2+\pi)+\frac 12\HA_I(t,\frac z2+\pi).$$
From the properties of $\HA$ and $\HA_I$, it is easy to check that $F_{j,t}$ is the linear space spanned by $\Xi_j(t,\cdot)$ for $j=1,2,3$.
For $j=1,2,3$, Define
$$J_j=\pa_t\Xi_j-\Xi_j''-\Xi_j'\HA,\quad C_j(t)=\frac12{\Res_{z=0} J_j(t,\cdot)} .$$
Fix $t>0$. Note that $0$ is a simple pole of both $\HA(t,\cdot)$ and $\Xi_1(t,\cdot)$ of residue $2$. It is easy to conclude that $0$ is also a simple pole of $J_1(t,\cdot)$. From that $\Xi_1(t,\cdot)\in F_{1,t}$, that $\HA(t,\cdot)$ has period $2\pi$,
and that $\HA(t,z+2\pi)=\HA(t,z)-2i$, it is easy to check that $J_1(t,\cdot)\in F_{1,t}$ as well. So
 $J_1(t,\cdot)=C_1(t)\Xi_1(t,\cdot)$. Thus, $\Xi_1$ solves  (\ref{Theta2}).
 Similarly, $\Xi_2$ and $\Xi_3$ both solve  (\ref{Theta2}).

\subsection{$\kappa=3$}
Let $\kappa=3$. Let $\Xi_j$, $j=1,2,3$, be as in the previous subsection. For $j=1,2,3$, let $\Gamma_j=\Xi_j$, and define
$$H_j=\pa_t{\Gamma_j}-\frac 32\Gamma_j''-\HA\Gamma_j'
-\frac 12\HA'\Gamma_j,\quad C_j(t)=\frac12{\Res_{z=0} H_j(t,\cdot)}.$$
Using the argument in the last subsection, we find that $H_j(t,\cdot)\in F_{j,t}$ for any $t>0$. So $H_j(t,\cdot)=C_j(t)\Gamma_j(t,\cdot)$.
Thus, $\Gamma_1,\Gamma_2,\Gamma_3$ solve (\ref{C'-2}).
For $j=4,5,6$, let $\Gamma_j(t,z)=\Gamma_{j-3}(t,z+it)$. Since $\HA_I(t,z)=\HA(t,z+it)+i$, $\Gamma_4,\Gamma_5,\Gamma_6$ solve (\ref{C'-1}).

For $j=2,3$, $\Gamma_j$ takes positive real values on $(0,2\pi)+4\pi\Z$, takes negative real values on $(-2\pi,0)+4\pi\Z$, and has antiperiod $2\pi$. So $\Lambda_j:=3\frac{\Gamma_j'}{\Gamma_j}$ is a chordal-type annulus drift function that solves (\ref{dot-Lambda=2}) for $\kappa=3$. It is worth to mention that the annulus SLE$(\kappa;\Lambda_j)$ process preserves the following local martingale, which resembles the $G(\Omega,a,b,z)$ in Proposition 11 of \cite{SS-3}. The proof uses the fact that $\Gamma_j$ solves (\ref{C'-2}) for $z\in\C\sem\{\mbox{poles}\}$.

\begin{Proposition}
Let $j\in\{2,3\}$ and $p>0$.  Let $x_0\in\R$ and $z_0\in\R\sem(x_0+2\pi\Z)$. Let $\xi(t)$, $0\le t<T$, be the driving function for the covering annulus SLE$(\kappa;\Lambda_j)$ process in $\St_p$ started from $x_0$ with marked point $z_0$. Let $\til g_t$, $0\le t<T$, be the covering annulus Loewner maps of modulus $p$ driven by $\xi$. Then for every $z\in\St_p$,
$$M_t(z):=\frac{\Gamma_j(p-t,\til g_t(z)-\xi(t))}{\Gamma_j(p-t,\til g_t(z_0)-\xi(t))} \cdot\frac{\til g_t'(z)^{1/2}}{\til g_t'(z_0)^{1/2}} $$
is a local martingale for $0\le t<T$. \label{3-martingale}
\end{Proposition}

For $j=1$, $\Gamma_1(t,\cdot)$ takes nonzero pure imaginary values on $\R_t$, the related function $\Gamma_4$ agrees with the solution given by Section \ref{Feynman} for $\kappa=3$ and $\sigma=\frac 4\kappa-1=\frac 13$ up to a pure imaginary multiplicative constant,
and $\Lambda_4:=3\frac{\Gamma_4'}{\Gamma_4}$ is a crossing annulus drift function that solves (\ref{dot-Lambda=}) for $\kappa=3$.
The annulus SLE$(\kappa;\Lambda_4)$ process also preserves a local martingale. In fact, Proposition \ref{3-martingale} holds with $z_0\in\R_p$, $\Lambda_j$ replaced by $\Lambda_4$, and $\Gamma_j$ replaced by $\Gamma_1$.

\subsection{$\kappa=0$}
Let $\kappa =0$. Let $L_t$ be as in Section \ref{Section-kappa2}. Let $\HA_2(t,z)=\HA(t,z/2)$. From (\ref{pat-HA}) we have
\BGE \pa_t \HA_2=4\HA_2''+2\HA_2'\HA_2.\label{HA-2-PDE}\EDE
Let $\Lambda_1=\HA-2\HA_2$.
Then for each $t>0$, $\Lambda_1(t,\cdot)$ is an odd analytic function on $\C\sem L_t$, and each
$z\in L_t$ is a simple pole of $\Lambda_1$. From $\HA(t,z+2\pi)=\HA(t,z)$ and $\HA(t,z+i2t)=\HA(t,z)-2i$ we see that
both $4\pi$ and $i4t$ are periods of $\Lambda_1(t,\cdot)$. Fix $t>0$, and define
$$J(z)=\frac{\Lambda_1(t,z)^2}2-2\Lambda_1'(t,z)+3\HA'(t,z).$$
Then $J$ is an even analytic function on $\C\sem L_t$ and has periods $4\pi$ and $i4t$. Fix any $z_0=2n_0\pi+i2k_0t\in L_t$
for some $n_0,k_0\in\Z$. Then
$2 z_0$ is a period of $J$, so $J_{z_0}(z):=J(z-z_0)$ is an even function. Thus, $\Res_{z=z_0} J(z)=0$.
The degree of $z_0$ as a pole of $J$ is at most $2$. The principal part of $J$ at $z_0$ is $\frac {C(z_0)}{(z-z_0)^2}$ for
some $C(z_0)\in \C$. Note that $\Res_{z_0} \HA(t,z)=2$ and $\Res_{z_0}\Lambda_1(t,z)=-6$ or $2$. In either case, we compute $C(z_0)=0$.
Thus, every $z_0\in L_t$ is a removable pole of $J$, which, together with the periods $4\pi$ and $i4t$, implies that $J$ is a constant depending only on $t$. Differentiating $J$ w.r.t.\ $z$, we conclude that
\BGE 2\Lambda_1''=\Lambda_1'\Lambda_1+3\HA''.\label{G-PDE}\EDE
From $\Lambda_1=\HA-2\HA_2$ we have $2\HA_2=\HA-\Lambda_1$. So from (\ref{HA-2-PDE}) and (\ref{G-PDE}), we have
$$ \pa_t\HA-\pa_t \Lambda_1=2\pa_t\HA_2=8\HA_2''+4\HA_2'\HA_2=4\HA''-4\Lambda_1''+(\HA'-\Lambda_1')(\HA-\Lambda_1)$$
$$= 4\HA''-2(\Lambda_1'\Lambda_1+3\HA'')+(\HA'-\Lambda_1')(\HA-\Lambda_1) =-2\HA''-\Lambda_1'\Lambda_1+\HA'\HA-\Lambda_1'\HA-\HA'\Lambda_1.$$
From the above formula and (\ref{pat-HA}), we have
\BGE \pa_t \Lambda_1=3\HA'' +\Lambda_1'\Lambda_1+\HA'\Lambda_1+\Lambda_1'\HA. \label{dot-G}\EDE
Thus, $\Lambda_1$ solves (\ref{dot-Lambda=2}). Note that $\HA_I(t,z/2)$ also satisfies (\ref{HA-2-PDE}).
Let $\Lambda_2(t,z):=\HA(t,z)-\HA_I(t,\frac z2)$. Then $\Lambda_2(t,\cdot)$ is also an odd analytic function on $\C\sem L_t$ and
has periods $4\pi$ and $i4t$. The principal part of $\Lambda_2(t,\cdot)$ at every $z_0\in L_t$
 is also either  $\frac{-6}{z-z_0}$ or $\frac{2}{z-z_0}$.
Using a similar argument, we conclude that $\Lambda_2$ also solves   (\ref{dot-Lambda=2}).

\subsection{$\kappa=16/3$}
Let $\kappa=16/3$. Let $\Lambda_1$ and $\Lambda_2$ be as in the last subsection. Let $\Lambda_3=-\Lambda_1/3$. From (\ref{G-PDE}) we have
$$ 0=\frac 83\Lambda_3''+4\Lambda_3'\Lambda_3+\frac43\HA''. $$ 
From (\ref{dot-G}) we have
$$ \pa_t \Lambda_3=-\HA''-3\Lambda_3'\Lambda_3+\HA'\Lambda_3+\Lambda_3'\HA.$$ 
Summing up the above two equalities, we get
$$ \pa_t \Lambda_3=\frac 83 \Lambda_3''+\frac 13 \HA''+\HA'\Lambda_3+\Lambda_3'\HA+\Lambda_3'\Lambda_3.$$ 
Thus, $\Lambda_3$ solves (\ref{dot-Lambda=2}). Similarly, $\Lambda_4:=-\Lambda_2/3$ also
solves (\ref{dot-Lambda=2}). Here $\Lambda_3$ and $\Lambda_4$ have period $4\pi$ instead of $2\pi$. If we want a solution to (\ref{dot-Lambda=}) with period $2\pi$, we may first restrict $\Lambda_3$ or $\Lambda_4$ to the interval $(0,2\pi)$ or $(-2\pi,0)$, and then extend it to $\R\sem\{2n\pi:n\in\Z\}$ so that the function has period $2\pi$. 


\begin{thebibliography}{99}
\bibitem{elliptic} K.\ Chandrasekharan. {\it Elliptic functions}. Springer-Verlag Berlin Heidelberg, 1985.
\bibitem{Ahl} Lars V.\ Ahlfors. {\it Conformal invariants: topics
in geometric function theory}. McGraw-Hill Book Co., New York, 1973.
\bibitem{BC} R.\ F.\ Bass and Z.-Q.\ Chen. One-dimensional stochastic differential equations with singular and degenerate coeffients. {\it Sankhya: The India Journal of Statistics}, 67:19-45, 2005.
\bibitem{dim-SLE} V.\ Beffara. The dimension of the SLE curves. {\it Ann. Probab.}, 36(4):1421-1452, 2008.
\bibitem{Julien-Comm} Julien Dub\'edat. Commutation relations for SLE,
{\it Comm. Pure Appl.\ Math.}, 60(12):1792-1847, 2007.
\bibitem{Law} Gregory F.\ Lawler. {\it Intersection of random walks}.
Birkh\"auser, Boston, 1991.
\bibitem{LawSLE} Gregory F.\ Lawler. {\it Conformally Invariant Processes in the Plane}.
Am.\ Math.\ Soc., Providence, RI, 2005.
\bibitem{Law-restriction} Gregory F.\ Lawler. Defining SLE in multiply connected domains with the Brownian loop measure, arXiv:1108.4364v1.
\bibitem{LSW1} Gregory F.\ Lawler, Oded Schramm and Wendelin Werner.
Values of Brownian intersection exponents I: Half-plane exponents.
{\it Acta Math.}, 187(2):237-273, 2001.
\bibitem{LSW2} Gregory F.\ Lawler, Oded Schramm and Wendelin Werner.
Values of Brownian intersection exponents II: Plane exponents. {\it
Acta Math.}, 187:275-308, 2001.
\bibitem{LSW-8/3} Gregory F.\ Lawler, Oded Schramm and Wendelin Werner.
Conformal restriction: the chordal case, {\it J.\ Amer.\ Math.\
Soc.}, 16(4): 917-955, 2003.
\bibitem{LSW-2} Gregory F.\ Lawler, Oded Schramm and Wendelin Werner.
Conformal invariance of planar loop-erased random walks and uniform
spanning trees. {\it Ann.\ Probab.}, 32(1B):939-995, 2004.
\bibitem{RY} Daniel Revuz and Marc Yor. {\it Continuous Martingales
and Brownian Motion}. Springer, Berlin, 1991.
\bibitem{MS1} Jason Miller and Scott Sheffield. Imaginary Geometry I: interacting SLE paths, arXiv:1201.1496.
\bibitem{MS4} Jason Miller and Scott Sheffield. Imaginary Geometry IV: interior rays, whole-plane reversibility, and space-filling trees, arXiv:1302.4738.
\bibitem{RS-basic} Steffen Rohde and Oded Schramm. Basic properties of
SLE. {\it Ann.\ Math.}, 161(2):883-924, 2005.
\bibitem{S-SLE} Oded Schramm. Scaling limits of loop-erased random walks
and uniform spanning trees. {\it Israel J.\ Math.}, 118:221-288, 2000.
\bibitem{SS-4} Oded Schramm and Scott Sheffield. Harmonic explorer and its convergence to SLE$_4$.
{\it Ann.\ Probab.}, 33:2127-2148, 2005.
\bibitem{SS} Oded Schramm and Scott Sheffield.
Contour lines of the two-dimensional discrete Gaussian free field. {\it Acta Math.}, 202(1):21-137, 2009.
\bibitem{SW} Oded Schramm and David B.\ Wilson. SLE coordinate changes. {\it New York Journal of Mathematics}, 11:659--669, 2005.
\bibitem{SS-6} Stanislav Smirnov. Critical percolation in the plane:
conformal invariance, Cardy's formula, scaling limits. {\it C.\ R.\
Acad.\ Sci.\ Paris S\'er.\ I Math.}, 333(3):239-244, 2001.
\bibitem{SS-3} Stanislav Smirnov. Towards conformal invariance of 2D lattice
models. {\it International Congress of Mathematicians.},
vol. II, 1421-1451, Eur.\ Math.\ Soc., Zurich (2006)
\bibitem{Zhan} Dapeng Zhan. Stochastic Loewner evolution in doubly connected domains.
{\it Probab. Theory Rel.}, 129(3):340-380, 2004.
\bibitem{ann-prop} Dapeng Zhan. Some properties of annulus SLE.
{\it Electron.\  J.\  Probab.}, 11, Paper 41:1069-1093, 2006.
\bibitem{LERW} Dapeng Zhan. The Scaling Limits of Planar LERW in Finitely Connected Domains. {\it Ann.\ Probab.}, 36(2):467-529, 2008.
\bibitem{reversibility} Dapeng Zhan. Reversibility of chordal SLE. {\it Ann.\ Probab.}, 36(4):1472-1494, 2008.
\bibitem{duality} Dapeng Zhan. Duality of chordal SLE.  {\it Invent.\ Math.}, 174(2):309-353, 2008.
\bibitem{duality2} Dapeng Zhan. Duality of chordal SLE, II. {\it Ann.\ I.\ H.\ Poincare-Pr.}, 46(3):740-759, 2010.
\bibitem{kappa-rho} Dapeng Zhan. Reversibility of some chordal SLE$(\kappa;\rho)$ processes. {\it Journals \ Stat.\ Phys.}, 139(6):1013-1032, 2010.
\bibitem{int-LERW} Dapeng Zhan. Continuous LERW started from interior points. {\it Stoch.\ Proc.\ Appl.}, 120:1267-1316, 2010.
\bibitem{LEBM} Dapeng Zhan. Loop-Erasure of Planar Brownian Motion. {\it Commun.\ Math.\ Phys.}, 303(3):709-720, 2011.
\bibitem{ann-restriction}  Dapeng Zhan. Restriction Properties of Annulus SLE. {\it J.\ Stat.\ Phys.}, 146(5):1026-1058, 2012.

\end{thebibliography}
\end{document}